\documentclass[11pt]{amsart}
\usepackage[margin=1in]{geometry}

\let\originalleft\left
\let\originalright\right
\renewcommand{\left}{\mathopen{}\mathclose\bgroup\originalleft}
\renewcommand{\right}{\aftergroup\egroup\originalright}

\date{\today}
\usepackage{todonotes}

\usepackage{mathtools}
\mathtoolsset{showonlyrefs}

\usepackage{amssymb}
\usepackage{amsmath}
\allowdisplaybreaks
\usepackage{mathtools}
\usepackage{xfrac}
\usepackage{amsthm}
\usepackage{bookmark}
\usepackage{hyperref}
\usepackage{xcolor}
\usepackage{todonotes}
\usepackage{enumitem}

\usepackage{tikz-cd}

\usepackage{changes}

\usepackage{ulem}

\let \Im \undefined 
\DeclareMathOperator{\Im}{Im}

\let \Re \undefined 
\DeclareMathOperator{\Re}{Re}

\numberwithin{equation}{section}
\newtheorem{theorem} {Theorem}
\newtheorem{lemma}[equation]{Lemma}
\newtheorem{corollary}[equation]{Corollary}
\newtheorem{proposition}[equation]{Proposition}
\newtheorem{definition}[equation]{Definition}
\newtheorem{remark}[equation]{Remark}

\newtheorem{theoremofothers}[equation]{Theorem}

\let \Re\undefined
\newcommand{\R}{\mathbb{R}}
\newcommand{\oM}{ M}
\DeclareMathOperator{\Re}{Re}
\DeclareMathOperator{\Diff}{Diff}
\DeclareMathOperator{\End}{End}

\newcommand{\Aphg}{\mathcal{A}_{\mathrm{phg}}}
\newcommand{\bSM}{{}^{b\!}{S^*{M}}}
\newcommand{\bbSM}{{}^b{\overline{S^*{M}}}}
\newcommand{\OSM}{{}^0\!SM}
\newcommand{\OTM}{{}^0TM}
\newcommand{\OTsM}{{}^0T^*M}
\newcommand{\bTM}{{}^bTM}
\newcommand{\bTsM}{{}^bT^*M}
\newcommand{\OSD}{{}^0\!{S\mathbb{D}}}
\renewcommand{\d}{\mathrm{d}}
\newcommand \Rm { \mathbb{R}}

\newcommand{\ocup}{\overline{\cup}}
\DeclareMathOperator{\grad}{grad}
\DeclareMathOperator{\supp}{supp}
\DeclareMathOperator{\sgn}{sgn}
\DeclareMathOperator{\id}{id}

\definecolor{HaimOrange}{RGB}{191,78,20}

\author[N. Eptaminitakis] {Nikolas Eptaminitakis}
\address{ Institut für Differentialgeometrie, 
Leibniz Universit\"at Hannover,
	Welfengarten 1,
	 30167 Hannover, 	
	Germany}

\author[H. Grebnev]{Haim Grebnev}
\address{Department of Mathematics, Purdue University, West Lafayette, IN 47907}

\title[The Non-Abelian X-Ray Transform on Simple AH Surfaces]{The Non-Abelian X-Ray Transform on Simple Asymptotically Hyperbolic Surfaces}

\begin{document}

\begin{abstract}
We prove a gauge equivalence result for the non-Abelian X-ray transform with connection and Higgs field on simple 2-dimensional asymptotically hyperbolic spaces.
Our approach is based on the construction of holomorphic integrating factors. 
To construct these, we analyze the  geodesic X-ray transform on functions on simple asymptotically hyperbolic manifolds of general dimension and prove polyhomogeneous mapping properties and a surjectivity result for the adjoint.  
\end{abstract}

\maketitle 

\tableofcontents 

\section{Introduction and Main Results} \label{IntroSec}

In this paper we consider the problem of recovering up to gauge a connection and Higgs field  on a vector bundle over a non-compact surface with prescribed behavior at infinity, called asymptotically hyperbolic. 
Our approach is based on energy identities and the construction of holomorphic integrating factors in this setting.
It necessitates studying in detail suitable compactifications of the unit (co)sphere bundle and using them to analyze the geodesic X-ray transform on functions, in particular the surjectivity of its adjoint (without dimensional restriction).
Before we state the problem and main results, we introduce the geometric setting.

\subsection{Setup}\label{sec:Setup}

Let $M^{n+1}$ be an oriented compact manifold with boundary.
A Riemannian metric on the interior $M^\circ$ will be called \textbf{asymptotically hyperbolic (AH)}, if for some (and hence any) boundary defining function (for short, ``bdf'')\footnote{Recall that this means $\rho>0$ on $M^\circ$, $\rho=0$ on $\partial M$ and $\d \rho \neq 0$ on $\partial M$.} $\rho\in C^\infty(M)$, the metric $\bar g\coloneqq \rho^{2}g$ extends from $M^{\circ}$ to a smooth Riemannian metric on $M$ such that $\|\d \rho\big|_{\partial M}\|^2_{\bar g}=1$. 
AH metrics are geodesically complete and their sectional curvature approaches $-\|\d \rho\big|_{\partial M}\|^2_{\bar g}=-1$ as $\rho\to 0$ (see \cite[Proposition 1.10]{Mazzeo1986}). 
Each AH metric $g$ induces a conformal class of metrics on $\partial M$, called the conformal infinity,  given by $[\bar g\big|_{T\partial M}]$ in the notation above.
In part due to this fact they have been of great interest in conformal geometry (\cite{Fefferman1985}, \cite{Fefferman2012}) and in theoretical physics in the context of the AdS-CFT correspondence (\cite{Maldacena1998}).

For each representative $h_0$ in the conformal infinity there exists a bdf $\rho_0$, called a \textbf{geodesic bdf}, such that $\rho_0^{2}g\big|_{T\partial M}=h_0$  and such {}that in a neighborhood of $\partial M$, one has  $\|\d \rho_0\|_{\rho_0^{2}g}^{2}=1$ ($\rho_0$ is unique in that neighborhood).
The flow of the $\rho_0^{2}g$-gradient of $\rho_{0}$ induces a diffeomorphism 
\begin{equation}\label{eq:diffeo_collar}
	\psi:[0,\epsilon)_{\rho_{0}}\times \partial M\to \mathcal{C}_\epsilon,
\end{equation} where $\mathcal{C}_\epsilon$ is a collar neighborhood  of $\partial M $, such that for $\rho_0>0$ over $C_\epsilon$
\begin{equation}\label{eq:metric_normal_form}
	\psi^{*}g=\frac{\d \rho_0^2+h_{\rho_0} }{\rho_0^2},
\end{equation}
where $h_{\rho_0}$ is a smooth 1-parameter family of metrics on $\partial M$, extending smoothly down to ${\rho_0}=0$.  
Throughout, the subscript $0$ on a bdf will indicate without further comment that the bdf is geodesic.

The geometric inverse problems we will be discussing in the next section is formulated with the help of transport equations on the unit sphere or cosphere bundle.
We let 
\begin{equation}
	SM^\circ \coloneqq \{(z,v)\in TM^\circ: \|v\|_g =1 \}, \quad S^{*}M^\circ \coloneqq \{(z,\zeta)\in T^*M^\circ: \|\zeta\|_g =1 \}.
\end{equation}
The musical isomorphism $\flat  : TM^\circ\rightarrow T^\ast M^\circ$ and its inverse $\sharp : T^\ast M^\circ\rightarrow TM^\circ$ induce diffeomorphisms between these spaces, which we will often not denote explicitly.
The closures of $SM^\circ$ and $S^*M^\circ$ in $TM$, resp. $T^*M$  are not smooth manifolds and the subsequent analysis will necessitate working on suitable compactifications of them.
We denote by $X$ the generator of the geodesic flow on either of them (it should be clear from context which one we mean); they are pushforwards of one another by means of the musical isomorphisms. 
As observed in \cite{Graham2019}, transport equations with respect to $X$ can be effectively analyzed in the b-cosphere bundle defined below, whereas fiberwise Fourier analysis is better behaved on the 0-sphere bundle introduced in \S\ref{ssec:unit_cosphere_cpt}, as demonstrated in \cite{Grebnev2026,Eptaminitakis2025}. 
In \S\ref{ssec:unit_cosphere_cpt}, we will show that a suitable compactification of the former is naturally diffeomorphic with
a
space
obtained by
blowing up a certain submanifold of the latter.

Recall that the b-tangent bundle $\bTM$ (see \cite{Melrose1993}) is a smooth rank $n+1$-bundle over $M$ whose smooth sections are spanned by $\mathcal{V}_b(M)$, the Lie algebra of vector fields which are tangent to $\partial M$. 
In local coordinates $(\rho,y^1,\dots ,y^n)$ near $\partial M$, where $y^j$ restrict to coordinates on $\partial M$ and $\rho $ is a bdf as before, they are spanned over $C^\infty(M)$ 
by $\{\rho\sfrac{\partial}{\partial\rho},\sfrac{\partial}{\partial{y^1}},\dots,\sfrac\partial{\partial{y^n}}\}$. 
The dual bundle is $\bTsM$, locally spanned by $\frac{\d \rho}{\rho}, \d y^1,\dots, \d y^n$.
Over $M^\circ,$ $TM^\circ $ is canonically isomorphic with $\bTM\big|_{M^\circ}$, and similarly for the duals. 
This identification will be used throughout without explicit comment.
We define the \textbf{$b$-cosphere bundle} $\bSM$ as the closure of $S^*M^{\circ}$, viewed as a subset of ${}^bT^*M$;
 it is a non-compact, smooth  submanifold with boundary of ${}^bT^*M$.
Using local coordinates $(\rho_0,y)$ with $\rho_0$ a geodesic bdf as  before, we can write it locally near $\partial M$ as
\begin{equation}\label{eq:bSM}
	\bSM=\big\{(z,\zeta)=(\rho_0,y,\bar{\xi}\frac{\d\rho_0}{\rho_0}+{\eta}\cdot \d y)\in {}^bT^*M: \bar{\xi}^2+h_{\rho_0} ^{ij}\rho_0^{2}\eta_i\eta_j=1\big\}.
\end{equation}
Here and throughout, we use the Einstein notation where Latin indices  related to the dimension of $M$ always run from 1 to $n$. 
If $\partial M$ is connected, the boundary $\partial \bSM$ of $\bSM$ has two connected components, which we write as
\begin{equation}
	\partial_\pm \bSM=\big\{(z,\zeta)=(\rho_0,y,\bar{\xi}\frac{\d\rho_0}{\rho_0}+{\eta}\cdot \d y)\in {}^bT^*M:\ \rho_0 =0,\  \bar{\xi}=\mp 1 \big\}.
\end{equation}
As noted in \cite{Graham2019}, these sets are well defined independently of coordinates.
We mention that the authors \cite{Graham2019} use the notation $\overline{S^\ast M^\circ}$ for what we denote here by $\bSM$.

A natural setting in the study of geometric inverse problems is  \textbf{simple manifolds}, namely nontrapping compact manifolds with strictly convex boundary and no conjugate points. 
For instance, this is the setting in which the important boundary rigidity conjecture was initially posed (\cite{Michel1981}).
The AH analog is a \textbf{simple AH manifold}.  By definition, simple AH metrics are nontrapping (i.e., every parametrized geodesic $\gamma(t)$ leaves any given compact subset of $M^\circ $  for $|t| $ sufficiently large), and there are no boundary conjugate points, in the sense that there exists no Jacobi field $Y(t)$ along any unit speed geodesic $\gamma(t)$ such that $\|Y(t)\|_g\to 0$ as $t\to \pm \infty$.
This condition also excludes interior conjugate points, see \cite[Section 5]{Graham2019}, where this definition was formulated.
As already observed in \cite[p. 2871]{Graham2019}, simple AH manifolds are simply connected. Moreover, if $\rho $ is any bdf for $\partial M$, the set $\{\rho \geq \epsilon\}$ is a simple manifold for $\epsilon >0$ sufficiently small.

\bigskip

\subsection {The non-Abelian X-ray Transform}\label{sec:NonAbelTrasfDef} 

Since AH manifolds are complete, integral curves of $X$ are defined for all time.
Per Lemma 2.1 in \cite{Graham2019} and its proof, $X$ extends smoothly to $\bSM$ and is of the form $X=\rho\overline{X}$ where $\overline{X}$ is a smooth vector field over $\bSM$ that is transversal to $\partial \bSM$ pointing inwards and outwards on $\partial _-\bSM$ and $\partial _+\bSM$ respectively.
Since $X$ and $\rho\overline{X}$ are positive scalar multiples of each other over $S^\ast M^\circ$, their integral curves match up to reparameterization there.
It follows easily that every $ \left(z,\zeta \right)\in\partial _\pm\bSM$ is the limit of an integral curve $\overline{\sigma}\left(\tau\right) : \left[-a,a\right]\rightarrow\bSM$ of $\overline{X}$ as $\tau\rightarrow\pm a$ and hence is also the limit of an integral curve $\sigma\left(t\right) :  \left(-\infty,\infty \right)\rightarrow S^\ast M^\circ$ of $X$ as $t\rightarrow\pm\infty$, where $\sigma$ is a reparameterization of $\overline{\sigma}$.
We will often let $\sigma$ also denote $\sharp\circ\sigma$ which is an integral curve of $X$ over $SM^\circ$ without explicit mention.

Suppose $(M^{\circ},g)$ is a nontrapping AH manifold.
In this paper we consider an inverse problem associated with the non-Abelian X-ray transform which we now formulate.
Consider the Hermitian vector bundle $M\times\mathbb{C}^d\rightarrow M$ with respect to the inner product $\langle e_1,e_2\rangle _{\mathbb{C}^d}=\sum_{k=1}^{d}{ \left(e_1 \right)^k \left(\overline{e_2} \right)^k}$ on the fibers.
We denote the set of $C^\infty$ sections of this bundle as $C^\infty\left(M;\mathbb{C}^d\right)$.
Let $\Phi\in \rho C^\infty\left(M;\mathbb{C}^{d\times d}\right)$ be a $C^\infty$ complex $d\times d$ matrix field over $M$ and let $\nabla$ denote a $C^\infty$ connection on $M\times\mathbb{C}^d$.
Let $A : TM\rightarrow\mathbb{C}^{d\times d}$ be the \textbf{connection matrix} of $\nabla$:
\begin{equation}\label{eq:ConnExplForm}
\nabla_vu=v \left(u \right)+A \left(v \right)u.
\end{equation}
Throughout this paper we will also use the \textbf{pullback} of $\nabla$ to $SM^{\circ }$ defined as follows.
Let $\pi:TM\rightarrow M$ denote the natural projection. If $u$ is instead of the form $u:SM^{\circ }\rightarrow {\mathbb{C}}^{d}$, then
$$\nabla_{\omega }u\coloneqq\omega \left(u\right)+A\left(d\pi \left(\omega \right)\right)u ,\quad \omega \in TSM^{\circ}.$$
In particular, we will quite often use $\nabla_{X}u(x,v)=Xu+A\left(x,v\right)u$.

As in \cite[Section 2.5]{Grebnev2026}, we define the operator
$$T^{ \left(\Phi,\nabla \right)} : \partial _-\bSM\times\mathbb{C}^d\rightarrow\partial _+\bSM\times\mathbb{C}^d$$
as follows.
Fix a point $ \left(z,\zeta \right)\in\partial _-\bSM$ and an $e\in\mathbb{C}^d$.
Let $\sigma :  \left(-\infty,\infty \right)\rightarrow S^\ast M^\circ$ be an integral curve of $X$ with $ \left(z,\zeta \right)=\lim_{t\rightarrow-\infty}{\sigma \left(t \right)}$.
Letting $\gamma :  \left(-\infty,\infty \right)\rightarrow M^\circ$ denote the $g$-geodesic that is the projection of $\sigma$ onto $M^\circ$, consider the solution $u :  \left(-\infty,\infty \right)\rightarrow\mathbb{C}^d$ to the following transport equation:
\begin{equation}\label{eq:TranspEqDefineT}
\frac{du}{dt}+ \left(A\circ\sigma \right)u+ \left(\Phi\circ\gamma \right)u=0\qquad\mathrm{with}\qquad \lim_{t\rightarrow-\infty}{u \left(t \right)}=e,
\end{equation}
and define
\begin{equation}\label{ScattOp}
T^{ \left(\Phi,\nabla \right)} \left(z,\zeta,e \right)= \big(\lim_{t\rightarrow\infty}{\sigma \left(t \right)},\lim_{t\rightarrow\infty}{u \left(t \right)} \big),
\end{equation}
where the existence of $u$ and the limits in \eqref{ScattOp} follow from  \cite[Lemma 2.6]{Grebnev2026}.
The definition of \eqref{ScattOp} does not depend on the $\sigma$ chosen since if $\widetilde{\sigma}$ was another integral curve of $X$ with the same limit in $\partial _-S^\ast M$ as $t\rightarrow-\infty$, then $\widetilde{\sigma} \left(t \right)=\sigma \left(t+T \right)$ for some fixed $T\in\mathbb{R}$.
The map $T^{(\Phi,\nabla)}$ is often called the \textbf{scattering data} of the pair $(\Phi,\nabla)$.
The \textbf{non-Abelian X-ray transform}  is  the map
\begin{equation}\label{NonAbelTransfDef}
\left(\Phi,\nabla \right)\mapsto T^{ \left(\Phi,\nabla \right)}.
\end{equation}
The inverse problem is whether $ \left(\Phi,\nabla \right)$ can be recovered from knowing $T^{ \left(\Phi,\nabla \right)}$ (i.e. is the map \eqref{NonAbelTransfDef} injective)?
As is known for the case of nontrapping compact manifolds with boundary, the answer here is also no. 
Indeed, if  $\nabla$ is a $C^\infty$ connection on $M\times\mathbb{C}^d$ and $\Phi\in\rho C^\infty \left(M;\mathbb{C}^{d\times d} \right)$, then for any
$Q\in C^\infty \left(M;\mathbb{C}^{d\times d} \right)$ with $Q=\mathrm{id}$ on $\partial M$ one has $T^{ \left(\Phi,\nabla \right)}=T^{ (\widetilde{\Phi},\widetilde{\nabla} )}$, where
\begin{equation}\label{eq:GaugeEq}
\widetilde{\nabla}\coloneqq Q^{-1}\nabla Q\qquad\text{and}\qquad\widetilde{\Phi}\coloneqq Q^{-1}\Phi Q.
\end{equation}
This can be checked with a straightforward computation, e.g., see \cite[p. 10]{Grebnev2026}.
We say that the pairs $ \left(\Phi,\nabla \right)$ and $ \big(\widetilde{\Phi},\widetilde{\nabla} \big)$ are \textbf{gauge equivalent} if there exists a $Q$ as above for which \eqref{eq:GaugeEq} holds.
The nontrivial question we address in this paper is whether this gauge equivalence is the only obstruction to the non-invertibility of the non-Abelian X-ray transform \eqref{NonAbelTransfDef}.
Under suitable growth and regularity assumptions, we settle this question in the affirmative in Theorem~\ref{thrm:GaugeExist} below.
We mention that the scalar case $d=1$ with $A=0$, resp. $\Phi=0$  corresponds to the usual (Abelian) X-ray transform on functions, resp. 1-forms, an object which has been studied extensively in the literature.

\subsection {Main Results}

Recall that a connection on $M\times \mathbb C^d$ is called \textbf{unitary} if it satisfies the product rule: $v\langle u,w\rangle _{\mathbb{C}^d}=\langle \nabla_vu,w\rangle _{\mathbb{C}^d}+\langle u,\nabla_vw\rangle _{\mathbb{C}^d}$ for $v\in TM$ and differentiable sections $u,w : M\rightarrow\mathbb{C}^d$.
This is equivalent to the connection matrix of $\nabla$ being skew-Hermitian (c.f. \eqref{eq:ConnExplForm}).
Our main result regarding the inverse problem for the non-Abelian X-ray transform is the following:
\begin{theorem}\label{thrm:GaugeExist}
Let $ \left(M^\circ,g \right)$  be a simple AH surface and suppose that $\nabla,\widetilde{\nabla}$ are $C^\infty$ unitary connections on $M\times\mathbb{C}^d$ and that $\Phi,\widetilde{\Phi}\in\rho C^\infty \left(M;\mathbb{C}^{d\times d} \right)$ where $\Phi,\widetilde{\Phi}$ are skew-Hermitian everywhere.
If $T^{ \left(\Phi,\nabla \right)}=T^{ (\widetilde{\Phi},\widetilde{\nabla} )}$, then there exists a $Q\in C^\infty \left(M;\mathbb{C}^{d\times d} \right)$ satisfying $Q=\mathrm{id}$ on $\partial M$ for which \eqref{eq:GaugeEq} holds.
\end{theorem}

\noindent We believe that the regularity assumed for the connections and Higgs fields in Theorem~\ref{thrm:GaugeExist} is not sharp, however some regularity and growth assumptions are needed to make sense of their non-Abelian X-ray transform, see \cite[Lemma 2.6]{Grebnev2026} and the discussion following it.
Further, it is reasonable to expect that some decay assumptions are essential and not only artifacts of the proof.
On $\R^{2}$ there are (polynomially decaying) connections whose scattering data agree with those of the flat connection but are not gauge equivalent to the latter (known as ``ghosts''), see  \cite{Novikov2002b}.
Our connection and Higgs field are polynomially weighted from the point of view of the compactified manifold, but they are exponentially weighted from the point of view of the AH geometry, that is, in terms of the geodesic distance from an interior point.

As a first step towards proving Theorem~\ref{thrm:GaugeExist} we prove boundary determination in Theorem~\ref{thm:ArbDecay} below, that is, we show that one can replace the pair $(\widetilde \nabla, \widetilde \Phi)$ with a gauge equivalent one for which \eqref{eq:GaugeEq} holds up to infinite order at $\partial M$.
A pseudolinearization argument (see \cite[\S 13.2]{Paternain2023} and \eqref{eq:TranspEqForW} below) then reduces the (nonlinear) problem of determining the injectivity of the non-Abelian X-ray transform  up to gauge to a \emph{linear} inverse problem for an attenuated X-ray transform, which is of independent interest.
To formulate our result regarding the latter, we  use vertical Fourier analysis on the fibers of $SM^{\circ}$.
As is well known, on a smooth oriented Riemannian surface $M^{\circ}$ one can decompose a function $u\in C^{\infty}(SM^{\circ}; \mathbb C^{d})$ in a Fourier series with respect to the fiber variable as 
\begin{equation}\label{eq:FourierDecomp}
	u=\sum_{k=-\infty}^{\infty} u_k,\quad u_{k}\in \Omega_k\coloneqq \{ u\in C^{\infty}(SM^{\circ}; \mathbb C^{d})\colon -i V u = ku \},
\end{equation}%
where $V$ is the generator of the counterclockwise rotation on the fibers, acting on $u$ component-wise -- see also \eqref{eq:Xconf}.
\begin{theorem}\label{MainTranspEqFourierThrm} Let $ \left(M^\circ,g \right)$ be a simple AH surface and suppose that $\nabla$ is a unitary $C^\infty$ connection in $M\times\mathbb{C}^d$, $\Phi\in\rho C^\infty \left(M;\mathbb{C}^{d\times d} \right)$ is a skew-Hermitian matrix field, and that $\phi$ and $B$ are such that $\phi\in\rho^NC^\infty(M;\mathbb{C}^d)$ and $B\in\rho^{N+1}C^\infty(M;\mathrm{Hom}(TM;\mathbb{C}^d))$ for all integers $N\geq0$.
Then any solution $u\in C^\infty \left(\bSM;\mathbb{C}^d \right)$ to
\begin{equation}\label{eq:TranpEqFourier}
\nabla_Xu+\Phi u=\phi+B
\end{equation}
satisfying $\left.u\right|_{\partial \bSM}=0$ only depends on position and not on velocity (i.e. $u\in\Omega_0$).
\end{theorem}

As will be explained in the proof, the assumptions in the above theorem imply that the right-hand side of the transport equation \eqref{eq:TranpEqFourier} is in $\Omega_{-1}\oplus\Omega_0\oplus\Omega_1$.
The proof of Theorem~\ref{MainTranspEqFourierThrm} is based on an energy identity (eq. \eqref{eq:PestovProofV}).
It belongs to a family of energy identities called \textbf{Pestov Identites}, which have been very successful in addressing a variety of geometric inverse problems (see, e.g. \cite{Paternain2023} for an overview of applications).
Crucial for the applicability of the Pestov identity is a detailed understanding of the regularity of certain solutions to transport equations on $\bSM$. This is addressed in  Proposition~\ref{prop:TranspEqReg} below, which improves a similar result proved in \cite{Grebnev2026}. 
Moreover, in the Pestov identity  we need to control the sign of various quantities. 
This can be done if assumptions on the curvature of $g$ or the connection $\nabla$ are available; to avoid making such assumptions, we make use of so-called (anti)holomorphic integrating factors.
A function $w\in C^{\infty}( SM^{\circ })$ is called \textbf{fiber-wise holomorphic} if $w\in\oplus_{m=0}^\infty\Omega_m$ and \textbf{fiber-wise anti-holomorphic} if $w\in\oplus_{m=-\infty}^{0}\Omega_m$.
For a differential operator $Lu=\nabla_Xu+\psi u$ on $SM^\circ$ where $\psi:SM^\circ\rightarrow\mathbb{C}$, a \textbf{(anti-)holomorphic integrating factor} is a fiberwise (anti-)holomorphic function $\mu$ such that $\mu Lu=\nabla_X(\mu u)$.
They were first constructed on simple surfaces  in \cite{Salo2011} and one of the main novelties of the present work is carrying out  out their construction in the AH setting.
The statement below  uses polyhomogeneous conormal functions and edge differential operators (see Appendix \ref{Sec:ConormPolyhomFuncSp}) on the 0-sphere bundle $\OSM$, a suitable compactification of $SM^\circ$ (see \S\ref{ssec:unit_cosphere_cpt}), as well as the Guillemin-Kazhdan operators $\eta_\pm$, see \S\ref{sec:vertical_fourier}.

\begin{theorem}\label{HolomIntFactExistThrm}
Let $\left(M^{\circ },g \right)$ be a simple AH surface and let $f\in \Aphg^E(M)$ be real valued  with $E=\overline{\{(0,1)\}}$, that is, assume that $f$ has an expansion near $\partial M$ of the form
\begin{equation}
	f\sim a_{01}(y)\log \rho+a_{00}(y)+a_{11}(y)\rho \log(\rho)+a_{10}(y)\rho+\dots, \qquad a_{ij}\in C^\infty(\partial M),
\end{equation}
 and set $\varphi \coloneqq -\eta_+ f+\eta_- f$.
Then there exists a fiber-wise holomorphic $C^\infty$ function $w : SM^{\circ}\rightarrow\mathbb{C}$ satisfying $X \left(w \right)=\varphi$.
For any $P\in \Diff_e^{*}(\OSM)$ there exists  an $m\in \mathbb{R}$ (depending on the order of $P$) such that  $\rho^{m}Pw\in L^2 (SM^{\circ})$ and $\ \rho^{m}e^{w}\in L^\infty (SM^{\circ})$.
The same statement holds with ``holomorphic'' replaced by ``anti-holomorphic''.

\end{theorem}

\begin{remark}\upshape
	In Theorem~\ref{HolomIntFactExistThrm} one could consider  $\varphi$'s of more general form, however proving boundedness of $\rho ^{m}e^{w}$ would require additional arguments and the present form suffices for our  application.
\end{remark}

A crucial step in the construction of holomorphic integrating factors in Theorem~\ref{HolomIntFactExistThrm} is a surjectivity result for the backprojection operator $I_0^{*}$.
Roughly speaking, for a function $h\in C_c^{\infty }(\partial_-\bSM)$, $I_0^{*}h\in C^\infty(M^\circ)$ is defined by extending $h$ to be constant along the orbits of the geodesic flow and then taking the fiberwise average.
It is, in a suitable sense, the formal adjoint of the geodesic X-ray transform $I_0$ on functions (see \S\ref{sub:the_scalar_x_ray_transform_and_its_mapping_properties}), hence the notation.
Many important developments in geometric inverse problems on simple surfaces in the last two decades, such as those regarding boundary rigidity,  attenuated X-ray transform, and tensor tomography, have relied crucially on the surjectivity of $I_0^*$, proved in that setting in \cite{Pestov2005} using microlocal methods.
See \cite{Paternain2023} for an overview.
The following theorem is a surjectivity result for $I_0^*$ on AH manifolds of any dimension $\geq 2$.
The notations on conormal/polyhomogeneous functions and index sets are explained in Appendix~\ref{sec:polyhomongeneity}.

\begin{theorem}\label{thm:surj_I0}  Let $(M^{\circ},g)$ be a simple AH manifold of dimension $n+1$, $n\geq 1$.
	Suppose $f\in \Aphg^E(M)$, for some smooth index set $E=\{(s_m,p_m)\colon m=1,2,\dots\}$.
Then there exists a function $h\in C^\infty({\partial_-\bSM})$ such that $I_0^{*}h=f$ and it admits a partial polyhomogeneous expansion in the sense that there exist $a_{m,k}\in C^{\infty}\big(\partial\big(\overline{\partial_-\bSM}\big)\big)$ 
such that
\begin{equation}\label{eq:first_terms}
	h- \sum_{\Re (s_m)<n}\sum_{k=0}^{p_m}a_{m,k}(y,\omega)\delta^{s_m}\log^{k}(\delta)\in \mathcal{A}_-^{n}(\overline{\partial_-\bSM}) ,
\end{equation}
where $\delta $ denotes a bdf for $\partial(\overline{\partial_-\bSM})$ in the compactification $\overline{\partial_-\bSM}$ of ${\partial_-\bSM}$ discussed in \S\ref{ssec:unit_cosphere_cpt}.
\end{theorem}
Our proof of Theorem~\ref{thm:surj_I0} relies partly on a stability result for the \textbf{normal operator} $I_0^{*}I_0 $ proved in \cite{Eptaminitakis2022} using the 0-calculus of Mazzeo-Melrose (\cite{Mazzeo1987}). 
This result can be used to show surjectivity of $I_0^{*}$ onto certain function spaces (see Proposition~\ref{prop:solvability}), which however are not large enough for our application, nor for applications in other geometric inverse problems.
In the AH setting, one typically needs surjectivity on a space which contains unbounded functions that grow logarithmically at $\partial M$.
To achieve this, we use Melrose's pullback and pushforward theorem (\cite{Melrose1992}) to study in detail polyhomogeneous mapping properties of $I_0$ and $I_0^{*}$ (Proposition~\ref{prop:phg_mapping})  and in addition we identify the leading order coefficient of $I_0^*h $ when $h$ is polyhomogeneous (Lemma~\ref{lm:surjectivity_terms}).
Then, given a polyhomogeneous function $f\in \Aphg^{E}(M)$ as in the statement above, we  explicitly construct a partial expansion for a polyhomogeneous function $ h_1$ on $\overline{\partial _-\bSM}$ in terms of the expansion of $f$, such that $f-I_0^{*}h_1$ becomes regular enough that the aforementioned surjectivity result in Proposition~\ref{prop:solvability} applies to it.

It is plausible that  the function $h$ in the statement of Theorem~\ref{thm:surj_I0} can be shown to be fully polyhomogeneous, using a method analogous to the one  in \cite[Proposition 6.65]{Eptaminitakis2022} to improve the conclusion of Lemma~\ref{lm:phg_normal_op} to polyhomogeneity of the function $h$ there.

\subsection{Motivation and Historical Discussion}

The non-Abelian X-ray transform for connections and Higgs fields is of interest not only on its own, but also due its applications in both the pure and the applied mathematics literature.
One of the earlier contexts in which it was considered was gauge theory, in particular in the study of certain reductions of the self-dual Yang-Mills equations in 2+2 dimensions  (\cite{Manakov1981}, \cite{Ward1988}). In fact the aforementioned ghosts constructed by Novikov on $\R^{2}$ were inspired by the gauge theory literature.
The non-Abelian transform for Hermitian connections also becomes relevant in inverse problems for the wave equation associated with the Schrödinger equation with external Yang-Mills potential, where it can be recovered from the hyperbolic Dirichlet-to-Neumann map (\cite{Finch2001}).
Apart from the more theoretical applications, Polarimetric Neutron Tomography, a tomographic technique in which one measures the effect of a magnetic field on the spin of a polarized neutron beam, is based on the non-Abelian X-ray transform for Higgs fields, see \cite{Desai2020}.

We now give an overview of the literature, big part of which has been concerned with the study of the non-Abelian X-ray transform on $\R^{n}$ and on compact manifolds with boundary. We focus on the non-Abelian case $d\geq 2$.
On $\R^{2}$, the paper \cite{Vertgeim1992} studied the matrix attenuated transform and a nonlinear problem of recovering a Higgs field from parallel transport data along families of curves. 
Finch and Uhlmann (\cite{Finch2001}) proved a gauge equivalence result for the non-Abelian X-ray transform for unitary connections in $\R^2$ under assumptions on the curvature of the connection, by adapting the energy method approach of \cite{Mukhometov1981}-\cite{Vertgeim1992}. 
Novikov (\cite{Novikov2002b}),  proved injectivity of the non-Abelian X-ray transform for connections modulo gauge on $\R^{n}$, $n\geq 3$ and locally near the flat connection on $\R^{2}$, with reconstructions. He also provided examples of ghosts on $\R^{2}$, as mentioned earlier.
Slightly later, \cite{Eskin2004} studied the non-Abelian X-ray transform for compactly supported but not necessarily small connections on $\R^{2}$. 

In non-constant curvature, Sharafutdinov (\cite{Sharafutdinov2000})  proved that on a simple manifold, two Hermitian connections with the same scattering data which are sufficiently close to a given one must be gauge equivalent, under \textit{a priori} curvature assumptions.
The authors in \cite{Paternain2012a} studied the attenuated X-ray transform in the presence of a unitary connection and Higgs field, as well as the nonlinear problem of injectivity modulo gauge for the non-Abelian X-ray transform, on simple surfaces. 
\cite{Guillarmou2016} studied the attenuated X-ray transform for unitary connections and Higgs fields on compact manifolds with boundary and on closed manifolds of general dimension under the assumption of negative curvature and proved a gauge equivalence for the non-Abelian X-ray transform in the former case.
More recently, Paternain and Salo addressed the inverse problem for the non-Abelian X-ray transform for arbitrary connections on simple surfaces in \cite{Paternain2024}.
In \cite{Paternain2019b}, the non-Abelian X-ray transform for arbitrary connections was studied on compact Riemannian manifolds with boundary of dimension $\geq3$ admitting a strictly convex function.
A range characterization for the non-Abelian X-ray transform of unitary connections and Higgs fields on simple surfaces was proved in \cite{Bohr2023}.
Stability estimates for the Higgs field case on simple surfaces as well as statistical aspects are discussed in \cite{Monard2021},
 later generalized in \cite{BohrStabGe3} to higher dimensions.

In the recent years, there has been considerable interest in geometric inverse problems in the setting of AH manifolds, revolving around the X-ray transform, tensor tomography, and boundary rigidity (\cite{Graham2019}, \cite{Eptaminitakis2021}, \cite{Eptaminitakis2022}, \cite{Eptaminitakis2025}, \cite{Lefeuvre2020}, \cite{humbert2026}).  This can be viewed in a broader context, where inverse problems are studied on non-compact manifolds with curvature assumptions or prescribed behavior at infinity, e.g. asymptotically conic/Euclidean (\cite{Guillarmou2022}, \cite{Vasy2024}, \cite{Jia2022}), Cartan-Hadamard (\cite{Lehtonen2017}), or gas giant \cite{Hoop2024}.
The work closest to ours is \cite{Grebnev2026}, where the second-named author proved injectivity up to gauge for the non-Abelian X-ray transform on AH manifolds of dimension $\geq 2$ under curvature assumptions.
In the present work we drop them in dimension 2 and relax the assumptions on the vanishing of the connection and Higgs field.
We finally mention that inverse problems for versions of the non-Abelian X-ray transform can  also be posed on closed manifolds, see  \cite{Paternain2013b}, \cite{Cekic2025}.

A recent survey on the history of the non-Abelian X-ray transform and its applications can be found in \cite{Novikov2019}.

\subsection{Outline} The rest of the paper is organized as follows.
In most of \S\ref{sec:the_unit_sphere_and_cosphere_bundle_of_an_ah_manifold} we work on AH manifolds of general dimension. We compactify $\bSM$ and relate this compactification with a natural compactification of $SM^{\circ}$ in \S\ref{ssec:unit_cosphere_cpt}. Subsequently we provide some background on measures in this space (\S\ref{sec:sasaki}) and use it to prove conormal and polyhomogeneous mapping properties for the X-ray transform on functions and its adjoint (\S\ref{sub:the_scalar_x_ray_transform_and_its_mapping_properties}). We close the section with some background on vertical Fourier analysis in dimension 2 (\S\ref{sec:vertical_fourier}).
In \S\ref{sec:the_non_abelian_x_ray_transform} we prove Theorems \ref{thrm:GaugeExist}
 and \ref{MainTranspEqFourierThrm} assuming Theorem~\ref{HolomIntFactExistThrm}.
 In \S\ref{sec:RegOfSol} we prove a result for the regularity of solutions to transport equations on $\bSM$.
  In \S\ref{sec:FourAnalyTranspEq} we prove Theorem~\ref{MainTranspEqFourierThrm}. We then turn to the non-Abelian X-ray transform and prove a boundary determination result in \S\ref{sec:ArbDecay}  (Theorem~\ref{thm:ArbDecay}) and then Theorem~\ref{thrm:GaugeExist} is proved in \S\ref{sec:GaugeEquiv}. 
 In \S\ref{sec:surj_i0*} we prove Theorem~\ref{thm:surj_I0} with the help of two intermediate results, proved in Sections \ref{sub:leading_order} and \ref{sub:surjectivity_on_sobolev_functions} respectively.
 Theorem~\ref{HolomIntFactExistThrm} is proved in \S\ref{sec:holo_int}.
 Finally, Appendix \ref{Sec:ConormPolyhomFuncSp} contains some background on geometric microlocal analysis. 

\subsection{Acknowledgments} 
The authors thank Robin Graham, François Monard, Álvaro Sánchez Hernandez, and Gunther Uhlmann for helpful discussions.
Part of this research was carried out during academic visits of N.E. at Purdue University in 2025 and 2026 supported by the NSF grants DMS-2154489 and DMS-2452757 of Plamen Stefanov, who he would like to thank for the invitation and the hospitality. N.E. also thanks the Graduate Academy of Leibniz University Hannover for travel support.

\subsection{AI Statement} No AI tools were used during the development or writing of this paper. 

\section{The Unit (Co)Sphere Bundle and the Geodesic X-ray Transform}
\label{sec:the_unit_sphere_and_cosphere_bundle_of_an_ah_manifold}

\subsection{The Unit (Co)sphere Bundle and its Compactifications}\label{ssec:unit_cosphere_cpt}

Let $(M^{\circ},g)$ be AH of dimension $n+1$, $n\geq 1$. Throughout, it is assumed for convenience that $\partial M$ is connected.
Apart from $\bSM$, we will also need a compactification of $SM^\circ$, which is obtained by viewing it as a subset of the 0-tangent bundle.
The discussion regarding the 0-tangent and cotangent bundles is completely parallel to the one regarding the b-tangent and cotangent bundles. 
The 0-tangent bundle $\OTM$ (see \cite{Mazzeo1987}) is a smooth bundle over $M$ whose smooth sections are spanned by $\mathcal{V}_0(M)$, the Lie algebra of vector fields that vanish on $\partial M$. 
In local coordinates $(\rho,y^1,\dots ,y^n)$  as before, they are spanned over $C^\infty(M)$ 
by $\{\rho\sfrac{\partial}{\partial\rho},\rho \sfrac{\partial}{\partial{y^1}},\dots,\rho \sfrac\partial{\partial{y^n}}\}$. 
The dual bundle is $\OTsM$, locally spanned by $\frac{\d \rho}{\rho}, \frac{\d y^1}{\rho},\dots, \frac{\d y^n}{\rho }$.
The 0-unit sphere bundle $\OSM$ is the closure of $S M^{\circ}$ viewed as a subset of ${}^0TM$.
It is a compact manifold with boundary.
In terms of coordinates $(\rho_0,y)$,
\begin{equation}\label{eq:0SM}
	\OSM=\big\{(z,v)=(\rho,y,a{\rho_0}\partial_{\rho_0}+b\cdot{\rho_0}\partial_ y)\in {}^0TM: a^2+(h_{\rho_0})_{ij} b^{i}b^{j}=1\big\}.
\end{equation}
It is easy to check that the inclusion map $\iota:{}^0\!SM|_{M^\circ}\rightarrow TM^\circ$ extends smoothly to $\OSM$ up to the boundary. 
From another point of view, it is not hard to check that the AH metric $g$ induces a \emph{smooth} fiber metric on the vector bundle ${}^0TM\to M$, up to $\partial M$, and $\OSM$ is exactly its unit sphere bundle.
The space $\OSM$ has a natural edge structure, see Remark~\ref{rmk:edge} and \S\ref{sub:edge_structures}.

We will explain in detail how the two spaces $\bSM$ and $\OSM$ are related.
Since over $M^\circ$, $\OTM\big|_{M^{\circ}}$,  resp. $\bTsM\big|_{M^{\circ}}$, is canonically isomorphic with $TM^{\circ}$, resp. $T^*M^{\circ}$,
 by raising an index, the metric $g$ induces a diffeomorphism
\begin{equation}\label{eq:diff0b}
 	\Psi_g=\sharp:(\bSM)^{\circ}\to (\OSM)^{\circ }.
 \end{equation}
In coordinates as in \eqref{eq:bSM} and \eqref{eq:0SM} it takes the form
\begin{equation}\label{eq:diff0b:coord}
	(\rho_0,y,\bar{\xi},\eta)\mapsto (\rho_0,y, \bar{\xi},\rho_0 h_{\rho_0}^{-1}\eta)=(\rho_0,y,a,b).
\end{equation}
As we explain below, there exists a suitable compactification of $\bSM$ and a space arising from blowing up a certain submanifold of $\OSM$  such that $\Psi_g$ extends to a diffeomorphism between them.

\begin{figure}[h]
\centering
	\includegraphics[scale=.9]{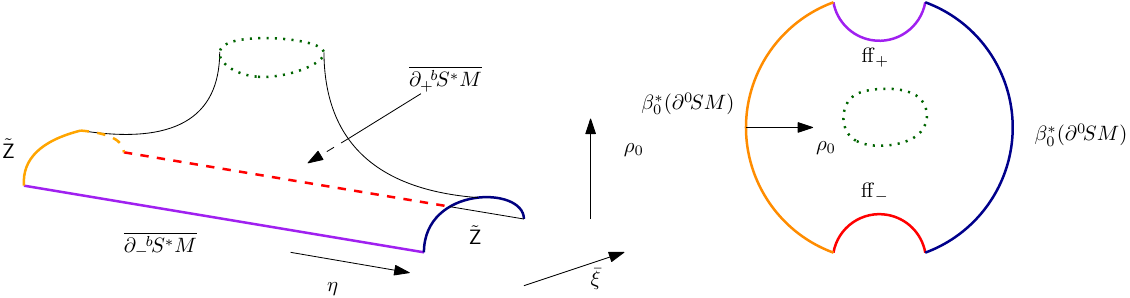}
	\caption{Case $n=1$. Left: $\bbSM$. Right: $[\OSM;\mathsf{Y}]$. $\partial M$-variable suppressed in both.
Colors indicate associated sets under the map $\overline{\Psi}_g$ in \eqref{eq:ext_diffeo}}\label{fig:left_bbsm_right_osm_y_}
\end{figure}

\begin{figure}[h!]
	\includegraphics[scale=1]{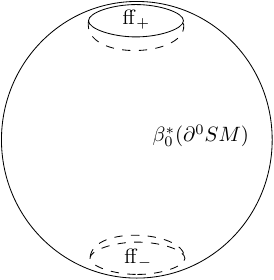}
	\caption{Case $n=2$. The boundary faces of $[\OSM;\mathsf{Y}]$, with the variables on $\partial M$ suppressed.}
	\label{fig:case_n_2_the_boundary_faces_of_osm_mathsf_y}
\end{figure}

Consider an arbitrary bdf $\rho$ for $\partial M$ (not necessarily geodesic).
Since $\bar{g}\coloneqq\rho^{2}g$ is a smooth metric on $M$,  it induces an outward pointing unit normal to $\partial M$; denote it by $\nu$. 
Extending it arbitrarily to  $M$ we obtain a  section of $TM$, which we still denote by $\nu$.
Then $\rho\nu\in \mathcal{V}_0(M)$ and in fact $\rho\nu\big|_{z}\in \OSM$ for every $z\in \partial M$. 
Moreover, one can check that if we start with a different bdf $\rho'=e^{h}\rho$, where $h\in C^\infty(\oM)$, inducing an outward pointing unit normal vector field $\nu'$, then $\rho'\nu'\big|_{\partial M}=\rho\nu\big|_{\partial M}$ as sections of ${}^0TM\big|_{\partial M}$.
Set
\begin{equation}
	\mathsf Y_{\pm }\coloneqq \{\pm \rho \nu \big|_{z}:z\in\partial M\}\subset \OSM.
\end{equation}
Near a point in $\mathsf Y_{\pm}$ we can parametrize $\OSM$ using the coordinates $(\rho_0,y,{b})$ in \eqref{eq:0SM} and in terms of them, $Y=\{\rho_0=0,{b}=0 \}$.
This shows that
\begin{equation}\label{eq:Y}
 	\mathsf{Y}\coloneqq \mathsf Y_{+}\cup \mathsf Y_{-}
 \end{equation}  is a codimension $n+1$ boundary $p$-submanifold of $\OSM$ with two connected components (see \S\ref{sub:manifolds_with_corners_p_submanifolds_blowups}).
We can therefore consider the blown-up space
\begin{equation}\label{eq:blowup}
	[\OSM;\mathsf Y],
\end{equation}
which is a manifold with corners up to  codimension $2$, with blow-down map $\beta_0$.
The front faces of this blow up we denote by $\mathrm{ff}_\pm = \beta_0^{*}\mathsf{Y}_\pm $. Here $\beta_0^{*}$ denotes lift of a submanifold, see \S\ref{sub:manifolds_with_corners_p_submanifolds_blowups}.
We also write $\mathrm{ff}=\mathrm{ff}_+\cup \mathrm{ff}_-$.
The involution $A_\sharp:SM^{\circ}\to SM^{\circ},$ $(z,v)\mapsto (z,-v)$ extends uniquely by continuity from the interior to a   diffeomorphism of $[\OSM;\mathsf Y]$ mapping $\mathrm{ff}_\pm\to \mathrm{ff}_\mp$. 
Since a collar neighborhood of  $\partial \OSM$ is diffeomorphic to $ \partial M\times [0,\epsilon) \times \mathbb{S}^{n} $ and $\mathsf{Y}_{\pm}\cong\partial M\times \{0\}\times \{P_\pm\}$, where $P_\pm$ denotes the north/south pole, each of the connected components of the front face is diffeomorphic to $\partial M\times \mathbb{S}^n_+$, where $\mathbb{S}^n_+=\{(v_{{0}},\dots,v_{{n}})\in\mathbb{S}^n:v_{{0}}\geq0\}$.
The remaining boundary hypersurface $\beta_0^{*}(\partial \OSM)$ of $[\OSM;\mathsf Y]$ is diffeomorphic to $\partial M\times  \mathsf{Z}, $ where $\mathsf{Z}\coloneqq[0,\pi]\times \mathbb{S}^{n-1}$ is a cylinder.
Note that in the 2-dimensional case, where $n=1$, $\mathsf{Z} $ has two connected components. See Figs.~\ref{fig:left_bbsm_right_osm_y_}-\ref{fig:case_n_2_the_boundary_faces_of_osm_mathsf_y}.

In the statement below, we write $\mathcal S^{*}(\partial M)\coloneqq (T^*\partial M\setminus 0)/\R^+ $ for the algebraic cosphere bundle, where $\R^+=(0,\infty)$ acts by scalar multiplication.

\begin{proposition}\label{lm:compactifications}

 The set
\begin{equation}
	\bbSM\coloneqq \bSM \sqcup \tilde{\mathsf{Z}},\quad \tilde{\mathsf{Z}}\coloneqq  [0,\pi]\times \mathcal S^{*}(\partial M),
\end{equation}
admits
a unique smooth structure  as a compact manifold with corners up to codimension 2 such that the  diffeomorphism $\Psi_g$ determined by $g$ in \eqref{eq:diff0b} extends to a diffeomorphism 
\begin{equation}\label{eq:ext_diffeo}
	\overline{\Psi}_g:\bbSM\to [\OSM;\mathsf Y].
\end{equation}
This diffeomorphism maps $\overline{\partial_\pm\bSM} $ diffeomorphically onto $\mathrm{ff}_{\pm }$.
The involution $A_\flat:S^*M^\circ\to S^*M^\circ$, $(z,\zeta)\mapsto (z,-\zeta)$ extends to a diffeomorphism of $\bbSM$ which satisfies $ A_{\flat}=\overline{\Psi}_g^{-1}\circ A_\sharp\circ \overline{\Psi}_g$, hence it interchanges $\overline{\partial_\pm \bSM}$.
\end{proposition}
\begin{proof}
Choose a conformal representative $h_0$ in the conformal infinity, fixing a geodesic bdf $\rho_0$ and a diffeomorphism $\psi$ as in \eqref{eq:diffeo_collar}.
If $\epsilon>0$ is sufficiently small, recall from \cite[Sec. 2.2]{Graham2019} that $\bSM$ near its boundary can be covered by the following types of neighborhoods:
\begin{equation}
\begin{aligned}
		U_1^{\pm }\coloneqq \big\{(\rho_0,y,\bar\xi,\eta)\in \bSM;\  \rho_0<\epsilon,\  |\eta|_{h_{\rho_0}}<M_0,\ \mp \bar\xi >0\big\}, \\
		U_2\coloneqq \Big\{(\rho_0,y,\bar\xi,\eta)\in \bSM; \ \rho_0<\epsilon, \ |\eta|_{h_{\rho_0}}>\frac12\Big\},\label{eq:U2}
\end{aligned}
\end{equation}
where $M_0\gg 1$. 
We write $U_1\coloneqq U_1^+\cup U_1^-$.
On $U_1$ one can use coordinates $(\rho_0,y,\eta)$ and on $U_2$ one can use $(\theta,y,\eta)$, where $\theta\in [0,\pi]$ is determined by $\sin(\theta) = \rho_0|\eta|_{h_{\rho_{0}}},$ $ \cos (\theta)=\bar{\xi}.$
Denote $\hat{U}_2\coloneqq U_2\sqcup \tilde{\mathsf{Z}}$ and define a bijection 
\begin{equation}
 	\varphi\colon \hat{U}	_2\to \mathcal{U}_2\coloneqq 
 	\{(\theta,\delta, y, \omega) :\theta\in [0,\pi],\delta \in [0,2), (y,\omega)\in  S^*_{h_{\delta \sin (\theta)}}\partial M\},\label{eq:U2hat}
 \end{equation} 
 where $S^*_{h_{\delta\sin(\theta)}}\partial M=\big\{(y,\eta\cdot \d y)\in T^{*}\partial M\colon |\eta|_{h_{\delta\sin(\theta)}}=1\big\}$, by
\begin{equation}\label{eq:local_ident}
	\mathcal{U}_2\ni (\theta,\delta,y,\omega)=\begin{cases}
		\big(\arccos(\bar\xi),|\eta|_{h_{\rho_0}}^{-1},y,\frac{\eta}{|\eta|_{h_{\rho_0}}}\big), \quad &(\rho_0,y,\bar \xi,\eta)\in U_2,
		\\ 
		(t,0,y,\frac{\eta}{|\eta|_{h_{0}}}),\quad &(t,y,[\eta])\in   \tilde{\mathsf Z}.
	\end{cases}
\end{equation}
In the second case, the definition is clearly unambiguous.
The map $\varphi $ is inverted by 
\begin{equation}\label{eq:U2Param}
	\varphi^{-1}(\theta,\delta,y,\omega) = \begin{cases}
		(\delta \sin(\theta),y,  \cos(\theta),\delta^{-1}{ \omega}), &\delta>0\\
		(\theta,y,[\omega]),&\delta=0.
	\end{cases}
\end{equation}
The map $\varphi$ endows $\hat{U}_2$ with a topology, and in terms of it a sufficiently small neighborhood of $\tilde{\mathsf {Z}}$ is disjoint from $U_1$.
It is now easy to check that the ``coordinates''\footnote{with some abuse of terminology since the target is not a subset of $\R^{2n+1}$} determined by $\varphi$ on $\hat U_2$  are smoothly compatible with the ones on $U_1$, since in the overlap $U_1\cap \hat U_2$, the latter corresponds to introducing  inverted polar coordinates in $\eta$.
A computation\footnote{This is a standard change of variables computation, for which one has to use that $\frac{\tilde \rho_0}{\rho_0}\partial_{\tilde \rho_0}\rho_0=1+O(\rho_0)$ and $\partial_{\tilde \rho_0}y=O(\rho_0)$, where the error terms are smooth. The former can be checked by writing $\tilde \rho_0=e^h \rho_0$ for some $h\in C^\infty(M)$.
 The latter, by writing $\partial_{\tilde \rho_0}y^{j}=\d y^{j}(\grad_{\tilde\rho_0^{2}g} \tilde \rho_0)=e^{-h}\d y^j(\grad_{\rho_0^{2}g}  \rho_0)+O(\rho_0)=e^{-h} \partial_{\rho_0 }y^j+O(\rho_0)=O(\rho_0)$ for $j=1,\dots,n$. 
} shows that if $\tilde \varphi$ is defined as above but with respect to different base coordinates $(\tilde \rho_0, \tilde y)$ on $M$, 
$\tilde{\varphi}\circ \varphi^{-1}$ is a smooth diffeomorphism.
Thus $U_1^{\pm}$, $\hat{U}_2$, together with a smooth charts for the region where $\rho_0>\epsilon/2$, define a smooth structure on $\bbSM$ as a manifold with corners up to codimension 2.

It remains to show that \eqref{eq:ext_diffeo} is a diffeomorphism. 
Near $\mathsf Y_\mp$,  $\OSM$ can be parametrized as
\begin{equation}
 	(\rho_0,y,b)\mapsto \Big(\rho_0,y,\pm ({1-|b|^2_{h_{\rho_0}}})^{1/2}\rho_0\partial_{\rho_0}+b\cdot\rho_0\partial_{y}\Big),
 \end{equation} and $[\OSM;\mathsf Y]$ arises as the blowup of $\rho_0=0$, $b=0$. Hence near the front face, we can consider two types of charts near each connected component of $\mathrm {ff}$:
in an open set $V_1\subset [\OSM;\mathsf{Y}]$ disjoint from  $\beta_0^{*}(\partial \OSM)$
we can use projective coordinates $(y,\rho_0,\tilde b)=(y,\rho_0,b/\rho_0)$, in terms of which $\rho_0=0$ corresponds to the front face $\mathrm{ff}$ and $\tilde b$ restrict to coordinates on it. 
In an open set $V_2\subset [\OSM;\mathsf{Y}]$ where $\frac{|(\rho_0,b)|}{|b|_{h_{\rho_0}}}$ is bounded above we can use ``coordinates'' $(y,|b|_{h_{\rho_0}},\rho_0/|b|_{h_{\rho_0}},\hat{b}=b/|b|_{h_{\rho_0}})$, where we obtain a true chart near each point by choosing $n-2$ of the $\hat{b}$. 
Now $|b|_{h_{\rho_0}}$ is a defining function for $ \mathrm{ff}$ and $\rho_0/|b|_{h_{\rho_0}}$ is a defining function for $\beta_0^{*}(\partial \OSM)$.

It follows from \eqref{eq:diff0b:coord} that on the image of $ U_1^{\circ}$ under  $\Psi_g$, $|b|_{h_{\rho_0}}/\rho_0$ is bounded, so the first coordinate system is valid there. 
We have, using the parametrization of $U_1^{\pm }$ mentioned above and \eqref{eq:diff0b:coord},
\begin{equation}\label{eq:Psi_diffeo_1}
	(y,\rho_0,\tilde b)\circ \Psi_g( \rho_0,y,\eta )=(y,\rho_0,h_{\rho_0}^{-1}\eta),
\end{equation}
which extends to be  smooth and smoothly invertible, down to $\rho_0=0$. 
This also shows that the inverse $\Psi_g^{-1}$ extends smoothly from the interior  all the way to   $\mathrm {ff}\subset [\OSM;\mathsf Y]$ away from the face $\beta_0^{*}(\partial \OSM)$. 
 On the other hand, the image of $U_2$ under $\Psi_g$ is contained in a set where $\frac{\rho_0}{|b|_{h_{\rho_0}}}$ is bounded above, and we have, using \eqref{eq:local_ident},
\begin{equation}\label{eq:Psi_diffeo_2}
(y,|b|_{h_{\rho_0}},\rho_0/|b|_{h_{\rho_0}},\hat{b})\circ \Psi_g\circ \varphi^{-1}(\theta,\delta,y,\omega)= 
	\Big(y,\sin(\theta),\delta,h_{\delta \sin(\theta)}^{-1}\omega\Big), \quad \delta>0,
\end{equation}
which extends smoothly down to $\delta=0$, in particular down to $\partial(\overline{\partial_\pm\bSM})$.
Moreover, from this formula we see that the inverse $\Psi_{g}^{-1}$ extends smoothly from the interior up to the boundary hypersurfaces of $[\OSM;\mathsf{Y}]$ on $V_2$.
Hence we obtain a smooth extension $\overline\Psi_g\colon \bbSM\to [\OSM;\mathsf Y]$ with smooth inverse, which implies that \eqref{eq:ext_diffeo} is a diffeomorphism.
In particular, $\overline\Psi_g\colon\overline{\partial _\pm \bSM}\to \mathrm{ff}_\pm $ is a diffeomorphism.
Since the smooth structure on $[\OSM;\mathsf{Y}]$ is canonical and $\overline\Psi_g\colon S^*M^{\circ}\to SM^{\circ}$ is a diffeomorphism, there can be at most one smooth structure on $\bbSM$ for which \eqref{eq:ext_diffeo} holds.
The statement regarding the antipodal involutions is easy to see.
\end{proof}

\subsection{Sasaki metric on Unit Tangent Bundle}\label{sec:sasaki}

Let $(M^\circ,g)$ be AH of dimension $n+1$.
As is well known $g$ induces a natural Riemannian metric on the tangent bundle $TM^\circ$, called the \textbf{Sasaki metric} and denoted by $G$ throughout. We use the same letter for its restriction to $SM^\circ$.
Its construction can be found in \cite{Paternain1999}. 
The Sasaki metric induces a volume form $\d V_G$ on $SM^{\circ}$.
Its pullback to $S^*M^\circ$ via the musical isomorphism induces the \textbf{Liouville measure} $\d \Sigma^{2n+1}$,  also given by $\d \Sigma^{2n+1}{=}\frac{1}{n!}|\alpha\wedge (\d \alpha)^{ n}|$, where $\alpha$ is the tautological 1-form on $T^*M^{\circ}$ (see e.g., \cite[Exercise 1.33]{Paternain1999}).
{In coordinates $(z^\alpha,\zeta_\beta\d z^\beta)\mapsto(z^\alpha,\zeta_\beta)$ of $T^*M^\circ$ we have that} $\d \Sigma^{2n+1}(z,\zeta)=|\d V_g(z)|\d \mu_z(\zeta)$, where $\d \mu_{z}$ is the measure induced on the fibers of $S^{*}_zM^{\circ}$ by the Lebesgue measure on  $T_z^*M^{{\circ}}$ -- see \cite[Section 3.6.2]{Paternain2023}. 
To identify the behavior of $\d \Sigma^{2n+1}$ on $\bbSM$ near $\tilde {\mathsf{Z}}$,
first note that in  terms of usual coordinates on {${}^bT^*M$}, $\d \mu_z$ is induced by restricting to $S^*_{(\rho_0,y)}M^{\circ}$ the form
\begin{align}
(\bar \xi &\partial_{\bar \xi}+\eta\cdot \partial_\eta)\lrcorner \Big(\rho_0^{n}(\det h_{\rho_0})^{-1/2} \d \bar \xi \wedge \d \eta_1\wedge\dots\wedge \d \eta_n\Big)\label{eq:measure1}\\*
&= (\det h_{\rho_0})^{-1/2}\rho^{n}_0\big(\bar\xi\d \eta_1\wedge\dots\wedge \d \eta_{n}+\sum_{j=1}^{n}(-1)^{j}\eta_j\d \bar\xi \wedge \d \eta_1\wedge\dots\wedge \widehat{\d \eta_j}\wedge \dots\wedge \d \eta_{n}\big),\label{eq:measure2}
\end{align}
where $\widehat{\cdot}$ denotes omission.
Now consider ``coordinates'' $(\theta,\delta,y,{\omega})$ for an interior point in  $\bbSM$ near $\tilde{\mathsf{Z}}$, as in \eqref{eq:local_ident}.
The fiber can parametrized by $(\theta, \omega)$ and with the help of \eqref{eq:measure1} we can see that if $n\geq 2$ the fiber measure becomes
\begin{equation}\label{eq:fiber_measure}
	\d \mu=\sin^{n-1}(\theta)\d \theta \d \lambda( \omega),
\end{equation} where $\d \lambda$ is the measure induced by the Lebesgue measure on the unit cosphere $S_{h_{\rho_0}}^*\partial M$ of $h_{\rho_0}$ on $T^*\partial M$ (this can be most easily checked by first changing $\eta\mapsto (r,\omega)$, $r=|\eta|_{h_{\rho_0}}$, so that $ \eta\cdot \partial_\eta = r\partial_r$ and $ (\det h_{\rho_0})^{-1/2} \d \eta = r^{n-1}\d r \d \lambda $, and then using that $r=\rho_0^{-1}\sin(\theta)$, $\bar{\xi}=\cos(\theta)$, with $\rho_0$ fixed).
If $n=1$, the part of  the fiber of $S^{*}M^{\circ}$ contained in a sufficiently small neighborhood  $\tilde {\mathsf Z}\subset \bbSM$ has two connected components, and  $\d \mu=\d \theta $ in each of them. 
Combining with the measure $(\det h_{\rho_0})^{1/2}\rho_0^{-n-1}\d \rho_0\d y$ of the base, in which we substitute $\rho_0=\delta \sin(\theta)$, we find
\begin{equation}\label{eq:liouville}
	\d \Sigma^{2n+1}=\frac{(\det h_{\rho_0})^{1/2}}{\delta^{n+1}\sin(\theta)}\d \theta\,\d y\,\d \delta\, \d \lambda( \omega),
\end{equation}
where $\d \lambda( \omega)$ is replaced by $1$ in the case $n=1$.

As discussed in \cite[{Lemma 2.2}]{Graham2019} {and the paragraph before it}, $\partial \bSM$ carries a canonical symplectic structure.
Moreover, $X\lrcorner(\frac{1}{n!}\alpha \wedge \d \alpha)$ extends smoothly to $\partial \bSM$, and the restriction of $X\lrcorner(\frac{1}{n!}\alpha \wedge \d \alpha)$ to $\partial \bSM\cong T^*\partial M$ induces the symplectic volume form there, whose associated measure we denote by $\d \Sigma_\partial$.

\subsection{The Scalar X-Ray Transform and its Mapping Properties}\label{sub:the_scalar_x_ray_transform_and_its_mapping_properties}

Let $(M^{\circ},g)$ be nontrapping AH, of dimension $n+1$. 
We recall some facts from \cite{Graham2019}.
As mentioned previously, the geodesic vector  field $X$ on $S^{*}M^{\circ }$ extends to a smooth vector field on $\bSM $ of the form $X=\rho \overline{X}$, where $ \overline{X}$ is smooth on $\bSM$ and transversal to $\partial \bSM$ and $\rho$ is a bdf. 
The flow of $\overline X$ in $\bSM$, denoted  by $\overline{\varphi}_\tau$, is incomplete and related to the flow $\varphi_t$ of $X$ over the interior of $\bSM$ via
\begin{equation}
	\overline{\varphi}_\tau(z,\zeta)=\varphi_{t(\tau,z,\zeta)}(z,\zeta) \text{ with }t(\tau,z,\zeta)\coloneqq \int_0^\tau \frac{1}{\rho(\overline{\varphi}_s(z,\zeta))}\d s.
\end{equation}
We denote by $0\leq \tau_\pm(z,\zeta)<\infty$ the smallest times  such that
$\overline{\varphi}_{\pm \tau_{\pm}(z,\zeta)}(z,\zeta)\in \partial_\pm \bSM$
 and define the X-ray transform  as the operator 
\begin{equation}
	I\colon \rho C^\infty(\bSM)\to C^\infty (\partial_-\bSM), \quad If(z,\zeta)=\int_{-\infty }^{\infty } f({\varphi }_{t}(z^\circ,\zeta^\circ))\d t=\int_0^{	\tau _+(z,\zeta)} \frac{f(\overline\varphi_\tau(z,\zeta))}{\rho(\overline\varphi_\tau(z,\zeta))}\d \tau,
\end{equation}
where $(z^\circ,\zeta^\circ)$ is any interior point on the curve $\tau\mapsto\overline\varphi_\tau(z,\zeta)$.
We let 
\begin{equation}
	I_0\colon  \rho C^\infty(M)\to C^\infty (\partial_-\bSM),\quad I_0\coloneqq I\circ \pi_b^{*}
\end{equation}
where $\pi_b:{}^b\!S^*M\rightarrow M$ is the projection.

Further{more}, let 
\begin{equation}
	B_\pm:\bSM\to \partial_\pm\bSM, \qquad B_\pm (z,\zeta)=\overline{\varphi}_{\pm \tau_{\pm}(z,\zeta)}(z,\zeta).
\end{equation}
The maps $B_{\pm}$ are smooth and $B_{\pm}=\mathrm{Id}$ on $\partial_\pm \bSM$. 
The \textbf{Scattering map} is given by 
\begin{equation}
	S_g\colon \partial \bSM\to \partial \bSM, \qquad S_g\big|_{\partial _\mp \bSM}\coloneqq B_\pm \big|_{\partial _\mp\bSM },
\end{equation}
and it is a diffeomorphism $\partial _\mp \bSM\to \partial _\pm \bSM$, see \cite[Proposition 2.6]{Graham2019}.

The formal adjoint of $I$ in the functional setting $L^2(\bSM,\d \Sigma^{2n+1})\to L^2(\partial _-\bSM , \d \Sigma_\partial)$  is the extension as a first integral:
\begin{equation}
	I^*: C^\infty (\partial_-\bSM)\to C^{\infty }(\bSM), \quad I^* u \coloneqq B_-^{*}u.
\end{equation}
This follows from \cite[Lemma 3.5]{Graham2019}. We finally let
\begin{equation}
	I_0^{*}\colon C^\infty (\partial_-\bSM)\to C^\infty(M^{\circ}), \quad I_0^*u(z)=
	\int_{S_z^{*} M^{\circ}} I^*u(z,\zeta)\d \mu_{z}(\zeta).
\end{equation}
It follows from the above and the local product splitting of the Sasaki volume form that the operator $I_0^*$ is the formal adjoint of $I_0$ in the functional setting $L^2(M,\d V_g)\to L^2(\partial_-\bSM, \d \Sigma_\partial)$, justifying the notation.

We would like to refine the mapping properties with the help of the compactified spaces $\bbSM$ and $\overline {\partial_-\bSM}$ of Proposition~\ref{lm:compactifications}.
In the rest of the section we will make heavy use of terminology and results presented in Appendix~\ref{sec:polyhomongeneity}.
We will need the following lemma.

\begin{lemma}\label{lm:bfibrations}
If $(M^{\circ},g)$ is nontrapping AH of dimension $n+1$, the maps $B_\pm :\bSM\to \partial_\pm \bSM$ and $\pi_b:\bSM \to M$ extend to smooth b-fibrations
	\begin{align}
			B_\pm \colon &\bbSM\to \overline{\partial_\pm \bSM}, \qquad 
			\pi_b \colon \bbSM\to  M.
	\end{align}
\end{lemma}
\begin{proof}We will check the requirements of Definition \ref{def:b-subm}. In both cases, condition \ref{item2} there is automatically satisfied since neither of the target spaces has  corners of codimension $\geq 2$, in other words it suffices to show that the two maps extend to b-submersions.

We check first that  $B_-$ satisfies the requirements of a b-submersion over compact subsets of $ \bSM$ (the case of $B_+$ is analogous).
If $K\subset \bSM$ is compact, $B_-(K)\subset \partial \bSM$ is also compact by continuity, thus disjoint from the boundary of $\overline{\partial_-\bSM}$.
Therefore \eqref{eq:bmap} is satisfied with $e_{B_-}(i,j)=0$ over compact subsets of $\bSM$.
Since $B_-\colon \bSM\to \partial_-\bSM$ is smooth and $B_-\big|_{{\partial_-\bSM}}=\mathrm{Id}$, near $\partial _-\bSM$ it can be written in terms of the coordinates $(\rho_0, y,\eta)$ on $U_1^{+}$, see \eqref{eq:U2}, (with $M_0$ there large enough) as $B_-(\rho_0,y,\eta)=(y,\eta) +O(\rho_0)$. 
Since in these coordinates the sections of ${}^bT(\bSM)$ are spanned by $\rho_0\partial_{\rho_0},$ $\partial _y$, $\partial _\eta$ and those of ${}^bT\big(\overline{\partial_-\bSM}\big) $ are spanned by $\partial _y$, $\partial _\eta$,
the matrix of $\d ^{b}B_-$ takes the form $\begin{bmatrix}
	0 & \mathrm {Id} & 0\\
	0 & 0 & \mathrm{Id}
\end{bmatrix}+O(\rho_0)$ .
It has full rank for small enough $\rho_0$, thus  $\d ^{b}B_-|_{(z,\zeta)}$ is surjective for $(z,\zeta)$ in a neighborhood of $\partial_-\bSM$ in $\bSM$. 
If ${(z,\zeta)}\in {\partial_+\bSM}$, then $\d^{b} B_-|_{(z,\zeta)}$ is surjective because $B_-|_{{\partial_+\bSM}}=S_g\colon  \partial_{+}\bSM\to \partial_- \bSM$ is a diffeomorphism.
Finally for arbitrary $(z,\zeta)\in (\bSM)^{\circ}$, we can write $B_- (z,\zeta)= B_-(\overline{\varphi}_{-\tau} (z,\zeta))$ for some $\tau< \tau_-(z,\zeta)$ sufficiently large to guarantee that $\d ^{b}B_-\big|_{\overline{\varphi}_{-\tau} (z,\zeta)}$ is surjective. 
Since $\overline{\varphi}_{-\tau}$ is a diffeomorphism, $\d ^{b}B_-\big|_{(z,\zeta)}$  must be surjective as well.
Hence $B_-$ satisfies the requirements for being a b-fibration  over compact regions of $\bSM$.

We now analyze $B_-$ near $\tilde{\mathsf{Z}}\subset \bbSM$,
 using ``coordinates'' $(\theta,\delta,y,\omega)$, see \eqref{eq:local_ident}, on the open set $\hat{U}_2\subset \bbSM$.
Consider $ \tilde{X}=|\eta|_{h_{\rho_{0}}}^{-1}\overline{X}=\delta \overline{X}$, which using \cite[eq. 2.22]{Graham2019} can be written as
\begin{equation}\label{eq:Xtilde}
	\tilde{X}=\big(1+\delta \sin(\theta)Q_1\big)\partial_\theta+ Q_2\delta^2\partial_\delta+
	P_1^{i}\delta\sin(\theta)\partial_{y^i}+ P_2^{j}\delta\partial_{\omega_{j}},
\end{equation}
where the $Q_j$, $P_j^i$ are all smooth up to $\theta=0,\pi$ and $\delta=0$.
In the last term, with some abuse of notation, we mean that near a  point of interest we choose $n-1 $ of the $\omega_{j}$'s that form a coordinate system together with $(\theta,y,\delta)$.
The vector field $\tilde{X}$
 is a rescaling of $\overline{X}$, so their integral curves $(\theta(s),\delta(s),y(s),\omega(s))$ agree up to reparametrization.
In $\hat U_2$, $\tilde{X}$ is smooth, transversal to $\overline{\partial\bSM}$ (which corresponds to $\{\theta=0\}\cup \{\theta=\pi\}$) and tangent to the boundary hypersurface $\tilde{\mathsf{Z}}$, which corresponds to $\delta=0$.
Hence it has a smooth incomplete flow $\tilde{\varphi}_s$, and the forward and backward exit times corresponding to an initial point $(\theta_*,\delta_*,y_*,\omega_*)\in\hat{U}_2$, defined by
\begin{equation}
	s_\pm(\theta_*,\delta_*,y_*,\omega_*)\coloneqq\sup \big\{s\geq 0:\tilde \varphi_{\pm s}(\theta_*,\delta_*,y_*,\omega_*)\in \bbSM\big\}<\infty,
\end{equation}
 are given implicitly by $\theta\circ \tilde{\varphi}_{ s_+ }(\theta_*,\delta_*,y_*,\omega_*)=\pi$, resp. $\theta\circ \tilde{\varphi}_{- s_- }(\theta_*,\delta_*,y_*,\omega_*)=0$. 
An application of the implicit function theorem shows that $s_\pm $ is smooth, up to $\delta_*=0$.
Since $B_-(\cdot)=\tilde{\varphi}_{-s_-(\cdot)}(\cdot)$, $B_-$ extends to be smooth on $\bbSM$, up to $\delta_*=0$.
We will analyze the flow $\tilde{\varphi}$ further near ${\delta}_*=0$, to show that 
\begin{equation}\label{eq:B-bmap}
 	B_-^{*}{\delta }={\delta }_*+O({\delta }_*^{2}),
 \end{equation}where the error term is smooth; since $\delta$ and $\delta_*$ are defining functions for $\partial (\overline {\partial_-\bSM})$ and $\tilde{\mathsf{Z}}\subset \bbSM$ respectively, \eqref{eq:B-bmap} will imply that $B_-$ is a b-map.

We  check that for the smooth map $F(s,\theta_*,\delta_*,y_*,\omega_*)\coloneqq{\delta}\circ\tilde{\varphi}_{s}(\theta_*,\delta_*,y_*,\omega_*) $ we have $F\big|_{\delta_* =0}=0$ and 
$\partial_{\delta_*}F\big|_{\delta_* =0}\neq 0$.
Indeed, for $n\geq 2$, solving the ODE system 
\begin{equation}\label{eq:explicit_system}
	\frac{\d \theta}{\d s}=1, \qquad \frac{\d {\delta}}{\d s}=0, \qquad \frac{\d y}{\d s}=0,\qquad \frac{\d \omega}{\d s}=0
\end{equation}we see that
\begin{equation}
	\tilde{\varphi}_s(\theta_*,0,y_*,\omega_*)=(\theta_*+s,0,y_*,\omega_*),\label{eq:sol_l_0}
\end{equation} so that $F(s,\theta_*,0,y_*,\omega_*)=0$.
Since $\omega\circ \tilde \varphi_{(\cdot)}(\cdot)$ is smooth in  all variables and equal to $\omega_*$ for $\delta_*=0$, we see that $\omega\circ \tilde \varphi_s(\theta_*,\delta_*,y_*,\omega_*) = \omega_*+O(\delta_{*})$, where the error is uniform in $s$ by the smoothness of $s_\pm$.
Hence if $\delta_{*}$ is sufficiently small, we can choose $n-1$ of the $n$ variables $\omega$, in terms of which the solution can be written for all $s\in [-s_-,s_+]$; we still keep the notation $\omega$.
If $n=1$, $\omega\in \{\pm h_{\rho_0}^{1/2}\}$, and similarly for $\omega_*$, and the discussion below carries over with the obvious modifications.  
By \cite[Ch. 1, Theorem 7.5]{Coddington1955} and the subsequent discussion, since $\tilde\varphi_s$ is given as the solution of a system of the form
\begin{equation}
 	\frac{\d }{\d s}\tilde\varphi_s=f(\tilde\varphi_s), \qquad \tilde\varphi_0(\theta_*,\delta_*,y_*,\omega_*)=(\theta_*,\delta_*,y_*,\omega_*),\quad f(u)=f(\theta,\delta,y,\omega)\text{ smooth},
 \end{equation} 
 its matrix of partial derivatives $\displaystyle \Phi_s\coloneqq \frac{\partial (\theta(s),\delta(s),y(s),\omega(s))}{\partial (\theta_*,\delta_*,y_*,\omega_*)}$
  satisfies the ODE system
 \begin{equation}
 	\frac{\d }{\d s}\Phi_s=\partial_uf(\tilde\varphi_s)\Phi_s
  \overset{ \eqref{eq:Xtilde}}=\left.\left(\begin{bmatrix}
  	0 & \star & 0 & 0\\ 
  	0 & 0 & 0 & 0\\ 
  	0 & \star & 0 & 0\\
  	0 & \star & 0 & 0
  \end{bmatrix}+O(\delta)\right)\right|_{{\tilde\varphi_s}}
  \Phi_s,\qquad \Phi_0=\mathrm{Id}.
 \end{equation}
Here and below, ``$\star$'' denotes smooth quantities whose exact value is not relevant for the argument. 
Since $\delta\circ \tilde \varphi_s\big|_{\delta_*=0}=0$, (implying that the second row of $\Phi_s\big|_{\delta_*=0}$ is of the form $\begin{bmatrix}
	0 & \star  & 0 & 0 
\end{bmatrix}$),
\begin{equation}\label{der_t=0}
	  \Phi_s\big|_{\delta_*=0}
  =\begin{bmatrix}
    	1 & \star & 0 & 0\\ 
    	0 & 1 & 0 & 0\\ 
    	0 & \star & \mathrm{Id} & 0\\
    	0 & \star & 0 & \mathrm{Id}
    \end{bmatrix}.
\end{equation}
From this we see that  $ \partial_{\delta_*}\delta \big|_{\delta_*=0}=1$, which means  $\partial_{{\delta}_*}F\big|_{{\delta}_*=0}=1$.
Now compute
\begin{equation}\label{eq:B_-der}
	\partial_{{\delta}_*}\big({\delta}\circ B_-(\theta_*,\delta_*,y_*,\omega_*)\big)\big|_{{\delta}_*=0}= 
	\partial_{{\delta}_*}\big(F(-s_-(\cdot ),\cdot )\big)\big|_{{\delta}_*=0}
	=\partial_s F\big |_{{\delta}_*=0} (-\partial_{{\delta}_*}s_-)\big |_{{\delta}_*=0}+\partial_{{\delta}_*}F\big |_{{\delta}_*=0}.
\end{equation}
Since $F\big |_{{\delta}_*=0}=0$, the first term vanishes and $\partial_{{\delta}_*}({\delta}\circ B_-(\theta_*,\delta_*,y_*,\omega_*))\big|_{{\delta}_*=0}=1$, proving \eqref{eq:B-bmap}.
We conclude that $B_-$ is a b-map.

We check that $B_-$ is a  b-submersion. The sections of ${}^{b}T (\bbSM)$ are locally spanned by  $\sin (\theta_*)\partial_{\theta_*}$, $\partial_{y_*}$,
$\delta _{*}\partial_{\delta_*}$, and $\partial_{\omega_{*}}$, whereas those of $ {}^bT\big(\overline{\partial _-\bSM}\big)$ are spanned by ${\delta}\partial_{{\delta}} $, $\partial_y$ and $\partial_{\omega}$.
It suffices to check that $\d^b B_-$ is surjective when $\delta_*=0$, because when $\delta_*>0$ the arguments in the second paragraph of this proof apply.
Using \eqref{eq:B-bmap}, \eqref{der_t=0}, and computations similar to \eqref{eq:B_-der},
the matrix of $\d ^{b}B_-$ in the bases $\{\delta_{*}\partial_{\delta_*},\partial_{y_*},\partial_{\omega_*}, \sin(\theta_*)\partial_{\theta_*}\}$ and $\{{\delta}\partial_{{\delta}} ,\partial_y,\partial_\omega\}$ (note the change in the order of the coordinates)
  becomes
  \begin{equation}
  	\begin{bmatrix}			
  		\frac{\d \delta}{\delta}(\d B_-({\delta}_*\partial_{{\delta}_*}))  & \star & \star &\star\, \\
  		{\d y}(\d B_-({\delta}_*\partial_{{\delta}_*})) & {\d y}(\d B_-(\partial_{y_*})) & \star  &\star\, \\
  		{\d \omega }(\d B_-({\delta}_*\partial_{{\delta}_*})) & {\d \omega }(\d B_-(\partial_{y_*})) & {\d \omega }(\d B_-(\partial_{\omega _*})) &\star\, 
  	\end{bmatrix}		=
  	\begin{bmatrix}
  		1 & \star &\star &\star \\
  		0 & \mathrm{Id} & \star & \star \\
  		0 & 0 & \mathrm{Id} & \star
  	\end{bmatrix}+O(\delta_*),
  \end{equation}
which has full rank for $\delta_* =0$, so $ \d ^{b}B_-$ is surjective for any point with ${\delta}_*=0$.
This concludes the proof of the fact that $ B_-$ is a b-submersion, hence also a b-fibration as already explained.

We now turn to  $\pi_b$. We know that $\rho_0$ is a bdf for $\partial M$. To avoid confusion we write it momentarily as $\tilde{\rho}_0$, and use $\rho_{0}$ for the function $\rho_0=\pi_b^{*}\tilde\rho_{0} $ on $\bSM$ appearing in \eqref{eq:U2}. 
We will write the coordinates in $ M$ as $(\tilde{\rho}_0,\tilde{y})$.
Since $\pi_b$ expressed in terms of the coordinates $(\rho_0, y, \eta)$ in $U_1$ as $\tilde\rho_{0}=\rho_{0}$, $\tilde y=y$, it is clearly a b-map there.
On $\hat{U}_2$ we have $ \pi_b^{*}\tilde{\rho}_{0}=\rho_{0} =\delta \sin(\theta)$.
Since $\sin(\theta)$ and $\delta $ are defining functions for $\overline{\partial \bSM}$ and $\tilde{\mathsf{Z}}$ respectively, $\pi_b$ extends to a b-map on $\bbSM$.
We check now that $\pi_b$ is a b-submersion. We use the coordinates $\delta, {y},\theta,\omega $  in $\hat U_2$, in this order (in the analysis of $B_-$ these had $*$ subscripts, which are now unnecessary).
As before, the sections of ${}^{b}T( \bbSM)$ are locally spanned by $\delta  \partial_{\delta }$, $\partial_{y}$, $\sin(\theta)\partial_\theta$ and   $\partial_\omega $, whereas those of ${}^{b}T M$ are spanned by $\tilde{\rho}_{0}\partial_{\tilde{\rho}_{0}},$ $\partial_{\tilde{y}}$ and the matrix of $\d^{b}\pi_b$ on $\hat U_2$ becomes
\begin{equation}
	\begin{bmatrix}
		\frac{\d \tilde{\rho}_0}{\tilde{\rho}_0}(\d \pi_b (\delta \partial_\delta )) & \star & \star & \star\ \\
		\d \tilde{y}(\d \pi_b (\delta \partial_\delta )) & \d \tilde{y}(\d \pi_b (\partial_y)) & \star & \star \
	\end{bmatrix}=\begin{bmatrix}
		1 & \star & \star &\star \ \\
		0 & \mathrm{Id}&  \star & \star\ 
	\end{bmatrix}.
\end{equation}
This matrix has full rank on $\hat U_2$, and it is not hard to check that the same is true for the matrix of $\d ^{b}\pi_b$ expressed in coordinates on $U_1^\pm$. We conclude that $\pi_b$ is a b-submersion, which  shows as before that it is also a b-fibration.
\end{proof}

\bigskip

Now that we have established the b-fibration property of $B_-$ and $\pi_b$, we can show mapping properties of $I$, $I_0$, $I^{*}$ and $I_0^* $, inspired by the double b-fibration picture:
\begin{center}
\begin{tikzcd}[column sep=1.5em, row sep=3em]
    & \bbSM \arrow[dl, "\pi_{b}"'] \arrow[dr, "B_-"] \\
    M & & \overline{\partial_-\bSM}
\end{tikzcd}.\label{eq:double_b_fibr}
\end{center}
The double fibration approach in the study of the geodesic X-ray transform is classical (see e.g., \cite{Guillemin1977}).
More recently, double b-fibrations were used in the study of the geodesic X-ray transform on compact manifolds with boundary (\cite{Mazzeo2021}).

To streamline the statements below, we introduce the following numbering of the boundary hypersurfaces of the spaces appearing in the double b-fibration picture above.
We write $\rho_1=\rho$, which is a bdf for $M$.
Assume for now that $n+1\geq 3$, in which case let $x_1,$ $x_2,$ $x_3\in C^\infty(\bbSM)$ be bdfs for $\tilde{\mathsf{Z}}$, $\overline{\partial_-\bSM}$, $\overline{\partial_+\bSM},$  respectively. 
Finally we let $x_{1,\pm }\in C^\infty(\overline{\partial\bSM})$ be defining functions for  $\partial(\overline{\partial_\pm \bSM})= \overline{\partial_\pm \bSM}\cap \tilde{\mathsf{Z}}$.
We can take $ x_{1,\pm }=x_1\big|_{\overline{\partial_\pm \bSM}}$.
All statements involving families of index sets on these manifolds with corners will respect this numbering of the boundary hypersurfaces.
Since $\pi_b^{*}\rho_1=\delta \sin(\theta)$ in the coordinates \eqref{eq:local_ident},
\begin{equation}\label{eq:pi_bmap}
	\pi_b^{*}\rho_1= e^{h}\prod_{j=	1}^{3} x_j, \quad h\in C^\infty(\bbSM),
\end{equation}
that is, $e_{\pi_b}(1,j)=1$ for $j=1,2,3$ in the notation of Definition \ref{def:b-subm}. %
Moreover, from \eqref{eq:B-bmap},
\begin{equation}\label{eq:B-lifting}
	B_-^*(x_{1,	-})= e^{h'}\prod_{j=1}^{3} x_j^{e_{B_-}(1,j)},\quad  h'\in C^\infty(\bbSM),\quad   e_{B_-}(1,1)=1,\  e_{B_-}(1,2)=e_{B_-}(1,3)=0.
\end{equation}

As remarked before, in the case $n+1=2$, both $\tilde{\mathsf{Z}}$ and $\partial(\overline{\partial_\pm \bSM})$ each have two connected components, and we can let $x_1$, resp. $x_{1,\pm }$ be a defining function for their union.
Then the statements above work in exactly the same way, and when we attach a single weight or index set to $\tilde{\mathsf{Z}}$ or $\partial(\overline{\partial_\pm \bSM})$, we will mean that the relevant conormal or polyhomogeneous functions have the same weight or index set on both connected components.

In Proposition~\ref{prop:phg_mapping} below we extend the X-ray transform to polyhomogeneous and conormal functions in the obvious way and show mapping properties on these spaces,
  analogous to Theorems 2.1 and 2.2 in \cite{Mazzeo2021}.
  There, sharp polyhomogeneous mapping properties are proved for the X-ray transform and the backprojection operator  on compact manifolds with boundary.
  Items \ref{mapping_item2} and \ref{mapping_item3} could have been formulated  a bit more generally in the case $n+1=2$ by allowing different weights and index sets in the connected components of $\tilde{\mathsf{Z}}$ and of $\partial(\overline{\partial_\pm \bSM})$ (with an analogous modification of the conclusion), but we had no concrete application for this and opted for a unified presentation.

\begin{proposition}\label{prop:phg_mapping}
Let $(M^\circ,g)$ be nontrapping, $n+1$-dimensional, with $g$ AH, $n\geq 1$.
\begin{enumerate}[leftmargin=0.7cm, label=(\roman*)]
	\item\label{mapping_item0} If $f\in \mathcal{A}_{-}^{\mathfrak{t}}(\bbSM)$ with $\mathfrak{t}=(t_1,t_2,t_3)$, $t_2,t_3>0$, then $If\in \mathcal{A}_{-}^{t_1}\big(\overline{\partial_-\bSM}\big)$.
	If in addition $f\in \Aphg^{\mathcal E}(M)$ for some index family $\mathcal {E}=(E_1,E_2,E_3)$ with $\Re(E_2),\, \Re(E_3)>0$, then $I f\in \Aphg^{E_1}\big(\overline{\partial_-\bSM}\big)$.

	\item\label{mapping_item1} If $\tilde f\in \mathcal{A}_{-}^{t}(M)$ with $t>0$, then $I_0\tilde f\in \mathcal{A}_{-}^{t}\big(\overline{\partial_-\bSM}\big)$.
	If in addition $\tilde f\in \Aphg^{E}(M)$ for some index set with $\Re(E)>0$, then $I_0 \tilde f\in \Aphg^{E}\big(\overline{\partial_-\bSM}\big)$.

	\item\label{mapping_item2} Let $h\in \mathcal{A}^{t}_{-}\big(\overline{\partial_-\bSM}\big)$, resp. $h\in \mathcal S^{t}\big(\overline{\partial_-\bSM}\big)$, for some  $t\in \R$. Then $I^{*}h\in \mathcal{A}_{-}^{\mathfrak{t}}(\bbSM)$, resp. $I^{*}h\in \mathcal S^{\mathfrak{t}}\big(\bbSM\big)$, where $\mathfrak{t}=(t,0,0)$. 
	If in addition $h\in \Aphg^{G}\big(\overline{\partial_-\bSM}\big)$ for some index set $G$, then $I^{*}h\in \Aphg^{\mathcal{G}}(\bbSM)$ with $\mathcal{G}=(G,0,0)$.

	\item\label{mapping_item3} 
	 If $h$ in $\mathcal{A}^{t}_-\big(\overline{\partial_-\bSM}\big)$ for some $t\in \R$, then $I_0^*h\in \mathcal{A}^{t'}_-(M)$, where $t'=\min\{t,n\}$.
	If $h\in \Aphg^{G}\big(\overline{\partial_-\bSM}\big)$ for some index set $G$, 
	then $I_0^{*}h\in \Aphg^{G'}(M)$ with $G'=G\,\overline{\cup }G\,\overline{\cup }\overline{\{(n,1)\}}$.
\end{enumerate}

\end{proposition}

\begin{proof}
The approach will be to realize $I$, $I_0$, $I^*$, $I_0^*$ as  (compositions of) pullback/pushforward operations and use Theorems~\ref{thm:pullback} and \ref{thm:pushforward}.
We will use some notations from the proof of Lemma~\ref{lm:bfibrations}.

For Part~\ref{mapping_item0},
using a partition of unity it suffices to prove the statement by assuming that $f$ is supported in various open sets whose union covers $\bbSM$.
 If $f$ is compactly supported in $S^*M^\circ$, then 
$If\in C_c^\infty(\partial_-\bSM)$ and there is nothing to show.
Then, if $\supp f $ is compact and contained in a sufficiently small neighborhood of a point on $\partial_\mp  \bSM$ in $\bSM$ we can use the coordinates on $U_1^{\mp }$.
By \cite[eq. (2.21)]{Graham2019}, 
\begin{equation}
	\overline X =\operatorname{sign}(\bar{\xi})\sqrt{1-\rho_0^{2}|\eta|^{2}_{h_{\rho_0}}}\partial_{\rho_0}+\rho_0 h_{\rho_0}^{ij}\eta_i\partial_{y^{j}}-\frac12  \rho_0\partial_{y^{k}}|\eta|^{2}_{h_{\rho_0}}\partial_{\eta_{k}}.
\end{equation}
First we work over $U_1^-$.
For any $(y_*,\eta_*)=(0,y_*,\eta_*)\in \partial_-\bSM$ we observe that we have\footnote{We use the convention of having the subscript $*$ in the variables of which we think as initial conditions; this is why here the subscripted variables are in $\overline{\partial_-\bSM}$ instead of $\bbSM$, as was the case in the proof of Lemma~\ref{lm:bfibrations}.} $\frac{\d }{\d \tau} (\rho_0,y,\eta)\circ \overline\varphi_\tau(y_*,\eta_*)\big|_{\tau =0}=(1,0,0), $ therefore $ (\rho_0,y,\eta)\circ \overline\varphi_\tau(y_*,\eta_*)=(\tau,y_*,\eta_*)+O(\tau^2)$ and the error term is smooth.
Then write 
\begin{equation}
	If(y_*,\eta_*)=\int_0^{\epsilon } \frac{f}{\rho_0}\circ \overline{\varphi}_\tau(y_*,\eta_*)\d \tau=\int_0^{\epsilon } \frac{f((\tau,y_*,\eta_*)+O(\tau^2))}{\tau+O(\tau^2)}\d \tau.
\end{equation}
Now it is not hard to check that if $f(\rho_{0},y,\eta)$ is in either of the spaces of the statement 
then $If$ is smooth.
On the other hand, if $f$ is compactly supported in $U_1^{+}$ we can use a change of variables to write
\begin{equation}\label{eq:pushf_Uplus}
	If(y_*,\eta_*)=\int_{\tau_+(y_*,\eta_*)-\epsilon}^{\tau_+(y_*,\eta_*)} \frac{f}{\rho_0}\circ \overline{\varphi}_{\tau}(y_*,\eta_*)\d \tau= \int_{0}^{\epsilon} \frac{f}{\rho_0}\circ \overline{\varphi}_{-u}(S_g(y_*,\eta_*))\d u.
\end{equation}
Since $S_g$ is smooth, smoothness of \eqref{eq:pushf_Uplus} follows similarly to the case over $U_1^-$.

So we can now assume that $f$ is supported near $\tilde {\mathsf {Z}}$, in $\hat U_2$.
Let
\begin{equation*}
	H(s,\delta_*,y_*,\omega_*)\coloneqq\tilde{\varphi}_s(\delta_*,y_*,\omega_*)=(\theta(s),\delta(s),y(s),\omega(s)), \quad (\delta_*,y_*,\omega_*)\in \overline{\partial_-\bSM},
\end{equation*} where as before $\tilde\varphi_s$ is the flow of $\tilde X$ in \eqref{eq:Xtilde}, expressed in the ``coordinates'' \eqref{eq:local_ident} on $\hat U_2$.
Note that 
\begin{equation*}
	B_-^{-1}(\{(\delta_*,y_*,\omega_*)\})=\big\{H(s,\delta_*,y_*,\omega_*):s\in [0,s_+(\delta_*,y_*,\omega_*)]\big\}.
\end{equation*}
By \eqref{eq:sol_l_0}, $H(s,0,y_*,\omega_*)=(s,0,y_*,\omega_*)$.
Using this and \eqref{der_t=0} we can check that
\begin{equation}
	\partial_{(s,\delta_*,y_*,\omega_*)}H=\begin{bmatrix}
    	1 & \star & 0 & 0\\ 
    	0 & 1 & 0 & 0\\ 
    	0 & \star & \mathrm{Id} & 0\\
    	0 & \star & 0 & \mathrm{Id}
    \end{bmatrix}+O(\delta_*)\implies \det \partial_{(s,\delta_*,y_*,\omega_*)}H=1+O(\delta_*)\label{eq:detF}.
\end{equation}

Let $\mu_*=\frac{|\d \delta_*\wedge \d y_*\wedge  \d \lambda(\omega_*)|}{\delta_*}$,  which spans the b-density bundle on $\overline{\partial_-\bSM}$ near its boundary (as before, $\d \lambda $ denotes the fiber measure on the fibers of $S_{h_0}^*\partial M$ and is replaced by 1 if $n=1$) and $\nu= \frac{|\d \theta\wedge  \d \delta\wedge \d y\wedge\d \lambda(\omega)|}{\delta \sin(\theta)}$, {which spans} the b-density bundle in $\bbSM$.
Note that since  $\tilde X=\delta \overline{X}$, we have $\tilde{\varphi}_s=\overline \varphi_\tau$, where $\frac{\d \tau}{\d s}=\delta\circ\tilde\varphi_s$, and recall that $\rho_0=\delta \sin \theta$.
Using these facts, for  $0<\delta_*\ll 1$ and for
$f\in \mathcal{A}_{-}^{\mathfrak{t}}(\bbSM)$ (resp. $f\in \Aphg^{\mathcal{E}}(\bbSM)$) as in the statement we compute
\begin{align*}
	If\cdot\mu_*= &(If)(\delta_*,y_*,\omega_*)\cdot \mu_*\\
	=&\int_{0}^{\tau_+(\delta_*,y_*,\omega_*)}f(\overline{\varphi}_s(\delta_*,y_*,\omega_*))\frac{\d \tau }{\rho_0\circ \overline{\varphi}_\tau (\delta_*,y_*,\omega_*))}\cdot  \mu_*\\
	=&\int_{0}^{s_+(\delta_*,y_*,\omega_*)}f(\tilde{\varphi}_s(\delta_*,y_*,\omega_*))\frac{\d s}{\sin(\theta \circ \tilde{\varphi}_s(\delta_*,y_*,\omega_*))}\cdot  \mu_*\\ 
	=&\int_{0}^{s_+(\delta_*,y_*,\omega_*)}\Big(\frac{ f}{\sin(\theta)}\Big)(H(s,\delta_*,y_*,\omega_*)) {\d s}\cdot  \frac{|\d \delta_*\wedge \d y_*\wedge \d \lambda(\omega_*)|}{\delta_*}\\
	=&\frac{1}{\delta_*}(B_-)_{*}\Big( \frac{\delta f}{|\det \partial_{(s,\delta_*,y_*,\omega_*)}H|}\cdot \nu\Big),
\end{align*}
where the last line follows by the change of variables $(s,\delta_*,y_*,\omega_*)\mapsto (\theta,\delta,y,\omega)$.
By \eqref{eq:detF},
\begin{equation}
 	\frac{\delta f}{|\det \partial_{(s,\delta_*,y_*,\omega_*)}H|}\cdot \nu\in\mathcal{A}_-^{(t_1+1,t_2,t_3)}(\bbSM;\Omega_b), \quad \text{resp.}\quad  \Aphg^{(E_1+1,E_2,E_3)}(\bbSM;\Omega_b).
 \end{equation}  
By Lemma~\ref{lm:bfibrations} $B_-$ is a b-fibration, and since $t_2,\, t_3>0$ (resp. $\Re E_2,\, \Re E_3>0$), the assumptions of Theorem~\ref{thm:pushforward} are satisfied by the weights (resp. index sets) corresponding to $x_2=0$ and $x_3=0$.
Moreover, $(t_1+1,t_2,t_3)\geq B_-^{\sharp}\big(t_1+1\big)=(t_1+1,0,0)$ and $(B_-)_{\sharp }(E_1+1,E_2,E_3)=E_1+1$.
Thus
\begin{equation*}
\begin{aligned}
If\cdot\mu_*\in \delta_*^{-1}\mathcal{A}_-^{{t_1+1}}\big(\overline{\partial_-\bSM};\Omega_b\big)=\mathcal{A}_-^{{t_1}}\big(\overline{\partial_-\bSM};\Omega_b\big),\\ 
 \text{resp.}\quad 
	If\cdot\mu_*\in \delta_*^{-1}\Aphg^{E_1+1}\big(\overline{\partial_-\bSM};\Omega_b\big)=\Aphg^{E_1}\big(\overline{\partial_-\bSM};\Omega_b\big).	
\end{aligned}
\end{equation*}

\smallskip

For Part~\ref{mapping_item1}, if $n\geq 2$, Theorem~\ref{thm:pullback} implies that $\pi_b^{*}f\in \mathcal{A}_-^{\pi_b^{\sharp}t}(\bbSM)$ (resp. $\pi_b^{*}f\in \Aphg^{\pi_b^\sharp(E)}(\bbSM)$), and by \eqref{eq:pi_bmap}, $ \pi_b^\sharp t=(t,t,t)$ and $ \pi_b^\sharp(E)=(E,E,E)$.
It is not hard to check that if $n=1$ the same conclusion still holds, under the interpretation explained before the statement of Proposition~\ref{prop:phg_mapping}.
The claim now follows from Part~\ref{mapping_item0}.

\smallskip

Part~\ref{mapping_item2} is immediate from the fact that $I^{*}=B_-^{*}$, together with \eqref{eq:B-lifting} and Theorem~\ref{thm:pullback}.

\smallskip

For Part~\ref{mapping_item3}, for $h\in \mathcal {A}_{-}^{t}\big(\overline{\partial_-\bSM}\big)$, resp. $h\in \Aphg^{G}\big(\overline{\partial_-\bSM}\big)$,  we have (identifying a density with the measure it induces)\begin{equation}
\begin{aligned}\label{eq:adjoint_density}
 		I_0^{*} h(z)\cdot | \d V_g|=\int_{S^*_zM^{\circ}}\!& ( B_-^{*}h)(z,\zeta)\d \mu_z(\zeta)|\d V_g|\\
& =\int_{\pi_{{b}}^{-1}(z)} \!( B_-^{*}h)\cdot \d \Sigma^{2n+1}=(\pi_b)_*\left( (B_-^{*}h)\cdot  \d \Sigma^{2n+1}\right), \quad z\in M^{\circ}.
\end{aligned}
\end{equation}
By Part~\ref{mapping_item2}, $B_-^*h\in \mathcal{A}_- ^{(t,0,0)}(\bbSM)$, resp. $B_-^*h\in \Aphg ^{(G,0,0)}(\bbSM)$.
Since by \eqref{eq:liouville}, $\d \Sigma^{2n+1}\in x_1^{-n}C^\infty(\bbSM;\Omega_b)$, it follows that ($B_-^{*}h)\cdot \d \Sigma^{2n+1}\in \mathcal{A} ^{(t-n,0,0)}(\bbSM;\Omega_b)$, resp.
 $ (B_-^{*}h)\cdot \d \Sigma^{2n+1}\in \Aphg ^{(G-n,0,0)}(\bbSM;\Omega_b)$.
By \eqref{eq:pi_bmap}, all boundary faces of $\bbSM$ are mapped to $\partial M$ by the b-fibration $\pi_b$.
So if $n\geq 2$, wee see from  Theorem~\ref{thm:pushforward} that
$(\pi_b)_{*} \left( (B_-^{*}h)\d \Sigma^{2n+1}\right)\in \mathcal{A}_{-}^{\tilde{t}}(M;\Omega_b)$, resp. $(\pi_b)_{*} \left( (B_-^{*}h)\d \Sigma^{2n+1}\right)\in \mathcal{A}_{\mathrm{phg}}^{\tilde{G}}(M;\Omega_b)$, where 
 \begin{equation}
		\tilde t=	\min\{t-n,0\},\qquad \text{resp.}\qquad \tilde{G}=(G-n)\overline{\cup}\overline{\{(0,0)\}} \overline{\cup}\overline{\{(0,0)\}}=(G-n)\overline{\cup}\overline{\{(0,1)\}}. 	
 \end{equation}
This statement also holds for $n=1$ in the conormal case, and in the polyhomogeneous one we have instead $\tilde G=(G-1)\overline{\cup }(G-1)\overline{\cup }\overline{\{(0,1)\}}$.
Using that $|\d V_g|\in \rho_{1}^{-n}C^\infty(M;\Omega_b)$, we can divide out $|\d V_g|$ in \eqref{eq:adjoint_density} to obtain the claim.
\end{proof}

\begin{remark}\upshape
	We make no claim that the index sets in Proposition~\ref{prop:phg_mapping} are sharp (and it is clear from the proof that $G'$ in \ref{mapping_item3} is not if $n\geq 2$).
	In fact, by \cite[Prop. 3.2 and 4.4]{Eptaminitakis2022}, if $f\in \Aphg^{E}({M})$ with $\Re E>0$ then $I_0^*I_0\in \Aphg^{E'}(M)$ with $E'=E\,\overline{\cup}\,\overline{\{(n,0)\}}$. 
	This is consistent with the statement of Proposition~\ref{prop:phg_mapping}, up to additional logarithmic factors that the latter allows. 
\end{remark}

\begin{remark}\upshape

	Proposition~\ref{prop:phg_mapping}\ref{mapping_item2} implies that if $h\in C^\infty (\overline{\partial_-\bSM})$, then  $Ih\in C^\infty(\bbSM)$. 
	So we need no compatibility conditions for the initial data of geodesically invariant functions on $\bbSM$, unlike the case of a compact manifold with boundary $N$, where the initial data of geodesically invariant functions on $N$ must satisfy a compatibility condition which depends on the scattering relation.
	One way to explain this discrepancy is that the space $\bbSM$ is ``resolved'', in the sense that the maps $B_-$, $\pi_b$ become b-fibrations in their respective targets. 
	In the compact case, it should be compared with the resolved space $\partial_+SN\times[0,1]$, see \cite{Mazzeo2021}.

\end{remark}

\subsection{Vertical Fourier Analysis in Dimension 2}\label{sec:vertical_fourier}
We close \S\ref{sec:the_unit_sphere_and_cosphere_bundle_of_an_ah_manifold} with some background on fiberwise Fourier analysis on $\OSM$, where
 $(M^\circ,g)$ is an 
oriented
AH surface (i.e. $n=1$).
For any tangent vector $v\in T_zM^\circ$ we let $v^\bot\in T_zM^\circ$ denote the vector obtained by rotating $v$ by ${90}^\circ$ counterclockwise and let $v_\bot=-v^\bot$ (i.e. the ${90}^\circ$ clockwise rotation of $v$).
We let $V$ denote the \textbf{vertical vector field} which is a tangent vector field over $SM^\circ$ defined as follows.
For any $(z,v)\in SM^\circ$, letting $\mathcal{I}_{ \left(z,v \right)} : T_zM^\circ\rightarrow T_{ \left(z,v \right)}T_zM^\circ$ denote the canonical identification we set
$$V \left(z,v \right)=\mathcal{I}_{ \left(z,v \right)} (v^\bot ).$$
Next, we let $X_\bot$ denote the following vector field over $SM^\circ$.
Fixing any point $ \left(z,v \right)\in SM^\circ$, let $\gamma_{x,v_\bot} :  \left(-\infty,\infty \right)\rightarrow M^\circ$ denote the unique geodesic with initial condition $\gamma_{x,v_\bot}^\prime \left(0 \right)=v_\bot$, and let $W :  \left(-\infty,\infty \right)\rightarrow TM$ denote the parallel transport of $v$ along $\gamma_{x,v_\bot}$.
We define
$$X_\bot \left(z,v \right)= \left(\gamma_{x,v_\bot},W \right)^\prime \left(0 \right).$$

\noindent It is well known that $\left\{X,X_\bot,V\right\}$ forms a $C^\infty$ orthonormal frame for $TSM^{\circ}$ (with respect to $G$).

Now  assume also that $M$ is simply connected.
Let $\tilde \rho$ be any bdf, so that $ \left({M},\tilde\rho^2g \right)$  is a compact simply connected, oriented Riemannian manifold with boundary.
Thus by \cite[Theorems 3.4.15 \& 3.4.16]{Paternain2023} there exists a global chart $ {M}\rightarrow{\mathbb{D}}$ such that in these coordinates the metric is of the isothermal form
$\tilde \rho^2g=e^{2{\lambda}} \left(
\d x^2+\d y^2 \right)$,
where
$\mathbb{D}=\{z\in\mathbb{C}:|z|\leq 1\}$ and
${\lambda}\in C^\infty({\mathbb{D}})$.
Let $\rho:=\tilde \rho e^{-\lambda } : \mathbb{D} \rightarrow \mathbb{R}$, so that the metric is written in terms of global isothermal coordinates on $\mathbb{D}$ as 
\begin{equation}
	g=\rho^{-2} \left(\d x^2+\d y^2 \right)=\rho^{-2}|\d z|^{2},
\end{equation}
where we write $z=x+iy$.
In the following discussion we identify $M$ with $\mathbb{D}$ by a global chart as above.
The isothermal coordinates induce a smooth global parametrization of $\OSD$ as
\begin{equation}\label{eq:global_coords}
	M\times (\R/2\pi \mathbb{Z})\ni \left(x,y,\theta \right)\longmapsto (x,y,\rho \left(\cos{\theta}\sfrac{\partial }{\partial  x}+\sin{\theta}\sfrac{\partial }{\partial  y} \right)).	
\end{equation}
In terms of these coordinates,\footnote{The variable $\theta$ is different from the one in \eqref{eq:U2hat}. It will be clear in each discussion which one we mean. } the three canonical vector fields $X$, $X_\perp$, $V$ are written as 
\begin{equation}
    \begin{aligned}
        &\quad X = \cos(\theta)\rho\partial_{x}+\sin(\theta)\rho\partial_{y}+(\partial_{x}\rho\sin(\theta)-\partial_{y}\rho \cos(\theta))\partial_\theta, \label{eq:Xconf}\\
        X_\perp &= \sin(\theta)\rho\partial_{x}-\cos(\theta)\rho\partial_{y}-(\partial_{x}\rho \cos(\theta)+\partial_{y}\rho \sin(\theta))\partial_\theta,\quad \text{and}\quad
        V=\partial_\theta,
    \end{aligned}    
\end{equation}
see \cite[Lemma 3.5.6]{Paternain2023}.
It is easy to see using \eqref{eq:Xconf} that $\Omega_k$ in \eqref{eq:FourierDecomp} is globally spanned over $C^{\infty}(M^{\circ};\mathbb{C}^{d})$ by
$e^{ik\theta}$.
Moreover, the $k$-th mode is given by
\begin{equation}\label{eq:FourierMode}
u_k(z)=\frac{1}{2\pi}\int_{0}^{2\pi}{u\left(z,\theta^\prime\right)e^{-ik\theta^\prime}\d\theta^\prime}\, e^{ik\theta}.
\end{equation}
This implies that if
$u\in C^\infty(\OSM;\mathbb{C}^d)$, then $u_k$ extends to a smooth $\mathbb C^{d}$-valued function on $\OSM$.

The statement of Theorem~\ref{HolomIntFactExistThrm} included the \textbf{Guillemin-Kazhdan operators} $\eta_\pm\coloneqq \frac12(X\pm iX_\perp)$, which  can be written as
\begin{equation}
    \eta_+= e^{i\theta}(\rho \partial_{{z}}-i\partial_z\rho\, \partial _\theta), \qquad 
    \eta_-= e^{-i\theta}(\rho\partial_{\bar{z}}+i\partial_{\bar{z}}\rho\, \partial _\theta).
\end{equation}
Therefore, they act on the $k$-eigenspace of $-iV$ as 
\begin{equation}\label{eq:GK_isothermal}
	\begin{aligned}
\eta_+(f(z)e^{ik\theta})=	e^{i(k+1)\theta}(\rho \partial_{{z}}+k\partial_{{z}}\rho)f=	e^{i(k+1)\theta}\rho^{-k+1}\partial_{{z}}(\rho^{k}f),\\%
    \eta_-(f(z)e^{ik\theta})=e^{i(k-1)\theta}(\rho \partial_{\bar{z}}-k\partial_{\bar{z}}\rho)f=e^{i(k-1)\theta}\rho^{k+1}\partial_{\bar{z}}(\rho^{-k}f),%
\end{aligned}	 
\end{equation}  
where $\partial_z\coloneq(1/2)(\partial_x-i\partial_y)$ and $\partial_{\bar{z}}\coloneq(1/2)(\partial_x+i\partial_y)$.
From here
one sees the mapping property $\eta_{\pm}\colon \Omega_{k}\to \Omega_{k\pm 1}$.

\begin{remark}\label{rmk:edge}\upshape
	 	 The 0-sphere bundle is in a natural way an edge space in the sense of \cite{Mazzeo1991}, see Appendix \ref{sub:edge_structures}, where the role of the fibration map is played by the projection $\pi_0:\OSM\to M$.
	 	 It follows from the coordinate expressions \eqref{eq:Xconf} that $X$, $X_\perp$, $V$ are tangent to the fibers of $\pi_0$ at $\partial \OSM$ (in other words, they are vertical for $\rho=0$), so they are edge vector fields. 
	 	 As edge vector fields they are also linearly independent for every $ (z,v)\in \OSM$, so they furnish a global frame for the edge tangent bundle ${}^{e}T\OSM$.
	 	 Similarly, $\eta_\pm\in \mathrm{Diff}_e^{1}(\OSM)$.
	 \end{remark}

It is well known that one-forms on $M^{\circ}$ can be viewed as functions on $SM^{\circ}$, by means of the map 
\begin{equation}
	\ell_1 :C^\infty (M^{\circ};T^*M^{\circ})\to C^\infty(SM^{\circ}), \qquad \omega\mapsto \ell_1\omega, 
\end{equation}
where $(\ell_1 \omega)(z,v)=\omega_z(v)$.
In isothermal coordinates, if $\omega = a_1(z) \frac{\d z}\rho+ a_{-1}(z)\frac{\d\bar{z}}\rho$ we have
\begin{equation}\label{eq:correspondence}
	\ell_1\omega (z,\theta)= a_1(z) e^{i\theta}+ a_{-1}(z)e^{-i\theta}\in \Omega_{-1}\oplus \Omega_1,
\end{equation}
where $\Omega_{\pm1}$, see \eqref{eq:FourierDecomp}, are here for $d=1$. 
This shows that $\ell_1\colon C^\infty (M^{\circ};T^*M^{\circ})\to \Omega_{-1}\oplus \Omega_1$ is an isomorphism, and also that
 it naturally extends to a map $C^\infty (M;{}^{0}T^*M)\to C^\infty(\OSM)$. 
For more details and a generalization to symmetric tensor fields see \cite[Section 3.1.3]{Eptaminitakis2025}.

\section{The Non-Abelian X-ray Transform}\label{sec:the_non_abelian_x_ray_transform}

In this section, we prove  Theorems \ref{thrm:GaugeExist} and \ref{MainTranspEqFourierThrm} assuming Theorem~\ref{HolomIntFactExistThrm}.

\subsection{Regularity of Solutions to Transport Equations}\label{sec:RegOfSol}

In this section we address the regularity of solutions to certain transport equations  over AH spaces.
First we define our regularity spaces:
\begin{definition}\label{def:Rkr}
Suppose that $ \left(M^\circ,g \right)$ is an AH manifold, $k\geq0$ is an integer, and that $r\in\mathbb{R}$.
We define the space $\mathcal{R}^{k,r} \left(SM^\circ;\mathbb{C}^d \right)$ as the set of all sections $u\in C^\infty \left(SM^\circ;\mathbb{C}^d \right)$ such that for any smooth vector fields $V_1,\ldots,V_k\in C^\infty \left(\OSM;T\OSM \right)$ over the 0-sphere bundle,
\begin{equation}\label{eq:TestIf_uInRkr}
u,\quad V_1u,\quad \ldots\quad V_1\ldots V_ku\in\rho^r L^\infty \left(SM^\circ \right).
\end{equation}
\end{definition}
\begin{remark}\label{rmrk:Rkr}\upshape
The above is a slight generalization of \cite[Definition 5.1]{Grebnev2026} where the spaces $\mathcal{R}^k$ defined therein are equal to $\mathcal{R}^{k,\sfrac{ \left(n+1 \right)}{2}}$ in our notation.
An inessential difference between \cite[Definition 5.1]{Grebnev2026} and our definition is that the former uses ``${}^0\!S^\ast M$" and we use ``${}^0\!SM$" -- the definitions are equivalent since ${}^0\!S^\ast M\cong{}^0\!SM$ are diffeomorphic (over coordinates $(z^i)$ of $M$ the diffeomorphism is given by $\alpha_\mu(\sfrac{\d z^\mu}{\rho})\mapsto\rho^{-2}g^{\mu\nu}\alpha_\mu(\rho\sfrac{\partial}{\partial_{z^\nu}})$ -- see \eqref{eq:0SM})
\end{remark}
The following result follows from a careful examination of the proof of Proposition 5.56 in \cite{Grebnev2026} and we explain the required changes for the convenience of the reader.
The stronger conclusion below was not proved in \cite{Grebnev2026} since it was not needed at the time.
\begin{proposition}\label{prop:TranspEqReg} 

Suppose that $ \left(M^\circ,g \right)$ is a nontrapping AH manifold of dimension $n+1$ and that $\nabla$ is a smooth connection in $M\times\mathbb{C}^d$.
Suppose also that $\Phi\in\rho C^\infty \left(M;\mathbb{C}^{d\times d} \right)$, $\phi\in\rho C^\infty \left(M;\mathbb{C}^d \right)$, and $A\in C^\infty \left(M;\mathrm{Hom} \left(TM;\mathbb{C}^d \right) \right)$.
Then for any $h\in C^\infty\left(\partial_-\bSM;\mathbb{C}^d\right)$ there exists a unique solution $u\in C^\infty \left(\bSM;\mathbb{C}^d \right)$ to
$$\nabla_Xu+\Phi u=\phi+A$$
satisfying $u=h$ on $\partial _-\bSM$.

Suppose furthermore that $h=0$ and $u=0$ on $\partial _+\bSM$ and fix any integer $l\geq0$ and any $r\in\mathbb{R}$.
Then for all sufficiently large $N\geq0$, if $\phi\in\rho^{N+1} C^\infty \left(...\right)$ and $A\in\rho^N C^\infty \left(...\right)$  then $u\in\mathcal{R}^{l,r} \left(SM^\circ;\mathbb{C}^d \right)$.
The minimum size of $N$ only depends on $ \left(M^\circ,g \right)$, $M$, $l$, and $r$.

\end{proposition}

\begin{proof}

The proof is a minor modification of the proof of \cite[Proposition 5.56]{Grebnev2026}.
To explain the notation difference, \cite{Grebnev2026} writes $\mathcal{E},\End{\mathcal{E}},\nabla^\mathcal{E},\mathbb{X},M,\overline{M}$ for what in our case/notation is denoted by $M\times\mathbb{C}^d,M\times\mathbb{C}^{d\times d},\nabla,\nabla_X,M^\circ,M$ respectively.
The existence and uniqueness of the solution $u$ directly follows from \cite[Proposition 5.56]{Grebnev2026}.
We explain the changes required for the second statement.
Suppose that $h=0$ and $u=0$ on $\partial _+\bSM$ and fix an integer $l\geq0$ and $r\in\mathbb{R}$.
The first modification is changing ``$\mathcal{R}^k$" to ``$\mathcal{R}^{k,r}$," which requires replacing ``$\sfrac{ \left(n+1 \right)}{2}"$ with ``$r$" every time it appears in the proof.
Related to this, one needs to change all instances of ${}^0\!S^*M$ to $\OSM$, which does not affect the proof since the two manifolds are diffeomorphic (see Remark~\ref{rmrk:Rkr}).
To describe the second modification, we observe that we satisfy all of the assumptions of \cite[Proposition 5.56]{Grebnev2026} except the hypotheses that $\Phi\in\rho^{N+1}C^\infty \left(M;\mathbb{C}^{d\times d} \right)$ and that the connection coefficients $\Gamma_{\mu j}^k$ of $\nabla$
are in $\rho^NC^\infty \left(M \right)$
with respect to frames that cover the boundary.
These regularity assumptions play a role in only one place of the proof there, as we explain below.

Letting $\varphi$ denote the flow of $X$, the proof in \cite{Grebnev2026} proceeds by defining
$$B_j^k=\Gamma_{\mu j}^kv^\mu+\Phi_j^k\quad\mathrm{and}\quad b^k=\phi^k+A_\mu^kv^\mu\quad\mathrm{for}\quad (z,v)\in SM^\circ$$
and showing that since $u$ satisfies the following ODE on $SM^\circ$:
$$\frac{\d (u^k\circ \varphi_t)}{\d t}+(B^k_j\circ \varphi_t)(u^j\circ \varphi_t)=b^k\circ \varphi_t,\qquad u\circ \varphi_t(z,v)\overset {t\to -\infty}\to 0 \quad \forall (z,v)\in SM^{\circ },$$
then because $|b|\leq C \rho^N$ one has
\begin{equation}\label{eq:Grebnev2023uEst}
\left|u\right|\le C_0^\prime\rho^{N_0}\quad\text{on\ } SM^{\circ}
\end{equation}
for constants $C_0^\prime\geq0$ and $N_0\geq r$ dependent on $N$.
Furthermore,
$N_0$ can be made arbitrarily large by increasing $N$.

To analyze first order derivatives of $u$, if $V_1\in C^\infty \left({}^0\!SM;T{}^0\!SM \right)$, 
\begin{equation}
\begin{aligned}
\frac{\d (V_1(u^k\circ \varphi_t))}{\d t }+(B^k_j\circ \varphi_{t})V_1(u^j\circ \varphi_t)=&V_1(b^k\circ \varphi_{t})-V_1(B^k_j\circ \varphi_{t})(u^j\circ \varphi_{t}),\label{eq:V_1u}\\
V_1\big(u^j\circ \varphi_t(z,v)\big)\overset {t\to -\infty}\to&  0 \quad \forall (z,v)\in SM^{\circ }.
	\end{aligned}
\end{equation}
In \cite{Grebnev2026}, the stronger decay assumptions on $B^k_j$ are used to obtain that the right-hand side of this ODE is bounded by $\tilde{C}\rho^{\tilde{N}}$, where $\tilde{N}\geq0$ can be made arbitrarily large. This leads to
$$\left|u\right|\le C_1^\prime\rho^{N_1}\quad\text{on\ } SM^{\circ}$$
for $C_1^\prime\geq0,N_1\geq r$ dependent on $N$, where $N_1$ can be made arbitrarily large by increasing $N$.
The difference in our case, since we have a weaker assumption on the decay of $B^k_j$, is to instead use the assumption $|b|\leq C \rho^N$ and \eqref{eq:Grebnev2023uEst} to analogously bound the right-hand side of \eqref{eq:V_1u} by $\tilde{C}\rho^{\tilde{N}}$.
The case of higher order derivatives of $u$ is similarly dealt with inductively.
\end{proof}

\subsection{Fourier Analysis of the Transport Equation -- Proof of Theorem~\ref{MainTranspEqFourierThrm}}\label{sec:FourAnalyTranspEq}

In this section we prove Theorem~\ref{MainTranspEqFourierThrm} on the Fourier content of solutions to transport equations assuming Theorem~\ref{HolomIntFactExistThrm}.
We assume throughout that $M$ is two-dimensional (i.e. $n=1$). 
We will use the vector fields \eqref{eq:Xconf} extensively here.
We begin with a few definitions and properties.
\begin{definition}\label{HilbertTransfDef}

Suppose that $f\in C^\infty \left(SM^{\circ } \right)$. 
Decomposing $f \left(z,\theta \right)=\sum_{k=-\infty}^{\infty}{f_k}$, $f_k\in \Omega_k$, the \textbf{fiberwise Hilbert transform} is defined as
\begin{equation}
	Hf (z,\theta )=\sum_{k=-\infty}^{\infty} (-i)
	\sgn(k)f_k (z),
	 \quad \text{where }\quad \sgn{k}=\begin{cases}
	 \  \pm  1, & \mathrm{if}\ \pm k>0,\\ \ \ 0, & \mathrm{if}\  k=0.
	 \end{cases}
\end{equation}
\end{definition}
\noindent One can readily check that $H(\overline{f})=\overline{Hf}$ and hence $H$ is a real operator. The Hilbert transform acts boundedly on edge spaces (see \S\ref{sub:edge_structures}):
\begin{proposition}\label{prop:HilbMapProp}
The Hilbert transform is  continuous  $ \rho^\delta H_e^k \left(\OSM;\mathbb{C}^d \right)\rightarrow\rho^\delta H_e^k \left(\OSM;\mathbb{C}^d \right)$.
\end{proposition}

\begin{proof}
We will prove the proposition for $k=0,1$ since the proof for larger $k$ is similar.
Consider the global coordinates of $\OSM$ given by the inverse of \eqref{eq:global_coords}, over which $\left\{\rho\partial _x,\rho\partial _y,\partial _\theta\right\}$ is a frame for $\mathcal{V}_e \left(\OSM \right)$.
Fix $u\in C_c^\infty \left(SM^\circ;\mathbb{C}^d \right)$ and $W$ equal to $\rho\partial _x$, $\rho\partial _y$, or $\partial _\theta$.
Since $C_c^\infty \left(SM ^{\circ}\right)$ is dense in $H_e^k \left(\OSM \right)$, the proposition for $k=0,1$ will follow if we prove that $\left\lVert Hu\right\rVert _{\rho^\delta L^2}\le\left\lVert u\right\rVert _{\rho^\delta L^2}$ and $\left\lVert WHu\right\rVert _{\rho^\delta L^2}\le\left\lVert Wu\right\rVert _{\rho^\delta L^2}$.
We have that
$$
\begin{aligned}
		\left\lVert Hu\right\rVert _{\rho^\delta L^2}^{2}=&\int_{M^\circ}{\sum_{k=-\infty}^{\infty}\lVert  \left(-i \right)\sgn{ \left(k \right)}u_k\rVert _{\rho^\delta L^2 \left(S_zM^\circ \right)}^{2}\d V_g \left(z \right)}\\
		&\qquad \le\int_{M^\circ}{\sum_{k=-\infty}^{\infty}\left\lVert u_k\right\rVert^{2} _{\rho^\delta L^2 \left(S_zM^\circ \right)}\d V_g \left(z \right)}=\left\lVert u\right\rVert _{\rho^\delta L^2}^{2}.
\end{aligned}$$
Next, writing $u_k={\widehat{u}}_k \left(z \right)e^{ik\theta}$ observe that $\rho\partial _i{\widehat{u}}_k={\widehat{ \left(\rho\partial _iu \right)}}_k$ for $i=x,y$.
Hence a similar calculation gives that $\left\lVert \partial _\theta H u\right\rVert _{\rho^\delta L^2}=\left\lVert \partial _\theta u\right\rVert _{\rho^\delta L^2}$ and $\left\lVert \rho\partial _iHu\right\rVert _{\rho^\delta L^2}{\leq}\left\lVert \rho\partial _iu\right\rVert _{\rho^\delta L^2}$ for $i=x,y$.
\end{proof}
\noindent Since $(\id+H^2)u=u_0$, one also gets
\begin{corollary}\label{cor:u_to_u0ContH0Map}
The map $u\mapsto u_0$ is continuous $\rho^\delta H_e^k \left(\OSM;\mathbb{C}^d \right)\rightarrow\rho^\delta H_e^k \left(\OSM;\mathbb{C}^d \right)$.
\end{corollary}
\begin{definition}\label{def:ConnCurvHodg}
Suppose $\nabla$ is a connection and let $B$ be its connection matrix.
The \textbf{curvature} of $\nabla$ is the tensor field $F_B
\in C^\infty
\left(M;\Lambda^2 \left(T^\ast M \right)\otimes\mathbb{C}^{d\times d} \right)$ defined by
\begin{equation}\label{ConnCurvDefEq}F_B=\d B+B\land B.
\end{equation}
Its \textbf{Hodge star} $\star F_B\in C^\infty \left(M^\circ;\mathbb{C}^{d\times d} \right)$ is defined as
\begin{equation}
\label{HodgeStarConnCurvDefEq}\star F_B \big(z \big)=F_B \big(v,v^\bot \big)\end{equation}
where $v\in S_zM^\circ$ is arbitrary (``$v^\perp$" defined in \S\ref{sec:vertical_fourier}).
In addition, we define the \textbf{Hodge star} of $B$ over $SM^\circ$ as
\begin{equation}\label{eq:HodgeStarConnMtrxDef}
\star B (z,v )=-V({\ell_1}B)(z,v)=-B \big(z,v^\bot \big).
\end{equation}
\end{definition}
Using \eqref{eq:global_coords} and \eqref{HodgeStarConnCurvDefEq}, a computation in the isothermal coordinates of \S\ref{sec:vertical_fourier} implies that for a $C^\infty$ unitary connection $\nabla$ on $M\times\mathbb{C}^d$,
\begin{equation}\label{eq:starFB}
	\star F_B\in \rho ^{2}C^{\infty} \big( M ;\End_{\mathrm{Sk}}{\mathbb{C}^d} \big).
\end{equation}
\begin{lemma}\label{lm:P_function}
There exists a real valued $f\in \Aphg^{E}(M)$ with $E=\overline{\{(0,1)\}}$  such that for $\varphi \coloneqq-i X_\perp f=  -\eta_+ f+\eta _- f$ one has  $i\,  \d(\ell_{1}^{-1}\varphi)=P\left(z\right)\d V_g$, where $P\in C^\infty\left(M^{\circ};\mathbb{R}\right)$ extends to a smooth, positive function on $M$.
Furthermore, $\varphi\in C^\infty(\OSM)$.
\end{lemma}
\begin{proof}
Identify $M $ with $\mathbb{D}$  using global isothermal coordinates \eqref{eq:global_coords} and define $f\in C^\infty( M^{\circ})$ by
\begin{equation}
	f=-\log (1-x^2-y^2)=-\log (1-z\bar z).
\end{equation}
Since $1-z\bar{z}$ is a smooth bdf for $\partial \mathbb D$, $\rho =e^{h }(1-z\bar{z})$ for some $h\in C^\infty(M)$; this shows that $ f\in \Aphg^E(M)$ with $E=\overline{\{(0,1)\}}$.
Then by \eqref{eq:GK_isothermal},
\begin{equation}
\varphi=
-\rho\partial _z f\, e^{i\theta}+\rho\partial _{\overline{z}} f\,e^{-i\theta}\overset{\eqref{eq:correspondence}}\implies \ell_1^{-1}\varphi =  -\partial_z f\, \d z+\partial_ {\bar{z}}f\,\d \bar z.
\end{equation}
We immediately obtain $ \varphi =2i e^{h}\Im (z e^{-i\theta})\in C^\infty(\OSM)$ and
\begin{equation}
	i \, \d (\ell_1^{-1}\varphi) = 2i \,\partial_{z\bar z}^{2}f\,\d z\wedge \d \bar z=\frac {4\d x\wedge\d y}{(1-x^2-y^2)^{2}}
	=4e^{2h}\frac {\d x\wedge\d y}{\rho^{2}}
	=4 e^{2h}\d V_g,
\end{equation}
which yields the first claim with $P=4 e^{2h}$.
\end{proof}
\begin{proof}[Proof of Theorem~\ref{MainTranspEqFourierThrm}]

In this proof we let $L^p=L^p(SM^\circ)$.
We follow closely the proof of Theorem 13.8.1 in \cite{Paternain2023}.
We will prove that $u$ is holomorphic.
Applying the same argument to $\bar u$ shows that $u$ is also anti-holomorphic, from which the theorem follows.
Throughout, we will identify $S^*M^\circ$ with $SM^\circ$ via the musical isomorphisms without denoting them explicitly.
First we observe that $(\phi+B)\in\Omega_{-1}\oplus\Omega_0\oplus\Omega_1$, which is most easily seen in the coordinates \eqref{eq:global_coords}.
Next, let us write $\nabla_v\psi=v \left(\psi \right)+A(v)\psi$. 
Let $f$, $\varphi$ and $P\in C^\infty(M)$  be as in Lemma~\ref{lm:P_function} and let $w:SM^\circ\rightarrow\mathbb{C}$ be a holomorphic function as in Theorem~\ref{HolomIntFactExistThrm} with respect to these $f$ and $\varphi$.
For every $s\in\mathbb{R}$, consider the unitary connection $\nabla_v^s\psi=v(\psi)+A_s(v)\psi$ where $A_s:=A-s\varphi \id$
and define $u_s:=e^{sw}u$.
Note that $u$ is holomorphic if and only if $u_s$ is.
 Routine calculations give that $u_s$ satisfies the new transport equation $\nabla_X^su_s+\Phi u_s=e^{sw}(\phi+B)$ 
 and that
\begin{equation}\label{eq:FAsFARelat}
F_{A_s}=F_A+is(P(z)\d V_g)\otimes\textrm{id}\quad\Longrightarrow\quad \star F_{A_s}=\star F_A+isP(z)(\textrm{id}).
\end{equation}

We claim that $u_s\in H_e^2(\OSM;\mathbb{C}^d)$ for all $s$.
Fix any $s$.
The claim will follow if for any $W_1,W_2\in\mathcal{V}_e(\OSM)$, the functions ${u_s},W_1{u_s}$, $W_1W_2{u_s}$ are all in $L^2$.
By Theorem~\ref{HolomIntFactExistThrm}, the functions $e^{sw}$, $W_1e^{sw}$, $W_1W_2e^{sw}$ are in $\rho^{-m_s}L^2(\OSM)$ for some $m_s\in\mathbb{R}$.
By the rapid decay assumption of
$(\phi+B)$
, Proposition~\ref{prop:TranspEqReg} implies that $u\in\mathcal{R}^{N,\sfrac{ \left(n+1 \right)}{2}+m_s}$ for any $N\geq0$ and hence indeed $u_s\in H_e^2(\OSM;\mathbb{C}^d)$.

Then, to show that $u_s$ is holomorphic, define
$${\widetilde{u}}_s:=\sum_{k=-\infty}^{-1} \left(u_s \right)_k=\frac{1}{2} \left(u_s-iH \left(u_s \right)- \left(u_s \right)_0 \right).$$
We will show that $\widetilde u_s$ vanishes.
Note that the operator $\psi\mapsto{\widetilde{\psi}}$ is a continuous $H_e^2 \left(\OSM;\mathbb{C}^d \right)\rightarrow H_e^2 \left(\OSM;\mathbb{C}^d \right)$ map by Proposition~\ref{prop:HilbMapProp} and Corollary \ref{cor:u_to_u0ContH0Map}.
Next, by \cite[Lemma 13.4.3]{Paternain2023} we have that $\nabla_X=X+A : \Omega_m\rightarrow\Omega_{m-1}\oplus\Omega_{m+1}$ for any $C^\infty$ connection $\nabla$.
Similarly, $u\mapsto\Phi u$ clearly maps $\Omega_m\rightarrow\Omega_m$.
Hence $\nabla_{X}^{s}{\widetilde{u}}_{s}+\Phi {\widetilde{u}}_{s}\in \oplus_{k=-\infty }^{0} {\Omega }_{k}$ and
$$\nabla_{X}^{s}{\widetilde{u}}_{s}+\Phi {\widetilde{u}}_{s}=\nabla_{X}^{s}u_s+\Phi u_s+\nabla_{X}^{s}\left({\widetilde{u}}_{s}-u_s\right)+\Phi \left({\widetilde{u}}_{s}-u_s\right)\in \left(e^{sw}{(\phi+B)}+\oplus_{k=-1}^{\infty } {\Omega }_{k}\right)\subseteq\oplus_{k=-1}^{\infty } {\Omega }_{k}.$$
So
\begin{equation}\label{eq:uTildeTrasnpEq}
\nabla_X^s\widetilde{u}_s+\Phi \widetilde{u}_s\in\Omega_{-1}\oplus\Omega_0.
\end{equation}

Now, we claim that
\begin{equation}
\begin{aligned}
\label{eq:PestovProofV}
\left\lVert  \left(X+A_s+\Phi \right) \left(V{\widetilde{u}}_s \right)\right\rVert _{L^2}^2-\langle KV{\widetilde{u}}_s,V{\widetilde{u}}_s\rangle _{L^2}+\left\lVert  \left(X+A_s+\Phi \right){\widetilde{u}}_s\right\rVert _{L^2}^2-\left\lVert V \left(X+A_s+\Phi \right){\widetilde{u}}_s\right\rVert _{L^2}^2\quad \\
-\langle \star F_{A_s}{\widetilde{u}}_s,V{\widetilde{u}}_s\rangle_{L^2} -\operatorname{Re}{\langle  \left(\star \d_{A_s}\Phi \right){\widetilde{u}}_s,V{\widetilde{u}}_s\rangle_{L^2} }-\operatorname{Re}{\langle \Phi{\widetilde{u}}_s, \left(X+A_s+\Phi \right){\widetilde{u}}_s\rangle_{L^2} }=0,
\end{aligned}
\end{equation}
where $\d_{A_s}\Phi:=\d\Phi+[A_s,\Phi]$ and $K$ is the Gauss curvature.
Eq. \eqref{eq:PestovProofV} is known as a \textbf{Pestov Identity}.
Note that both $\star F_{A_s}$ and $\star \d_{A_s}\Phi$ are smooth over $\OSM$ -- see comment after \eqref{eq:0SM}.
We remind that in \eqref{eq:PestovProofV}, $K$ is bounded (see \S\ref{sec:Setup}) and $X,V$ are edge operators and hence map $H_e^m\rightarrow H_e^{m-1}$.
Furthermore, by the construction of $\varphi$ in Lemma~\ref{lm:P_function}, $A_s,V(A_s) \in C^\infty(\OSM)$.
Thus all operators inside the $L^2$ norms in \eqref{eq:PestovProofV} (e.g. ``$\left(X+A_s+\Phi \right)V$" in the first norm) are continuous $H_e^2 \left(\OSM;\mathbb{C}^d \right)\rightarrow H_e^0 \left(\OSM;\mathbb{C}^d \right)$ maps.

The equality \eqref{eq:PestovProofV} for a compact surface with boundary $N$ can be found in \cite[Eq. (13.20)]{Paternain2023} 
(the authors write ``$v$" instead of ``${\widetilde{u}}_s$") 
when ${u}_s\in C^\infty(SN)$ and vanishes on the boundary and $A_s$ and $\Phi$ are skew-Hermitian.
Hence \eqref{eq:PestovProofV} holds in our context for any $ u_s\in C_c^\infty \left(SM^\circ;\mathbb{C}^d \right)$ by applying their result to a compact 2 dimensional submanifold with boundary of $(M^\circ,g)$ containing $\supp{u_s}$.
Thus, \eqref{eq:PestovProofV} holds for $ u_s\in H_e^2 \left(\OSM;\mathbb{C}^d \right)$ by density.

Defining $h_s:=\left[ \left(X+A_s+\Phi \right){\widetilde{u}}_s\right]_0$, the following hold on compact surfaces with boundary by (13.21), (13.22), (13.24), and (13.25) in \cite{Paternain2023}:
\begin{align}
\left\lVert  \left(X+A_s+\Phi \right) \left(V{\widetilde{u}}_s \right)\right\rVert _{L^2}^2-\langle KV{\widetilde{u}}_s,V{\widetilde{u}}_s\rangle _{L^2}\geq0, \label{eqs:FourierInd} \\
\left\lVert  \left(X+A_s+\Phi \right){\widetilde{u}}_s\right\rVert _{L^2}^2-\left\lVert V \left(X+A_s+\Phi \right){\widetilde{u}}_s\right\rVert _{L^2}^2=\left\lVert h_s\right\rVert _{L^2}^2, \label{eqs:Fourier0Mode}\\
\operatorname{Re}{\langle  \left(\star \d_A\Phi \right){\widetilde{u}}_s,V{\widetilde{u}}_s\rangle _{L^2}}\le C_{A,\Phi}\sum_{k=-\infty}^{-1}{\left|k\right|\left\lVert  \left({\widetilde{u}}_s \right)_k\right\rVert _{L^2}^2}, \label{eqs:FourierdAsPhi}\\
\operatorname{Re}{\langle \Phi{\widetilde{u}}_s, \left(X+A_s+\Phi \right){\widetilde{u}}_s\rangle _{L^2}}\le C_{A,\Phi}\sum_{k=-\infty}^{-1}\left\lVert  \left({\widetilde{u}}_s \right)_k\right\rVert _{L^2}^2,\label{eqs:FourierPhi}
\end{align}
where
\eqref{eqs:FourierInd} requires the surface to be simple (in our case, this holds for the sets $\{\rho\geq \epsilon\}$ for $0<\epsilon \ll 1$) and 
\eqref{eqs:Fourier0Mode} requires $\widetilde{u}_s$ to satisfy a transport equation with right-hand side in $\Omega_{-1}\oplus\Omega_0\oplus\Omega_1$, see \eqref{eq:uTildeTrasnpEq}.
A careful look at their proof reveals that the constant $C_{A,\Phi}\geq0$ only depends on the $L^\infty(SM^\circ)$ norm of $\eta_-\Phi$, $\Phi$, and $A_{\pm1}, b_{\pm1}\in\Omega_{\pm1}$ in the
decompositions $A=A_1+A_{-1}$ and
$\d_A\Phi=b_1+b_{-1}$, all of whose norms are finite by \eqref{eq:global_coords}, \eqref{eq:GK_isothermal}, and \eqref{eq:FourierMode}.
In particular $C_{A,\Phi}$ does not depend on $s$ nor on $\supp u_s$.
A similar density argument as above implies that these also hold for $ u_s\in H_e^2 \left(\OSM;\mathbb{C}^d \right)$ where one justifies interchanging the limit and the above ``$\sum_{k=-\infty}^{-1} \ldots$ " sums by the fact that for some constant $C>0$ the sums in the limit are bounded by $C \sum_{k=-\infty}^{-1}\left\|k\left(\widetilde{u}_{s}\right)\right\|_{L^{2}}^{2}=\left\|V \widetilde{u}_{s}\right\|_{L^{2}}^{2}<\infty$. 

We need one more inequality.
Let $\psi \in C_{c}^{\infty}\left(S M^{\circ};\mathbb C^d\right)$.
By \eqref{eq:FAsFARelat}
\begin{equation}
\begin{aligned}
		-\big\langle\star F_{A_{s}} \widetilde{\psi}, V \widetilde{\psi}\big\rangle_{L^{2}}=&-\sum_{k=-\infty}^{-1} \int_{M^{\circ}}\big\langle\left(\star F_{A}+i s P(z)\right) \widetilde{\psi}_{k}, i k \widetilde{\psi}_{k}\big\rangle_{L^{2}\left(S_{z} M^{\circ}\right)} \d V_{g}	 \\
		=&\sum_{k=-\infty}^{-1}|k| \int_{M^{\circ}}\big\langle\left(s P(z)-i \star F_{A}\right) \widetilde{\psi}_{k}, \widetilde{\psi}_{k}\big\rangle_{L^{2}\left(S_{z} M^{\circ}\right)} \d V_{g}\\
		\geq&\Big(s\inf_{M^\circ}{\left\{P\right\}}-\left\lVert\star F_A\right\rVert_{L^\infty}\Big)\sum_{k=-\infty}^{-1}{\left|k\right|\big\lVert{\widetilde{\psi}}_k\big\rVert_{L^2}^2}.
\end{aligned}
\label{eqs:FourierFAs}
\end{equation}

Finally, plugging \eqref{eqs:FourierInd}, \eqref{eqs:Fourier0Mode}, \eqref{eqs:FourierdAsPhi}, \eqref{eqs:FourierPhi}, and \eqref{eqs:FourierFAs} (with $w=u_{s}$) into \eqref{eq:PestovProofV} and observing that $\d_{A_{s}} \Phi=\d_{A}\Phi$ (since $\left[A_{s}, \Phi\right]=[A, \Phi]$) gives 
$$0 \geq\|h_s\|_{L^{2}}^{2}+\Big(s\inf_{M^\circ}{\left\{P\right\}}-\left\lVert\star F_A\right\rVert_{L^\infty}-2C_{A,\Phi}\Big)\sum_{k=-\infty}^{-1}{\left|k\right|\left\lVert{\widetilde{u}}_k\right\rVert_{L^2}^2}.$$
Hence choosing $s>0$ big enough implies that the right-hand side is positive and thus for $k\leq-1$, $ \left({\widetilde{u}}_s \right)_k= \left(u_s \right)_k=0$, implying that $u_s$ is fiberwise holomorphic.
So $u=e^{-sw}u_s$ is fiberwise holomorphic as well, as mentioned before.
\end{proof}

\subsection{Boundary Determination}\label{sec:ArbDecay}

The main result of this section (Theorem~\ref{thm:ArbDecay}) demonstrates that if two connection-Higgs field pairs $ \left(\nabla,\Phi \right)$ and $(\widetilde{\nabla},\widetilde{\Phi})$ have the same non-Abelian X-ray transform, then one can modify one of these pairs in a gauge equivalent way so that the resulting pairs agree at $\partial M$ up to infinite order.
We begin with a preliminary result:
\begin{proposition}\label{prop:ArbDecayRhoDir}
Suppose that $ \left(M^\circ,g \right)$ is an AH manifold and that $\nabla,\widetilde{\nabla}$ are $C^\infty$ connections in $M\times\mathbb{C}^d$.
Let $\rho$ be a boundary defining function (not necessarily geodesic) and consider the diffeomorphism $[0,\epsilon)_\rho\times\partial M\rightarrow\mathcal{C}_\epsilon$ induced by the flow of $\grad_{\rho^2g}{\rho}$ where $\mathcal{C}_\epsilon$ is a collar neighborhood of $\partial M$.
Then there exists a $Q\in C^\infty(M;\mathbb{C}^{d\times d})$ with $Q|_{\partial M}=\mathrm{id}$ such that the connection $\widetilde{\nabla}=Q^{-1}\nabla Q$ satisfies
$(\widetilde{\nabla}-\widehat{\nabla}) \left(\partial _\rho \right)=0$ in a collar neighborhood of $\partial M$.
\end{proposition}
\begin{proof}
Let $B,\widetilde{B},\widehat{B}$ be the connection matrices of $\nabla,\widetilde{\nabla},\widehat{\nabla}$ respectively (see \ref{eq:ConnExplForm}).
We let $B_0$ denote $B \left(\partial _\rho \right)$ and similarly for ${\widetilde{B}}_0,{\widehat{B}}_0$.
Let $Q$ be the $\mathbb{C}^{d\times d}$-valued solution to the following PDE on $\mathcal C_\epsilon$:
\begin{equation}\label{eq:DiffEqForQ}
\partial _\rho Q-Q{\widetilde{B}}_0+B_0Q=0,\qquad\left.Q\right|_{\partial  M}=\mathrm{id}.
\end{equation}
Define $\widehat{\nabla}=Q^{-1}\nabla Q$ on $\mathcal{C}^\prime$, which is equivalent to $\widehat{B} \left(v \right)=Q^{-1}v \left(Q \right)+Q^{-1}BQ$.
Using this and \eqref{eq:DiffEqForQ} we get that on $\mathcal C_\epsilon$,
$$ \big(\widetilde{\nabla}-\widehat{\nabla} \big) \left(\partial _\rho \right)={\widetilde{B}}_0-{\widehat{B}}_0=Q^{-1} \big(Q{\widetilde{B}}_0-\partial _\rho Q-B_0Q \big)=0.$$
Next, for every $0<\epsilon'<\epsilon$, let $\psi_{\epsilon^\prime} : M\rightarrow\left[0,1\right]$ be a $C^\infty$ function that is equal to $1$ on $\left\{\rho\le\sfrac{{\epsilon^\prime}}{2}\right\}$ and 0 on $\left\{\rho>0.9{\epsilon^\prime}\right\}$.
Define $Q_{\epsilon^\prime}=\mathrm{id}+\psi_{\epsilon^\prime} \left(Q-\mathrm{id} \right)$ which still satisfies \eqref{eq:DiffEqForQ} on $\{\rho<{\epsilon^\prime}\}$.
Since $\left.Q\right|_{\partial  M}=\mathrm{id}$ and $\partial M$ is compact, $Q_{\epsilon^\prime}$ has nonzero determinant everywhere for small enough ${\epsilon^\prime}>0$.
Fixing such an ${\epsilon^\prime}$, redefining $Q=Q_{\epsilon^\prime}$ gives us the desired $Q$.
\end{proof}
Proposition~\ref{prop:ArbDecayRhoDir} did not assume equality of scattering data. The latter will be needed for boundary determination. For that, as well as to reduce the proof of Theorem~\ref{thrm:GaugeExist} to Theorem~\ref{MainTranspEqFourierThrm}, we will make use of the following pseudolinearization identity. If $(M^\circ,g)$ is nontrapping AH and $\left(\Phi,\nabla \right)$,  $(\widetilde{\Phi},\widetilde{\nabla} )$ are pairs with $\nabla$, $\widetilde \nabla $ smooth and $\Phi,\widetilde{\Phi}\in\rho C^\infty(M;\mathbb{C}^{d\times d})$ such that $T^{ \left(\Phi,\nabla \right)}=T^{ (\widetilde{\Phi},\widetilde{\nabla} )}$, there exists $W\in C^\infty(\bSM;\mathbb C^{d\times d})$ such that
\begin{equation}\label{eq:TranspEqForW}
XW+ \left(B+\Phi \right)W-W \big(\widetilde{B}+\widetilde{\Phi} \big)=\widetilde{B}+\widetilde{\Phi}-B-\Phi,\qquad W=0\ \mathrm{on}\ \partial \bSM,
\end{equation}
where $B,\widetilde{B}$ denote the connection matrices of $\nabla,\widetilde{\nabla}$ respectively.
Indeed, let $U,\widetilde{U}\in C^\infty \big(\bSM;\mathbb{C}^{d\times d} \big)$ be the solutions to
$$XU+BU+\Phi U=0,\qquad X\widetilde{U}+\widetilde{B}\widetilde{U}+\widetilde{\Phi}\widetilde{U}=0,\qquad\left.U\right|_{\partial _-\bSM}=\widetilde{U}\big|_{\partial _-\bSM}=\mathrm{id},$$
whose existence is justified by applying Proposition~\ref{prop:TranspEqReg} to the columns of $U$ and $\widetilde{U}$ (e.g. write $\nabla_Xu=Xu+Bu$).
Observe that for any fixed $e\in\mathbb{C}^d$, $Ue$ satisfies \eqref{eq:TranspEqDefineT} with $A$ replaced by $B$ and similarly for $\widetilde{U}$ with $A,\Phi$ replaced by $\widetilde{B},\widetilde{\Phi}$ respectively.
By assumption, we also have that for any integral curve $\sigma :  \left(-\infty,\infty \right)\rightarrow S^\ast M^\circ$ of $X$ and $ \left(x_-,\zeta_- \right)=\lim_{t\rightarrow-\infty}{\sigma \left(t \right)}\in\partial_-{}^bS^*M$,
$$\lim_{t\rightarrow\infty}{ \left(U\circ\sigma \right)e}=T^{ (\Phi,\nabla )} \left(x_-,\zeta_-,e \right)=T^{ (\widetilde{\Phi},\widetilde{\nabla} )} \big(x_-,\zeta_-,e \big)=\lim_{t\rightarrow\infty}{ \big(\widetilde{U}\circ\sigma \big)e}.$$
Per the discussion in \S\ref{sec:NonAbelTrasfDef}, any $ \left(x_+,\zeta_+ \right)\in\partial _+\bSM$ is the limit of such a curve $\sigma \left(t \right)$ as $t\rightarrow\infty$.
This implies that $U=\widetilde{U}$ on $\partial _+\bSM$ as well, so setting  $W=U{\widetilde{U}}^{-1}-\mathrm{id}$ yields \eqref{eq:TranspEqForW}.

\begin{theorem}\label{thm:ArbDecay}

Assume the same hypotheses
and notations
as in Proposition~\ref{prop:ArbDecayRhoDir}
and suppose in addition that $(M^\circ,g)$ is nontrapping.
Suppose also that we have $\Phi,\widetilde{\Phi}\in\rho C^\infty(M;\mathbb{C}^{d\times d})$ and the equality of the non-Abelian X-ray transforms: $T^{ \left(\Phi,\nabla \right)}=T^{ (\widetilde{\Phi},\widetilde{\nabla} )}$.
Then one can modify $ \left(\Phi,\nabla \right)$ up to gauge so that for any integer $N\geq0$,
\begin{equation}\label{eq:ArbDecayConcl}
\big(\widetilde{\nabla}-\nabla \big)\in\rho^NC^\infty \big(M;\mathrm{Hom} \big(TM;\mathbb{C}^{d\times d} \big) \big)\qquad\mathrm{and}\qquad\big(\widetilde{\Phi}-\Phi \big)\in\rho^{N+1}C^\infty \big(M;\mathbb{C}^{d\times d} \big).
\end{equation}
\end{theorem}

\begin{proof}

Let $Q$ be as in Proposition~\ref{prop:ArbDecayRhoDir} and redefine $ \left(\Phi,\nabla \right)=(Q^{-1}\Phi Q,Q^{-1}\nabla Q)$.
By the discussion in \S\ref{sec:NonAbelTrasfDef}, the non-Abelian X-ray transform does not change under this gauge transformation and so we still have that $T^{ \left(\Phi,\nabla \right)}=T^{ (\widetilde{\Phi},\widetilde{\nabla} )}$.
We prove the theorem by demonstrating that the new $ \left(\Phi,\nabla \right)$ satisfies \eqref{eq:ArbDecayConcl}.
Fix coordinates $(y^i)$ of $\partial M$ which generate
boundary coordinates $ \left(x^\mu \right)= \left(\rho,y^i \right)$
on $\mathcal{C}_\epsilon$.
It will suffice to
show by induction on $l=0,1,2\ldots$ that
\begin{equation}\label{eq:ArbDecayIDir}
\big(\widetilde{\nabla}-\nabla \big) \big(\partial _{x^\mu} \big)\in\rho^lC^\infty \big(M;\mathbb{C}^{d\times d} \big)\ \mathrm{and}\ \big(\widetilde{\Phi}-\Phi \big)\in\rho^{l+1}C^\infty \big(M;\mathbb{C}^{d\times d} \big)\quad\mathrm{for}\ \mu=0,\ldots,n.
\end{equation}
The base case $l=0$ holds by
assumption.
Due to the properties of $Q$, we see that the new $\nabla$ satisfies $(\widetilde{\nabla}-\nabla) \left(\partial _\rho \right)=0$ (i.e. the case $\mu=0$).
Now suppose \eqref{eq:ArbDecayIDir} holds for some integer $l\geq0$; we will prove it for $l+1$.

Let $B,\widetilde{B}$ be the connection matrices of $\nabla,\widetilde{\nabla}$ respectively.
Recalling that $\overline{X}=\rho^{-1}X$, \eqref{eq:TranspEqForW}
can be rewritten as the following PDE on $\bSM$:
\begin{equation}\label{eq:TranspEqForWExplic}
\overline{X}W+\rho^{-1} \left(B+\Phi \right)W-W\rho^{-1} \big(\widetilde{B}+\widetilde{\Phi} \big)=\rho^{-1} \big(\widetilde{B}+\widetilde{\Phi}-B-\Phi \big),\quad W=0\ \mathrm{on} \ \partial\bSM.
\end{equation}
Since the restriction of vectors $v_\mu\partial_{x^\mu}\in TM$ to $\bSM$ satisfy $v^0=O(\rho)$ and $v^i=O(\rho^2)$ for $i\geq1$ (see \eqref{eq:vmu}), all of the coefficients in this PDE extend smoothly to $\partial {}^b\!S^* M$ and hence $W$ can be expressed as solutions to an ODE along integral curves of $\overline{X}$ with zero initial conditions. 

Letting $E_i^j\in\mathbb{C}^{d\times d}$ denote the $d\times d$ matrix with a one in the $i^{\mathrm{th}}$ row and $j^{\mathrm{th}}$ column and zeros otherwise, and letting $e_i\in\mathbb{R}^{ (d^2 )}$ denote the column vector with a one in the $i^{\mathrm{th}}$ row and zeros otherwise, consider the isomorphism $\mathcal{I} : \mathbb{C}^{d\times d}\rightarrow\mathbb{C}^{d^2}$ given by $\mathcal{I}(E_i^j)=e_{ \left(i-1 \right)d+j}$ (we also write $E_i^j\cong e_{ \left(i-1 \right)d+j}$).
Then there exists $\mathcal L\in C^{\infty}(\bSM;\mathbb{C}^{d^2\times d^2})$ such that \eqref{eq:TranspEqForWExplic} takes the standard form
\begin{equation}\label{eq:with_linear_isom}
	\overline X (\mathcal I W)+\mathcal L (\mathcal I W)=\rho^{-1}\mathcal I\big(\widetilde{B}+\widetilde{\Phi}-B-\Phi \big), \quad\mathcal I W=0\ \mathrm{on} \ \partial\bSM.
\end{equation}
For $i=1,\ldots,d^2$, let $\mathcal{G}_i\in C^\infty \left({}^b\!S^* M;\mathbb{C}^{d\times d} \right)$ be the solution to the homogeneous version of \eqref{eq:TranspEqForWExplic}:
\begin{equation}\label{eq:XBarTranspEqForG}
\overline{X}\mathcal{G}_i+\rho^{-1} \left(B+\Phi \right)\mathcal{G}_i-\mathcal{G}_i\rho^{-1} \big(\widetilde{B}+\widetilde{\Phi} \big)=0\quad\mathrm{on} \ M\quad\mathrm{satisfying}\quad\left.\mathcal{G}_i\right|_{\partial _-{}^b\!S^* M}\cong e_i,
\end{equation}
in particular 
\begin{equation}\label{eq:std_form_homog}
	\overline X \mathcal I (\mathcal G_i)+\mathcal L \mathcal I (\mathcal G_i) = 0.
\end{equation}
The vectors $\mathcal{G}_1,\ldots,\mathcal{G}_{d^2}$ are linearly independent since they are solutions to an ODE along integral curves of $\overline{X}$ with linearly independent initial conditions.
Let $\mathcal{G}\in C^\infty \big({}^b\!S^* M;\mathbb{C}^{d^2\times d^2} \big)$ be the matrix field with columns $\mathcal{I}(\mathcal{G}_1),\ldots,\mathcal{I}(\mathcal{G}_{d^2})$, which is a fundamental matrix for \eqref{eq:std_form_homog}.
Recalling the notation $\varphi$, $\overline{\varphi}$, and $\tau_\pm$ from \S\ref{sub:the_scalar_x_ray_transform_and_its_mapping_properties},
by the variation of parameters formula for ODEs (\cite[Chapter 3, Theorem 3.1]{Coddington1955})
\begin{equation}\label{eq:W_int}
W \left(z,\zeta \right)=\mathcal{I}^{-1}\mathcal{G}\int_{-\tau_- \left(z,\zeta  \right)}^{0}{\rho^{-1}\mathcal{G}^{-1}\mathcal{I} \left(\widetilde{B}+\widetilde{\Phi}-B-\Phi \right)\circ{\overline{\varphi}}_{\tau} \left(z,\zeta \right)\d \tau}.
\end{equation}
This can also be directly verified by showing that the right-hand side of \eqref{eq:W_int} satisfies \eqref{eq:with_linear_isom}. 
For brevity, we drop writing $\mathcal{I}$ and $\mathcal{I}^{-1}$.

Since $W=0$ on $\partial _+{}^b\!S^* M$, applying $\mathcal{G}^{-1}$ to \eqref{eq:W_int} shows that the integral of $\rho^{-1}\mathcal{G}^{-1} \big(\widetilde{B}+\widetilde{\Phi}-B-\Phi \big)$ over maximal integral curves of $\overline{X}$ is zero.
Reparametrizing the integral over integral curves of $X=\rho\overline{X}$ gives that for any $z\in S^*M^\circ$ 
(c.f. \cite[(3.6)]{Graham2019})
\begin{equation}\label{eq:NonAbelIntegralEq0}
\int_{-\infty}^{\infty}{\mathcal{G}^{-1} \big(\widetilde{B}+\widetilde{\Phi}-B-\Phi \big)\circ \varphi_t \left(z,\zeta \right)}\d t=0.
\end{equation}
Without loss of generality, assume $\rho=\rho_0$ is a geodesic boundary defining function (hence \eqref{eq:metric_normal_form} holds).
Since we cite \cite{Graham2019} often below, we mention that \cite{Graham2019} use $\delta_0,y_0,\eta_0,\omega_0$
instead of $\delta_*,y_*,\eta_*,{{{\omega}}}_*$
respectively, introduced below. %
Fix any $q_*\in\partial  M$ in our coordinates $ \left(\rho_0,y^i \right)$.
By the proof of Lemma 2.8 in \cite{Graham2019}, there exists a neighborhood $\mathcal{U}\subseteq\partial  M$ of $q_*$ 
with closure contained in $\operatorname{dom}{ \left(\rho_0,y^i \right)}$ and $R_0>0$ such that for any $y_*\in\mathcal{U}$ and any $\eta_*\in\partial _-{}^b\!S^* M$ with $\left|\eta_*\right|_{h_0}\geq R_0$,
any integral curve $\sigma_{ \left(y_*,\eta_* \right)}$ of $X$ with initial conditions $ \left(y_*,\eta_* \right)$ at $t\rightarrow-\infty$ is in $\operatorname{dom}{ \left(\rho_0,y^i \right)}$ for all time.
We will prove the inductive step by demonstrating \eqref{eq:ArbDecayIDir} with ``$l$" replaced by ``$l+1$" in a neighborhood of $\mathcal{U}$.

Let $t$ and $\tau$ denote the parametrization variables of integral curves of $X$ and $\overline{X}$ respectively.
Consider the coordinates $ \left(\theta,y,\eta \right)$ on $U_2$ (see remark after \eqref{eq:U2}), consider the notation $\delta_*=\left|\eta_*\right|_{h_0}^{-1}$ and ${{{\omega}}}_*=\delta_*\eta_*$, and write $\sigma_{ \left(y_*,\delta_*,\omega_* \right)} \left(t \right)= \left(\rho_0 \left(t \right),y \left(t \right),\overline{\xi} \left(t \right),\eta \left(t \right) \right)= \left(\theta \left(t \right),y \left(t \right),\eta \left(t \right) \right)$.
We begin with two estimates.
First, by the comment after \cite[(2.24)]{Graham2019}
\begin{equation}\label{eq:RhoDeltaEstimate}
\rho_0 \left(t \right)=\delta_*\sin{ \left(\theta \left(t \right) \right)} \left(1+O_{y_*,{{{\omega}}}_*,\sigma,t} \left(\delta_* \right) \right)
\end{equation}
where ``$O_{\ldots} \left(\delta_*^{k} \right)$" denotes a quantity satisfying $-C\delta_*^{k}\le O_{\ldots} \left(\delta_*^{k} \right)\le C\delta_*^{k}$ for small $\delta_*>0$ where $C\geq0$ is a constant independent of ``...".
We use the same ``..." for $O$ below and hence omit writing it
(we may implicitly reparametrize $t$ below).
Next, by the middle of page 2891 in \cite{Graham2019} (here ${{{\omega}}}_*^\sharp=h_0^{-1}{{{\omega}}}_*$)
\begin{equation}\label{eq:yDotDeltaEstimate}
\dot{y} \left(t \right)=\delta_*\sin^2{ \left(\theta \right)}{{{\omega}}}_*^\sharp+O \left(\delta_*^2 \right).
\end{equation}

Now, let $\Gamma_{\mu j}^k$ and ${\widetilde{\Gamma}}_{\mu j}^k$ be the connection coefficients of $\nabla$ and $\widetilde{\nabla}$ respectively in our coordinates (e.g. $B=\Gamma_{0j}^kE_k^j\d\rho_0+\Gamma_{ij}^kE_k^j\d y^i$) and observe that $(\widetilde{\nabla}-\nabla)(\partial_{x^\mu})=(\widetilde{\Gamma}_{\mu j}^k-\Gamma_{\mu j}^k)E_k^j$.
It is explained at the bottom of page 2890 in \cite{Graham2019} that for every integral curve $\sigma \left(t \right)$ of $X$, the component $\theta : \mathbb{R}\rightarrow \left(0,\pi  \right)$ is a diffeomorphism for $\delta_*>0$ small enough and that $\sfrac{\d \theta}{\d t}=\sin{ \left(\theta \right)} \left(1+O \left(\delta_* \right) \right)$.
Hence, changing variables in \eqref{eq:NonAbelIntegralEq0} gives (in \eqref{eq:NonAbelIntegralEq0} we set $\varphi_t(z,\zeta)=\sigma_{ \left(y_*,\delta_*,\omega_* \right)} \left(t \right)$ which we omit writing right afterwards; in the second equality below we \textbf{1)} use that $ ({\widetilde{\Gamma}}_{0j}^k-\Gamma_{0j}^k)=0$, \textbf{2)} write $ ({\widetilde{\Gamma}}_{ij}^k-\Gamma_{ij}^k ){\dot{y}}^i=\rho_0^l\rho_0^{-l} ({\widetilde{\Gamma}}_{ij}^k-\Gamma_{ij}^k ){\dot{y}}^i$ and $ ({\widetilde{\Phi}}_j^k-\Phi_j^k )=\rho_0^{l+1}\rho_0^{- (l+1 )} ({\widetilde{\Phi}}_j^k-\Phi_j^k )$, and \textbf{3)} apply \eqref{eq:RhoDeltaEstimate} and \eqref{eq:yDotDeltaEstimate}) 
\begin{equation}
	\begin{aligned}
		0&=\int_{0}^{\pi }{\mathcal{G}^{-1} \left((\widetilde{B}-B)+(\widetilde{\Phi}-\Phi) \right)\circ\sigma_{ \left(y_*,\delta_*,\omega_* \right)} \left(\theta \right)\frac{\d\theta}{\sin{ \left(\theta \right)} \left(1+O \left(\delta_* \right) \right)}} \\
		&=\int_{0}^{\pi }{\left[\delta_*\sin{ \left(\theta \right)} \left(1+O \left(\delta_* \right) \right)\right]^l\mathcal{G}^{-1} \left(\rho_0^{-l} \big({\widetilde{\Gamma}}_{ij}^k-\Gamma_{ij}^k \big) \big(\delta_*\sin^2{ \left(\theta \right)}{{{\omega}}}_*^{\sharp i}+O \left(\delta_*^2 \right) \big)E_k^j \right)}\frac{\d\theta}{\sin{ \left(\theta \right)} \left(1+O \left(\delta_* \right) \right)}\\*
		&\qquad +\int_{0}^{\pi }{\left[\delta_*\sin{ \left(\theta \right)} \left(1+O \left(\delta_* \right) \right)\right]^{l+1}\mathcal{G}^{-1} \Big(\rho_0^{-l-1} \big({\widetilde{\Phi}}_j^k-\Phi_j^k \big)E_k^j \Big)}\frac{\d\theta}{\sin{ \left(\theta \right)} \left(1+O \left(\delta_* \right) \right)}.
	\end{aligned}
\end{equation}
We now multiply through the above equation by $\delta_*^{-l-1}$ and let $\delta_*\rightarrow0^+$.
We claim that for fixed $y_*$ and ${{{\omega}}}_*$,
\begin{equation}\label{eq:limits}
	\begin{aligned}
			\mathcal{G}^{-1}\circ\sigma_{ \left(y_*,\delta_*,\omega_* \right)} \left(\theta \right)\rightarrow \mathrm{id},\qquad\rho_0^{-l} \big({\widetilde{\Gamma}}_{ij}^k-\Gamma_{ij}^k \big)&\circ\sigma_{ \left(y_*,\delta_*,\omega_* \right)} \left(\theta \right)\rightarrow\left.\left[\rho_0^{-l} \big({\widetilde{\Gamma}}_{ij}^k-\Gamma_{ij}^k \big)\right]\right|_{y_*}, \\
			\rho_0^{-l-1} \big({\widetilde{\Phi}}_j^k-\Phi_j^k \big)\circ\sigma_{ \left(y_*,\delta_*,\omega_* \right)} \left(\theta \right)\rightarrow&\left.\left[\rho_0^{-l-1} \big({\widetilde{\Phi}}_j^k-\Phi_j^k \big)\right]\right|_{y_*},
	\end{aligned}
\end{equation}
all uniformly over $\theta\in\left[0,\pi \right]$ as $\delta_*\rightarrow0^+$.
Assuming \eqref{eq:limits} for the moment, taking the limit as $\delta_*\to 0^{+}$ yields for all $j,k=1,\dots,d$
$$0=\left.\left[\rho_0^{-l} \left({\widetilde{\Gamma}}_{ij}^k-\Gamma_{ij}^k \right)\right]\right|_{y_*}\omega_*^{\sharp i}\int_{0}^{\pi }{\sin^{l+1}{ \left(\theta \right)}d\theta}+\left.\left[\rho_0^{-l-1} \big({\widetilde{\Phi}}_j^k-\Phi_j^k \big)\right]\right|_{y_*}\int_{0}^{\pi }{\sin^l{ \left(\theta \right)}d\theta}.$$
Varying $y_*$ and ${{{\omega}}}_*^{\sharp i}$ yields that $\rho_0^{-l-1} ({\widetilde{\Phi}}_j^k-\Phi_j^k )$ and each $\rho_0^{-l} ({\widetilde{\Gamma}}_{ij}^k-\Gamma_{ij}^k )$ are zero on $\mathcal{U}$ for all $j,k=1,\dots,d$ and hence $ ({\widetilde{\Phi}}_j^k-\Phi_j^k )\in\rho_0^{l+2}C^\infty (\hat{\mathcal{U}} )$ and $ ({\widetilde{\Gamma}}_{ij}^k-\Gamma_{ij}^k )\in\rho_0^{l+1}C^\infty (\hat{\mathcal{U}} )$ on some neighborhood $\hat{\mathcal{U}}\subseteq\operatorname{dom}{ (\rho,y^i )}$ of $\mathcal{U}$.
This completes the induction.

We return to proving \eqref{eq:limits}.
For the first limit, since matrix inversion is continuous, we can instead prove this with ``$\mathcal{G}^{-1}$" replaced by ``$\mathcal{G}$." Recalling the vector field $\widetilde{X}=\left|\eta\right|_{h_{\rho_0}}^{-1}\overline{X}$ (c.f. \eqref{eq:Xtilde}), writing $\mathcal{G}_i= \left(\mathcal{G}_i \right)_j^kE_k^j$, multiplying \eqref{eq:XBarTranspEqForG} though by $\left|\eta\right|_{h_{\rho_0}}^{-1}$, and pushing the equation to ${}^b\!S^* M$ gives
\begin{equation}
    0=\widetilde{X} \left(\mathcal{G}_i \right)_j^k+\rho_0^{-1}\left|\eta\right|_{h_{\rho_0}}^{-1}\left[ \left(\Phi_r^k \left(\mathcal{G}_i \right)_j^r- \left(\mathcal{G}_i \right)_r^k{\widetilde{\Phi}}_j^r \right)+ \left(\Gamma_{\mu r}^k \left(\mathcal{G}_i \right)_j^r- \left(\mathcal{G}_i \right)_r^k{\widetilde{\Gamma}}_{\mu j}^r \right) v^\mu\right],
\end{equation}
where 
\begin{equation}\label{eq:vmu}
    v^\mu = \begin{cases}\rho_0\overline{\xi}&\mathrm{if} \ \mu=0\\\rho_0^2h^{\mu m}\eta_{m}&\mathrm{if} \ \mu>0\end{cases}
    \quad =\quad 
    \begin{cases}\delta \sin(\theta)\cos(\theta)&\mathrm{if} \ \mu=0\\\delta \sin^2(\theta)h^{\mu m}\omega_m&\mathrm{if} \ \mu>0\end{cases} 
\end{equation} 
 in  the coordinates \eqref{eq:local_ident}-\eqref{eq:U2Param} (recall that $\rho_0=\delta \sin(\theta)$).
 Hence, in terms of them we have  
\begin{equation}
    0=\widetilde{X} \left(\mathcal{G}_i \right)_j^k+(\sin(\theta))^{-1}\left[ \Big(\Phi_r^k \left(\mathcal{G}_i \right)_j^r- \big(\mathcal{G}_i \big)_r^k{\widetilde{\Phi}}_j^r \Big)+ \Big(\Gamma_{\mu r}^k \big(\mathcal{G}_i \big)_j^r- \big(\mathcal{G}_i \big)_r^k{\widetilde{\Gamma}}_{\mu j}^r \Big) v^\mu \right].
\end{equation}
Since  $\Phi,\widetilde{\Phi}\in\rho C^\infty(M;\mathbb{C}^{d\times d})$, all of the coefficients of this PDE for $\mathcal{G}$ are smooth in $(\theta,\delta,y,\omega)$.

Let us reparametrize $\theta\rightarrow s\in\left[0,s_+ \left(y_*,{{{\omega}}}_*,\delta_* \right)\right]$
so that $\sigma_{ \left(y_*,\delta_*,\omega_* \right)} \left(s \right)$ is an integral curve of $\widetilde{X}$.
It follows from \eqref{eq:sol_l_0} and the discussion afterwards that $\theta = s+O(\delta_*)$ (because now $\theta_*=0$), $\delta = O(\delta_*)$, $y=y_*+O(\delta_*)$, and $\omega = \omega_*+O(\delta_*)$.
Hence, setting $\delta_*=0$ we see that $\mathcal{G}\circ\sigma_{ \left(y_*,0,{{{\omega}}}_* \right)} \left(s \right)$ satisfies the ODE $0=\frac{\d}{\d s} \left(\mathcal{G}_i \right)_j^k$ with initial conditions $\mathcal{G}_i \left(0 \right)\cong e_i$.
Thus $\mathcal{G}\circ\sigma_{ \left(y_*,0{{{\omega}}}_* \right)} \left(s \right)=\textrm{id}$ for $s\in\left[0,s_+ \left(y_*,\omega_*,0 \right)\right]$.
Due to the compactness of this $s$-interval at $\delta_*=0$ and the continuity of $\mathcal{G}\circ\sigma_{ \left(y_*,0{{{\omega}}}_* \right)} \left(s \right)$, this shows that $\mathcal{G}\circ\sigma_{ \left(y_*,\delta_*,\omega_* \right)} \left(s \right)$ tends uniformly in $s$ to $\textrm{id}$ as $\delta_*\rightarrow0^+$ for fixed $y_*$ and $\omega_*$.
This implies our claim. 
The second and third limit in \eqref{eq:limits} also follow from $ y=y_*+O(\delta_*)$ and $\rho=O(\delta_*)$.
\end{proof}

\subsection{Proof of Gauge Equivalence}\label{sec:GaugeEquiv}

We prove Theorem~\ref{thrm:GaugeExist} on the recoverability of the non-Abelian X-ray transform up to gauge.
Since $T^{ \left(\Phi,\nabla \right)}=T^{ (\widetilde{\Phi},\widetilde{\nabla} )}$, there exists  $W\in C^\infty(\bSM;\mathbb C^{d\times d})$ satisfying \eqref{eq:TranspEqForW}.
Theorem~\ref{thm:ArbDecay} implies that, by modifying $ \left(\Phi,\nabla \right)$ up to gauge if necessary, we can assume without loss of generality that in \eqref{eq:TranspEqForW} we have $(\widetilde{B}-B)=\widetilde{\nabla}-\nabla\in\rho^NC^\infty \big(M;\mathrm{Hom} \big(TM;\mathbb{C}^{d\times d} \big) \big)$ and $(\widetilde{\Phi}-\Phi)\in\rho^{N+1}C^\infty \left(M;\mathbb{C}^{d\times d} \right)$ for any integer $N\geq0$.
Next, with respect to the matrix inner product $\langle U,V\rangle =\mathrm{trace} \left(UV^\ast \right)$, the linear map $W\mapsto (B+\Phi)W-W(\widetilde{B}+\widetilde{\Phi})$ is skew-Hermitian, e.g., per \cite[Exercise 13.2.4]{Paternain2023}.
Hence by choosing an appropriate isomorphism $\mathbb{C}^{d\times d}\cong\mathbb{C}^{ (d^2 )}$, \eqref{eq:TranspEqForW} can be rewritten into the form \eqref{eq:TranpEqFourier} in Theorem~\ref{MainTranspEqFourierThrm}.
Applying this theorem here implies that $W\in C^\infty \left(M;\mathbb{C}^{d\times d} \right)$ (i.e., $W$ only depends on position).
Let $Q=W+\mathrm{id}$, which observe equals ``$\mathrm{id}$" on $\partial M$.
Plugging this into \eqref{eq:TranspEqForW} and equating elements of $\Omega_0$ and $\Omega_{-1}\oplus\Omega_1$ gives respectively
$$\Phi Q-Q\widetilde{\Phi}-\Phi+\widetilde{\Phi}=\widetilde{\Phi}-\Phi\qquad\mathrm{and}\qquad XQ+BQ-Q\widetilde{B}-B+\widetilde{B}=\widetilde{B}-B.$$
The first equation gives $\widetilde{\Phi}=Q^{-1}\Phi Q$.
Regarding the second equation, since $Q$ only depends on position we have that $XQ=vQ$ (e.g. see \eqref{eq:Xconf}) and so $vQ+BQ-Q\widetilde{B}=0$.
A quick computation shows that this implies that ${\widetilde{\nabla}}_vu=Q^{-1}\nabla_v \left(Qu \right)$ for any $v\in TM$ and any $u\in C^\infty \left(M;\mathbb{C}^d \right)$ .

\section{Surjectivity of \texorpdfstring{$I_0^*$}{I0*} }\label{sec:surj_i0*}

We showed Theorems \ref{thrm:GaugeExist} and \ref{MainTranspEqFourierThrm} assuming Theorem~\ref{HolomIntFactExistThrm}. This in turn relies on Theorem~\ref{thm:surj_I0}. We will show the latter now assuming the next two results, which will be proved in \S\ref{sub:leading_order} and \S\ref{sub:surjectivity_on_sobolev_functions} respectively.
For the 0-Sobolev spaces in Proposition~\ref{prop:solvability} see \S\ref{sub:edge_structures}.

\begin{lemma}\label{lm:surjectivity_terms}Let $(M^{\circ},g)$ be a nontrapping AH manifold of dimension $n+1$ and fix a geodesic bdf $\rho_0$ for $\partial M$.
Let $h\in C^\infty(\partial M)$,
 $s\in \mathbb{C}$ with $\Re s<n$, and $k\in \mathbb{N}_0$, and set $$\tilde h(y,\eta)\coloneqq h(y)\langle\eta\rangle_{h_{0}}^{-s}\log^k(\langle\eta\rangle_{h_{0}})\in C^\infty(\partial_-\bSM),$$ where $\langle\eta\rangle_{h_{\rho_0}}=({1+|\eta|_{h_{\rho_0}}^2})^{1/2}$ for $(y,\eta)\in T^*\partial M\cong \partial_-\bSM$.
Then
\begin{equation}\label{eq:lemma_surjectivity_terms}
	\dfrac{I_0^*\tilde  h(\rho_0,y)-\sigma_{n-1} B\big(\frac{1}{2},\frac{n-s}{2}\big)h(y)\rho_0^{s}|\log (\rho_0)|^k}{ \rho_0^{\Re s}|\log (\rho_0)|^k }\overset{\rho_0\to 0}\to 0\quad \text{pointwise},
\end{equation}
where $\sigma_{n-1}$ is the volume of $\mathbb{S}^{n-1}$  and $B$ is the Beta function.
\end{lemma}

\begin{proposition}\label{prop:solvability}
Let $(M^{\circ},g)$ be simple AH of dimension $n+1$.
Fix $\alpha \in (-n/2,n/2)$, $s\geq 0$, and let $f\in\rho^{\alpha} H_0^{s+1}(M)$. 
Then there exists a unique $u\in \rho^{\alpha} H_0^{s}(M)$ such that $I_0^*I_0 u=f$.
If in addition $f\in \Aphg^{E}(M)$ (necessarily with $\Re E>0$), then
$u\in \mathcal{A}_{-}^{t}(M)$, with $t=\min\{\inf\Re (E),n\}$.
\end{proposition}

We now have:

\begin{proof}[Proof of Theorem~\ref{thm:surj_I0}.] 
First suppose that $\inf \Re E\geq n$.
By Proposition~\ref{prop:solvability}, there exists $u\in \mathcal{A}_-^{n}(M)$ such that $I_0^* I_0 u= f$.
By Proposition~\ref{prop:phg_mapping}\ref{mapping_item1}, $h\coloneqq I_0 u\in \mathcal{A}_-^{n}\big(\overline{\partial_- \bSM}\big)$ and satisfies $I_0^* h= f$.

Now let $f\in\Aphg^E(M)$, where $E$ satisfies $\inf \Re E< n$. Then
\begin{equation}
	f\sim \sum_{\Re (s_m)\to \infty}\sum_{k=0}^{p_m}a_{m,k}(y)\rho_0^{s_m}|\log(\rho_0)|^{k},\quad a_{m,k}\in C^{\infty}(\partial M), \quad (s_m,p_m)\in E,\quad m\geq 1,
\end{equation}
where we may assume upon rearranging that $\Re s_m $ is non-decreasing in $m$. 
Suppose that for $m=1,\dots, J$ we have $\Re s_m= s$ for some $s<n$ and $\Re s_m> s$ for all $m >J$.
In other words, the corresponding terms are the leading order terms, whose corresponding exponents have the same real part but differ in imaginary parts.
Then, as in the statement of Lemma~\ref{lm:surjectivity_terms}, set  
\begin{align}
	\tilde a_{m,p_m}(y,\eta)\coloneqq \frac{1}{\sigma_{n-1} B\big(\frac{1}{2},\frac{n-s_m}{2}\big)}a_{m,p_m}(y)\langle\eta\rangle_{h_{0}}^{-s_m}\log^{p_m}(\langle\eta\rangle_{h_{0}}), \qquad m=1,\dots, J.
\end{align}
Since $\langle\eta\rangle_{h_{0}}^{-1} =\delta (1+\delta^{2})^{-1/2} $ is a bdf for $\partial\big(\overline{\partial_-\bSM}\big)$, $\tilde a_{m,p_m}(y,\eta)\in \Aphg^{E_m}\big(\overline{\partial_-\bSM}\big) $ with $ E_ m \coloneqq \overline{\{(s_m,p_m)\}}$.
By Proposition~\ref{prop:phg_mapping}, $I_0^{*}\tilde a_{m,p_m}\in \Aphg^{E_m'}(M)$ with $ E_m'=E_m \overline\cup E_m\overline\cup\overline{\{(n,1)\}}$, so 
\begin{align}
		I_0^*\big(\tilde a_{m,p_m}\big )\sim	\sum_{\ell=0}^{\infty }\sum_{k=0}^{2p_m+1}b_{m,\ell, k}(y)\rho_0^{s_m+\ell }|\log(\rho_0)|^{k}+ \sum_{j=n}^{\infty}\sum_{t=0}^{T_m} c_{m,j,t}(y)\rho_0^{j}|\log(\rho_0)|^{t},
\end{align}
where  $b_{m,\ell,k}, c_{m,j,t}\in C^{\infty}(\partial M)$.
Above, $T_{m}=1$, unless $s_m$ is an integer, in which case $T_m=2p_m+3$. 
Subtracting $ a_{m,p_m}(y)\rho_{0}^{s_m}|\log(\rho_0)|^{p_m}$, dividing by $\rho_{0}^{s}|\log(\rho_0)|^{p_m} $, and taking a limit as $\rho_{0}\to 0$, Lemma~\ref{lm:surjectivity_terms} implies that $b_{m,0,k}=0$ if $k>p_m$ (that is, no additional higher logarithmic powers are generated by $I_0^*$) and $b_{m,0,p_m}=a_{m,p_m}$.
Hence
\begin{align}
		I_0^*\big(\sum _{m=1}^{J}\tilde a_{m,p_m}\big )\sim& \sum_{m=1}^{J}	\sum_{\ell=0}^{\infty }\sum_{k=0}^{p_m}b_{m,\ell, k}(y)\rho_0^{s_m+\ell }|\log(\rho_0)|^{k}+ \sum_{m=1}^{J}\sum_{j=n}^{\infty}\sum_{t=0}^{T_m} c_{m,j,t}(y)\rho_0^{j}|\log(\rho_0)|^{t}\\*
		\sim& \sum_{m=1}^{J}a_ {m,p_m}(y)\rho_0^{s_m}|\log(\rho_0)|^{p_m}+\tilde f,
\end{align}
where $\tilde f$ is a polyhomogeneous function whose leading order terms correspond to elements of the form $ (s_m,p_m-1)\in E$, provided $p_m-1\geq 0$, i.e. they are of lower order than the leading terms of $f$, which are given by the first sum.
In other words, we have ``matched'' the leading order asymptotic terms in $f$ with the leading order of $I_0^*h$ for a suitable function $h$ on $\overline{\partial_-\bSM}$.
Then in the same way construct a function on $\overline{\partial_-\bSM} $ to match the sum 
\begin{equation}
	\sum_{\{1\leq m\leq J:\ p_m\geq 1\}}(a_{m,p_m-1}(y)-b_{m,0,p_m-1}(y))\rho_0^{s_m}|\log(\rho_0)|^{p_m-1}.
\end{equation}
Proceeding in the same fashion, all terms in the expansion of $f$ corresponding to powers $s_m$ with $\Re s_m=s$ can be matched, and we can continue to match those with $\Re s_m>s$.
Since by definition of an index set there are only finitely many $(s,p)\in E$ with $\Re (s)< n$, in finitely many steps we can   construct a function $h_1\in \Aphg ^{E'}(\overline{\partial_- \bSM})$, where the elements of $E$ and $E'$ in the strip $\{\inf \Re E\leq \Re (s)< n\}$ agree, such that $I_0^{*}h_1-f=f_1\in \Aphg^G(M) $ with $\inf \Re G\geq n$.
Then
we can find $h_2\in \mathcal{A}_-^{n}(\overline{\partial_- \bSM})$ such that $I_0^*h_2=-f_1$ as in the beginning of the proof and set $h=h_1+h_2$ to obtain the result.
\end{proof}

\subsection{Leading Order Asymptotics of \texorpdfstring{$I_0^{*}$}{I0star} -- Proof of Lemma~\ref{lm:surjectivity_terms}} %
\label{sub:leading_order}

\begin{proof}[Proof of Lemma~\ref{lm:surjectivity_terms}]
The geodesic bdf $\rho_0$ together with coordinates on $\partial M$ induce coordinates $(\rho_0,y,\bar\xi,\eta)$ on ${}^bT^*M$ as usual.
Fix $L\gg 0$. For $(\rho_0,y)\in M$ with $0<\rho_0<1/(2L)$,
\begin{equation}\label{eq:decomp}
	I_0^*\tilde  h(\rho_0,y)=\int_{\bSM\big|_{{(\rho_0,y)}}} (B^*_-\tilde h)(\rho_0,y,\bar{\xi},\eta )\d \mu_{(\rho_0,y)}(\bar{\xi},\eta)=\mathsf{I}_1^{+}+\mathsf{I}_1^{-}+\mathsf{I}_2,
\end{equation}
where $\mathsf I_1^{\pm }$, resp. $\mathsf I_2$ is the integral over the region where $|\eta|_{h_{\rho}}\leq L$ and $\pm \bar{\xi}>0$, resp. $|\eta|_{h_{\rho_{0}}}\geq L$.

We first analyze $\mathsf I_1^{\pm }$. In the domain of integration we can parametrize the fiber with the coordinate $\eta$, so that $\bar{\xi}=\pm ({1-\rho_0^{2}|\eta|^{2}_{h_{\rho_0}}})^{1/2}$ and the fiber measure becomes $\d \mu_{(\rho_0,y)}= \frac{(\det h_{\rho_0})^{-1/2}\rho^{n}_0|\d \eta|}{ ({1-\rho_0^{2}|\eta|^{2}_{h_{\rho_0}}} )^{1/2}}$ (this can be checked using \eqref{eq:measure2}).
Since $\tilde{h}$ is bounded on compact subsets of $\partial_-\bSM$, $B_-^{*}\tilde h$ is bounded on $ |\eta|_{h_{\rho_0}}\leq L$ and dominated convergence applies to show that if $\Re s<n$ we have
\begin{align}\label{eq:first_lim}
		 \frac{1}{ \rho_0^{\Re s}|\log (\rho_0)|^k}\mathsf I_1^{\pm }= \frac{1}{ \rho_0^{\Re s}|\log (\rho_0)|^k}\int_{|\eta|_{h_{\rho_0}}\leq L}
		 B_-^*\tilde{h}(\rho_0,y,\eta)\frac{(\det h_{\rho_0})^{-1/2}\rho^{n}_0|\d \eta|}{ ({1-\rho_0^{2}|\eta|^{2}_{h_{\rho_0}}})^{1/2} }\overset{\rho_0\to 0}\to 0.
\end{align}

For $\mathsf I_2$,  parameterize the domain of integration in a fiber corresponding to a fixed base point using $(\theta,\omega)$, where $\bar{\xi}=\cos(\theta)$ and $\eta=\rho_0^{-1}\sin (\theta)  \omega$, $ \omega \in S^{*}_{h_{\rho_0}}\partial M$, with the fiber measure given by \eqref{eq:fiber_measure}. 
In the domain of integration $|\eta|_{h_{\rho_0}}\geq L$, $ \arcsin(L\rho_0)\leq \theta \leq \pi-\arcsin(L\rho_0)$ and note that 
\begin{equation}\label{eq:omega_n}
\begin{aligned}
	\int_{|\omega|_{h_{\rho_0}}=1}\int_{\arcsin(L\rho_0)}^{\pi-\arcsin(L\rho_0)}  \sin^{n-1-s}(\theta)\, \d \theta\, \d \lambda ( \omega)\overset{\rho_0\to 0} \to \int_{|\omega|_{h_{\rho_0}}=1}\int_{0}^{\pi}  \sin^{n-1-s}(\theta)\, \d \theta\, \d \lambda ( \omega)
	\\
	 =2\sigma_{n-1} \int_{0}^{\pi/2}  \sin^{n-1-s}(\theta)\, \d \theta=\sigma_{n-1}B\Big(\frac{1}{2},\frac{n-s}{2}\Big).
	\end{aligned}
\end{equation}
Motivated by this, consider the quantity
\begin{align}
		\frac{ 1}{ \rho_0^{\Re s}|\log (\rho_0)|^k }&\Big(\mathsf I _2(\rho_0,y)-
		h(y)\rho_0^s |\log(\rho_0)|^k
		\int_{|\omega|_{h_{\rho_0}}=1}\int_{\arcsin(L\rho_0)}^{\pi-\arcsin(L\rho_0)}  \sin^{n-1-s}(\theta)\, \d \theta\, \d \lambda ( \omega)\Big)\\* 
		&=\frac{ 1}{ \rho_0^{\Re s}|\log (\rho_0)|^k }\Big(\int_{|\omega|_{h_{\rho_0}}=1}\int_{\arcsin(L\rho_0)}^{\pi-\arcsin(L\rho_0)} (B_-^{*}\tilde h)(\rho_0,y,\theta, {\omega})\sin^{n-1}(\theta)\, \d \theta\, \d \lambda ( \omega)\\* 
		&\hspace{1.3 in} -
		\int_{|\omega|_{h_{\rho_0}}=1}\int_{\arcsin(L\rho_0)}^{\pi-\arcsin(L\rho_0)}   h(y)\rho_0^s |\log(\rho_0)|^k \sin^{n-1-s}(\theta)\, \d \theta\, \d \lambda ( \omega)\Big)\\ 
		&=\int_{|\omega|_{h_{\rho_0}}=1}\int_{\arcsin(L\rho_0)}^{\pi-\arcsin(L\rho_0)} 
			F(\rho_0,y,\theta,\omega)
		 \, \d \theta\, \d \lambda ( \omega)\label{eq:I2},
\end{align}
where
\begin{equation}\label{eq:I2'}
	F(\rho_0,y,\theta,\omega)=\dfrac{(B_-^{*} \tilde h)(\rho_0,y,\theta, {\omega})
	- h(y)\rho_0^{ s} |\log(\rho_0)|^k (\sin(\theta))^{-s}}
		 {\rho_0^{\Re s}|\log (\rho_0)|^k }\sin^{n-1}(\theta).
\end{equation}
Since $\langle\eta\rangle_{h_{0}}^{-1} =\delta (1+\delta^{2})^{-1/2} $, using
\eqref{eq:B-bmap} it follows that $B_-^*(\langle \eta\rangle_{h_0}^{-1})=\frac{\rho_0}{\sin(\theta)}H$, where $H$ is smooth and positive on $\bbSM$ down to $\tilde {\mathsf{Z}}$
and equal to 1 on $\tilde {\mathsf{Z}}$
(in particular
 is
bounded).
Hence for $|\eta|_{h_{\rho_0}}\geq L$, implying $\sin(\theta)\geq L\rho_0$, we estimate 
\begin{align}
		|F(\rho_0&,y,\theta,\omega) |\\*
		\leq &\frac{\sup |h|\Big(\rho_0^{\Re s}|H|^{\Re s}\big|\log (\rho_0)-\log(\sin(\theta))+\log(H)\big|^k+\rho^{\Re s}_0|\log (\rho_0)|^k\Big)}{\rho_0^{\Re s}|\log (\rho_0)|^k }(\sin(\theta))^{n-1-\Re s}\\
		\leq &C\Big(\Big|1+\frac{|\log(\sin(\theta))|}{|\log(\rho_0)|}\Big|^k+1\Big)(\sin(\theta))^{n-1-\Re s}\leq C'(\sin(\theta))^{n-1-\Re s},
\end{align}
where $C'$ is uniform all variables.
A similar computation using \eqref{eq:I2'}
and that $B_-^*|_{\partial_-{}^b\!S^*M}=\textrm{id}$
shows that $F(\rho_0,y,\theta,\omega)\overset{\rho_0\to 0}\to 0$ pointwise.
Since $(\sin(\theta))^{n-1-\Re s}\in L^1([0,\pi])$ for $\Re s<n$,
 by dominated convergence \eqref{eq:I2} vanishes in the limit as $\rho_0\to 0$. 
(To make the bounds of inner integral $\rho_0$-independent one can use a characteristic function as usual,
and for the outer one one can apply the change of variables $ \omega \mapsto \omega = \frac{ \omega}{| \omega|}\in \mathbb{S}^{n-1}$. This introduces a smooth and bounded determinant factor in the integrand which does not affect the limiting process.)
Now the claim follows upon substituting \eqref{eq:decomp} in \eqref{eq:lemma_surjectivity_terms}  and using \eqref{eq:first_lim}, \eqref{eq:omega_n}, as well as the vanishing of the limit of \eqref{eq:I2} that we just showed.
\end{proof}

\subsection{Surjectivity on Sobolev Functions -- Proof of Proposition~\ref{prop:solvability}} %
\label{sub:surjectivity_on_sobolev_functions}
Proposition~\ref{prop:solvability} follows from the following two statements.
The various notations related to 0-(pseudo)differential operators are explained in \S\ref{sub:edge_structures} and \S\ref{ssec:0pdo}.
\begin{lemma}\label{lem:solvability}
Let $(M^{\circ},g)$ be simple AH of dimension $n+1$.
Let $\alpha \in (-n/2,n/2)$, $s\geq 0$. Then $I_0^*I_0$ extends from $C_c^\infty(M^{\circ})$ to a homeomorphism $ \rho^{\alpha} H_0^{s}(M)\to \rho^{\alpha} H_0^{s+1}(M)$.
\end{lemma}
\begin{proof}By 
 \cite[Prop. 4.4 and Prop. 3.3]{Eptaminitakis2022} and \cite[Theorem 1]{Eptaminitakis2022} respectively one has
\begin{align}
\|I_0^*I_0 u\|_{\rho^{\alpha}H_0^{s+1}(M)}\lesssim\| u\|_{\rho^{\alpha}H_0^{s}(M)},&\qquad u\in \rho^{\alpha}H_0^{s}(M),\quad s\in \Rm,\quad \alpha\in(-n/2,n/2)\label{eq:sharp_cont},\\ 
\| u\|_{\rho^{\alpha}H_0^{s}(M)}\lesssim \|I_0^*I_0 u\|_{\rho^{\alpha}H_0^{s+1}(M)},&\qquad u\in \rho^{\alpha}H_0^s(M)\quad s\geq 0\quad \alpha\in(-n/2,n/2).\label{eq:sharp_stab} 
\end{align}
Recall that $L^2_0(M)= L^2(M,\d V_g)$.
These estimates show that for $s\geq 0$ and $\alpha\in(-n/2,n/2)$ the operator
\begin{equation}\label{eq:surj_Sob}
	I_0^*I_0:\rho^{\alpha}H_0^{s}(M)\to \rho^{\alpha}H_0^{s+1}(M)
\end{equation}
is bounded, injective, and has closed range. We will show that it is surjective.

First consider the case $s=0$. 
Since $I_0^*I_0$ has closed range, if it is not surjective there exists $0\neq w\in \rho^{\alpha}H_0^{1}(M)$ such that $\langle w,I_0^*I_0u\rangle_{\rho^{\alpha}H_0^{1}(M)}=0$ for all $u\in \rho^{\alpha}L_0^2(M)=\rho^{\alpha}H_0^{0}(M)$, where we equip $\rho^{\alpha}H_0^{1}(M)$ with the inner product 
\begin{equation}
	\langle\cdot, \cdot  \rangle_{\rho^{\alpha}H_0^{1}(M)}=\langle\cdot, \cdot  \rangle_{ \rho^{\alpha}L_0^2(M)}+\langle\d \,\cdot,\d \, \cdot  \rangle_{ \rho^{\alpha}L_0^2(M)}.
\end{equation}
Using that $I_0^{*} $ is the formal adjoint of $I_0$ in the functional setting $L^2(M,\d V_g)\to L^2(\partial_-\bSM, \d \Sigma_\partial)$ (see e.g., \cite[eq. (2-5)]{Eptaminitakis2022}),
for all $\varphi,\psi\in C_c^\infty(M^{\circ})$ we have 
\begin{equation}
\begin{aligned}
	\langle \psi,I_0^*I_0\varphi\rangle_{\rho^{\alpha}H_0^{1}(M)}
	=&\langle \psi,I_0^*I_0\varphi\rangle_{\rho^{\alpha}L_0^2(M)}+\langle \d  \psi,\d  I_0^*I_0\varphi\rangle_{\rho^{\alpha}L_0^2(M)}\\ 
	=&\langle \rho ^{-2\alpha}\psi,I_0^*I_0\varphi\rangle_{L_0^2(M)}+\langle\rho^{-2\alpha} \d  \psi,\d  I_0^*I_0\varphi\rangle_{L_0^2(M)}\\ 
	=&\langle (\rho ^{-2\alpha}-\mathrm{div} \rho^{-2\alpha} \d  )\psi,I_0^*I_0\varphi\rangle_{L_0^2(M)}\\
		=&\langle  \rho^{2\alpha }I_0^*I_0\rho ^{-2\alpha}Q\psi,\varphi\rangle_{\rho^{\alpha}L_0^2(M)} ,
\end{aligned}\label{eq:IBP}	
\end{equation}
where $Q\coloneqq 1-\rho^{2\alpha}\mathrm{div} \rho^{-2\alpha}\d \in \Diff_0^{2}(M)$.
Using \eqref{eq:sharp_cont} we have the boundedness properties
\begin{equation*}
	\rho^{\alpha}H^{1}_0(M)\overset{Q }{\to }\rho^{\alpha}H^{-1}_0(M)\overset{\rho^{-2\alpha}}{\to }\rho^{-\alpha}H^{-1}_0(M)\overset{I_0^*I_0}{\to }\rho^{-\alpha}L_0^2(M)\overset{\rho^{2\alpha}}{\to }\rho^{\alpha}L_0^2(M),\quad \alpha\in(-n/2,n/2).
\end{equation*}
Since $C_c^{\infty}(M^{\circ})$ is dense in $\rho^{\alpha }H_0^{s}(M) $ for $s\geq 0$, this and \eqref{eq:IBP} implies that 

\begin{equation*}
	0=\langle w,I_0^*I_0u\rangle_{\rho^{\alpha}H_0^{1}(M)}=\langle  \rho^{2\alpha }I_0^*I_0\rho ^{-2\alpha}Q w,u\rangle_{\rho^{\alpha}L_0^2(M)}.
\end{equation*}
Since $u\in \rho^{\alpha}L_0^2(M)$ was arbitrary, we conclude that 
\begin{equation}\label{eq:I*I}
	I_0^*I_0\rho ^{-2\alpha}Q w=0,\quad \text{where}\quad \rho ^{-2\alpha}Q w\in \rho^{-\alpha}H_0^{-1}(M).
\end{equation}
By the proof of \cite[Prop. 6.1]{Eptaminitakis2022},  there exist operators $B\in \Psi_0^{1,n+1,n+1}(M)$ and $K\in \Psi_0^{-\infty,\mathcal{F}}(M)$, where $F_f\geq 1$, $F_\ell$, $F_r\geq n$, such that 
\begin{equation}
	BI_0^{*}I_0=\mathrm{Id}-K\text{ on }\rho^{\alpha}L^2(M), \quad \alpha \in (-n/2,n/2).\label{eq:parametrix1}
\end{equation}
For any $m,N\in \Rm$ and $\alpha\in(-n/2,n/2)$ they are bounded\footnote{The boundedness properties are stated in \cite[Prop. 6.1]{Eptaminitakis2022} only for nonnegative order Sobolev spaces, however they hold more generally using the spaces of 0-pseudodifferential operators in which they live and Prop. 3.3 there.} 
\begin{equation}\label{eq:parametrix2}
	B\colon \rho^{\alpha}H_0^{m}(M)\to \rho^{\alpha}H_0^{m-1}(M)\quad \text{and}\quad K\colon  \rho^{\alpha}H_0^{m}(M)\to \rho^{\alpha}H_0^{N}(M).
\end{equation}
Approximating $w$ by elements of $C_c^{\infty}(M^{\circ})$ and using \eqref{eq:I*I}-\eqref{eq:parametrix2}
 imply that $\rho ^{-2\alpha}Qw\in \rho^{-\alpha}L_0^2(M)$, so
$Qw$
vanishes by \eqref{eq:sharp_stab} and \eqref{eq:I*I}.
Thus
\begin{equation*}
	0=\langle Qw,w\rangle _{\rho^{\alpha}L_0^2(M)}=\|w\|^{2} _{\rho^{\alpha}L_0^2(M)}+\|\d w\|^{2} _{\rho^{\alpha}L_0^2(M)},
\end{equation*}
hence $w=0$, which is a contradiction and shows surjectivity of the map \eqref{eq:surj_Sob} for $s=0$.
For the general case $s\geq 0,$ given $f\in\rho^{\alpha }H^{s+1}_0(M)\subset \rho^{\alpha }H_0^{1}(M)$, there exists a unique $u\in \rho^{\alpha }L_0^2(M)$ such that $I_0^*I_0u=f$; using the operators $B$ and $K$ above it follows that $u\in\rho^{\alpha } H^s_0(M)$.
\end{proof}

The proof of the following lemma uses induction and follows almost verbatim the proof of \cite[Lemma 6.2]{Eptaminitakis2022}, so we will only provide a sketch explaining the idea. We  refer the reader there for the details and the more formal inductive statements.

\begin{lemma}\label{lm:phg_normal_op}
	Let $(M^{\circ},g)$ be simple AH, of dimension $n+1$. Let $f\in \Aphg^{E}(M)$ with $\Re E>0$,
which implies that there exists $\alpha \in (-n/2,n/2)$ such that $f\in\rho^{\alpha} H_0^{1}(M)$.
	For any such $\alpha,$ let $h \in \rho^{\alpha}L^2(M)$ be the unique solution of $I_0^*I_0 h=f$, as in Lemma~\ref{lem:solvability}.
	Then $h\in \mathcal{A}_{-}^{t}(M)$, with $t\coloneqq\min\{\inf\Re (E),n\}.$
\end{lemma}

\begin{proof}[Proof (sketch).]
Let $\alpha$ be as in the statement. 
By the uniqueness in the statement of Lemma~\ref{lem:solvability} and the fact that $\alpha''<\alpha'\implies \rho^{\alpha'}L^2_0(M)\subset \rho^{\alpha''}L^2_0(M)$,  $h$ does not depend on $\alpha$.
In other words, $h\in \rho ^{\alpha }L^2_0(M)$ for all $\alpha\in(-\infty,n/2)\cap (-\infty,\inf \Re E-n/2)$.
For any such $\alpha,$ we will show that $h\in \rho^{\alpha+n/2}H_b^{\infty}(\oM)$, which amounts to showing that 
\begin{equation}\label{eq:conormality}
	V _m\cdots V_1h\in \rho^{\alpha}L^2_0(M)=\rho^{\alpha+n/2}L^2_b(\oM)\quad \forall\quad V_j\in \mathcal{V}_b(M).
\end{equation}
This will imply that $h\in \rho^{\alpha'}H_b^{\infty}(\oM)$ for all $\alpha'<\min\{n,\inf \Re E\}$, that is,  $h\in \mathcal{A}_-^{t}(\oM)$.

We now show \eqref{eq:conormality}. Applying  \eqref{eq:parametrix1} to $h$ we find
	$h=Bf+Kh.$
By \cite[Proposition 3.2]{Eptaminitakis2022}, $Bf\in \Aphg^{F}(M) $, where $ F=\overline{\{(n+1,0)\}}\ocup E$. 
This implies
\begin{equation}\label{eq:conormal_Bf}
 	V _m\cdots V_1Bf\in \Aphg^{F}(M)\subset \rho^{\alpha}L^2_0(M) \quad \forall\quad  V_j\in \mathcal{V}_b(M).
 \end{equation} 

As explained in the proof of \cite[Lemma 6.2]{Eptaminitakis2022},  if $V \in \mathcal{V}_b(M)$, $VK\in \Psi_0^{-\infty,\mathcal{F'}}(M)$, where $\mathcal{F}'=(F_\ell,F_r,F_f-1)$.
Thus, since $F_f\geq 1$, by \cite[Proposition 3.3]{Eptaminitakis2022} one sees that for any $V_1\in \mathcal{V}_b(M)$,
\begin{equation}\label{eq:V1}
	V_1 h=V_1 Bf+V_1K h\in \rho^{\alpha}L^{2}_0(M) .
\end{equation}
Subsequently for any $V_2\in \mathcal{V}_b(M)$,
\begin{equation*}
	V_2V_1 h=V_2V_1 Bf+V_2V_1Kh=V_2V_1 Bf+(V_1K)V_2h-[V_1K,V_2]h. 
\end{equation*}
The first two terms are in $\rho^{\alpha}L^{2}_0(M) $ by \eqref{eq:conormal_Bf} and \eqref{eq:V1} respectively.
The third is  in $\rho^{\alpha}L^{2}_0(M) $ again by \cite[Proposition 3.3]{Eptaminitakis2022} and the fact that from \cite[Proposition 3.30]{Mazzeo1991}, $Q\in \Psi_0^{-\infty,\mathcal{E}}(M)$ and $V\in \mathcal{V}_{b}(M)\implies  [Q,V]\in \Psi_0^{-\infty,\mathcal{E}}(M)$. Proceeding inductively, we obtain the claim.
\end{proof}

\section{Holomorphic Integrating Factors}\label{sec:holo_int}
In this section we prove Theorem~\ref{HolomIntFactExistThrm}.

\begin{proof}[Proof of Theorem~\ref{HolomIntFactExistThrm}]
	
By Theorem~\ref{thm:surj_I0}, there exists a real valued $h\in C^{\infty}(\partial_-\bSM)$ 
such that $I_0^*h= 2\pi f$ and it is of the form
\begin{equation}
	h=b_{01}(y)\log \delta+b_{00}(y)+\tilde h(y,\delta), \quad \tilde h\in \mathcal A_-^1(\overline {\partial_- \bSM}), \quad \delta \ll 1.
\end{equation}
We can now extend $h$
to be constant
along the flow of the normalized geodesic vector field $\tilde{X}$, see \eqref{eq:Xtilde}, that is, consider $\tilde p= I^{*} h$.
 Proposition~\ref{prop:phg_mapping}\ref{mapping_item2}, for any $t_0\in (0,1)$ we have  $ \tilde p\in \Aphg^{\mathcal{H}}(\bbSM)+\mathcal{S}^{\mathfrak{t}}(\bbSM)$, where $ \mathcal{H}=( \overline{\{(0,1)\}},0,0 )$ and $\mathfrak{t}=(t_0,0,0)$, and it is real valued.
By construction we have $X\tilde p=0$ over the interior.
Using the notations introduced immediately before the statement of Proposition~\ref{prop:phg_mapping}, by the properties of the index sets it follows that  there exists a constant $B$ such that
	$|\tilde p|\leq B |\log (x_1)|$, where without loss of generality, $x_1$ takes values in $[0,1/2)$.

Now consider $p\coloneqq \tilde p \circ \overline \Psi_g^{-1}\in \Aphg^{\mathcal{H}}( [\OSM;\mathsf Y] )+\mathcal{S}^{\mathfrak{t}}( [\OSM;\mathsf Y] )$, with the index sets transferred to the boundary hypersurfaces of $ [\OSM;\mathsf Y] $ by means of $\overline\Psi_g$.
Since $\overline{\Psi}_g:\bbSM\to [\OSM;\mathsf Y]$ is a diffeomorphism, $\tilde x_j\coloneqq x_j \circ \overline{\Psi}_g^{-1}$ are defining functions for the corresponding boundary hypersurfaces of $[\OSM;\mathsf Y]$ and therefore 
\begin{equation}\label{eq:p_estimate}
	|p|\leq \tilde B|\log(\tilde x_1)|
\end{equation} for some $\tilde B>0$.
Letting $p=\sum_{k=-\infty}^{\infty}p_k$
 denote $p$'s vertical Fourier decomposition, over $SM^\circ $
\begin{equation}
\label{XqEqual0Fouirer}
0=Xp=\sum_{k=-\infty}^{\infty}\underbrace{\left(\eta_+p_{k-1}+\eta_-p_{k+1} \right)}_{\in\Omega_k}\quad \Longrightarrow\quad\text{for each } k,\quad \eta_+p_{k-1}+\eta_-p_{k+1}=0,
\end{equation}
recall \S\ref{sec:vertical_fourier} for $\eta_\pm $.
We also have that 
\begin{equation}
	p_0=\frac{1}{2\pi}\int_{0}^{2\pi }p \left(x,y,\theta \right)d\theta
=\frac{1}{2\pi}I_0^*h=f.
\end{equation}

Now, let $w :[\OSM;\mathsf Y]\rightarrow\mathbb{C}$ be given by
\begin{equation}\label{eq:defn_w}
	w=p_0+2\sum_{k=1}^{\infty} p_{2k}=(\id +i H)\hat{p}, \qquad \hat p= \frac{p+p\circ A_\sharp}2,
\end{equation}
where recall that $A_\sharp:[\OSM;\mathsf Y]\to [\OSM;\mathsf Y]$ is the antipodal involution
and $H$ is the fiberwise Hilbert transform in Definition~\ref{HilbertTransfDef}.
Since $  A_\sharp$ interchanges the two components of the front face, as well as the two remaining boundary hypersurfaces, $\hat p$ still
 satisfies \eqref{eq:p_estimate}
(possibly upon enlarging the constant).
It  is also real.
Since the Hilbert transform is a real operator, $iH\hat p$ is purely imaginary, so that if $m\geq \tilde B$
\begin{equation}
	\rho^{m} |e^{w}|\overset{\eqref{eq:pi_bmap}}\leq C\tilde x_1^{m}\tilde x_2^{m} \tilde x_3^{m} e^{\hat p}\leq  C \tilde x_{1}^{m-B}\tilde x_2^{m}\tilde x_3^{m}<\infty.
\end{equation}
Finally, we check using \eqref{XqEqual0Fouirer} that 
\begin{equation}
 	Xw=\eta_-p_2+\eta_-p_0 =-\eta _+p_0+\eta_-p_0=-\eta _+f+\eta_-f=\varphi.
 \end{equation}

We now check the claim regarding $P w$ for $P\in \Diff_e^{*}(\OSM)$. 
By Proposition~\ref{prop:HilbMapProp}, \eqref{eq:defn_w} and the definition of Sobolev spaces it suffices to show it for $\hat p$ instead of $w$. 
First, note that $ \hat p\in \Aphg^{\mathcal{H}}( [\OSM;\mathsf Y] )+\mathcal{S}^{\mathfrak{t}}( [\OSM;\mathsf Y] )$, so there exists $m$ such that $(\tilde x_1\tilde x_2\tilde x_3)^{m}\hat p \in L^\infty([\OSM;\mathsf Y])$.
Thus by \eqref{eq:pi_bmap}, $\rho^{m}\hat p\in L^{\infty}(\OSM)$ and $ \rho^{m'}\hat p\in L^2(\OSM)$ for some $m'$ (here we still
use
$\hat p$ for the almost everywhere defined function on $\OSM$ obtained upon pulling back by the inverse of $\beta_0\big|_{[\OSM;\mathsf Y]^{\circ}}$, which is a diffeomorphism).
Now let  $W\in \mathcal V_e(\OSM)$.
In general, $W$ is not tangent to $\mathsf{Y}$. 
In terms of global isothermal coordinates
as in \S\ref{sec:vertical_fourier},
we can consider a  function $\chi \in C^{\infty }(\OSM)$ given by
$\chi = \sin^{2}(\omega-\theta)+\rho^{2}$ near $\partial \OSM$ and extended to be positive elsewhere, where $\omega =\arg(x+iy)\in \R/2\pi \mathbb{Z}$: since locally $\mathsf Y_+=\{\rho=0,\omega=\theta \}$ and $\mathsf Y_-=\{\rho=0,\omega=\theta +\pi \}$,
$\chi=0$ exactly on $\mathsf{Y}$ and $\beta_0^{*}\chi $ is the square of a defining function for the front face in $[\OSM;\mathsf Y]$.
Hence, in the notation of the proof of Proposition~\ref{lm:compactifications}, it is locally a smooth non-vanishing multiple of $\rho_0^{2}+|b|^{2}_{h_{\rho_0}}$, and one  can use 
\eqref{eq:Psi_diffeo_1}-\eqref{eq:Psi_diffeo_2} to check that 
$\beta_0^{*}\chi =e^{\tilde f}\tilde x_2^{2}\tilde x_3^{2}$
for some $\tilde f\in C^\infty([\OSM;\mathsf Y])$.
In particular, $\chi W $ is tangent to $\mathsf Y$ (since it vanishes there), so as mentioned in \S\ref{sub:manifolds_with_corners_p_submanifolds_blowups} it lifts to a tangent vector field on $[\OSM;\mathsf{Y}]$.
Therefore, $\beta_0^{*} (W\hat p\big|_{SM^{\circ }})=e^{-\tilde f}\tilde x_2^{-2}\tilde x_3^{-2}\tilde W  \hat p$, where $\tilde W \coloneqq \beta_0^{*}(\chi W)\in \mathcal{V}_b([\OSM;\mathsf{Y}])$.
Since $\hat p\in \Aphg^{\mathcal{H}}( [\OSM;\mathsf Y] )+\mathcal{S}^{\mathfrak{t}}( [\OSM;\mathsf Y] )$ and $\mathcal V_b$ preserves polyhomogeneous and conormal spaces, $\beta_0^{*}  W\hat p\in \Aphg^{\mathcal{H}'}( [\OSM;\mathsf Y] )+\mathcal{S}^{\mathfrak{t}'}( [\OSM;\mathsf Y] )$, where $\mathcal{H}'=\mathcal{H}-(0,2,2)$ and $\mathfrak{t}'=\mathfrak{t}-(0,2,2)$.
In particular, $(\tilde x_1\tilde x_2\tilde x_3)^{m+2}\beta_0^{*} W\hat p \in L^\infty([\OSM;\mathsf Y])$, and $\rho^{m'+2}W\hat p\in L^{2}(\OSM)$ as before.
Applying more edge vector fields shows the result by induction.

For the anti-holomorphic case, one can follow the construction above for $f$ replaced by $-f$ and consider the complex conjugate of the resulting holomorphic integrating factor
$w$.
\end{proof}

\appendix

\section{Edge Structures, Function Spaces, Blow-ups, 0-\texorpdfstring{$\Psi$}{Ps}DO}
\label{Sec:ConormPolyhomFuncSp}

For the convenience of the reader and to establish conventions, we provide background material that is by now largely standard in geometric microlocal analysis.

\subsection{Edge Structures} \label{sub:edge_structures}A standard reference for this section is \cite{Mazzeo1991}. An \textbf{edge space}  is a smooth compact manifold with boundary $M^{k+q+1}$ together with a fibration
 $\pi \colon \partial M\to Y^{k}$ of the boundary 
  with fiber $Z^{q}$.
 The associated \textbf{edge structure}
  is the Lie algebra $\mathcal{V}_e(M)$ of vector fields which are tangent to the fibers of $\pi$ at the boundary. 
Near $\partial M$ we can use coordinates $(x,y,z)$, where $x$ (here and throughout) is a bdf, the $y^j$'s  result from  extending to $M$  lifts of coordinates for $Y$ to $\partial M$, and the $z^{i}$'s restrict to coordinates on each  fiber at $\partial M$.  
Then $\mathcal{V}_e(M)$ is locally spanned over $C^\infty(M)$ by $x\partial_x,x\partial_{y^{j}}, \partial _{z^{i}}$, $j=1,\dots,k$, $i=1,\dots,q$. 
Associated with the edge structure there is a vector bundle ${}^eTM$ of rank $k+q+1$ which is canonically isomorphic with $TM$ over $M^{\circ}$ and whose smooth sections are exactly the vector fields in $\mathcal{V}_e(M)$. 
The sections of the corresponding dual bundle ${}^eT^*M$ are locally spanned by $ x^{-1}\d x$, $x^{-1}\d y^{j}$, $\d z^{i}$ in the coordinates above. An \textbf{edge metric} is a positive definite section of $\mathrm{Sym}^{2}({}^e T^{*}M)$.
As mentioned in Remark~\ref{rmk:edge}, the 0-sphere bundle of an AH surface is a typical example of an edge space 
and the Lie algebra generated by the vector fields in \eqref{eq:Xconf}  is an edge structure. 
These vector fields are orthonormal for the Sasaki metric (see \cite{Paternain2023}, so the latter is an example of an edge metric.

Given an edge space $(M,\pi)$, denote by $L_e^{2}(M)$ the space of measurable functions which are square integrable with respect to the measure induced by any edge metric. 
The edge Sobolev spaces\footnote{Our weight convention differs slightly from \cite{Mazzeo1991}; our space corresponds to $x^{\delta+ \frac{k+1}{2}}H_e^{\ell }(M)$ there.} for $\ell\in \mathbb{N}_0$ and $\delta\in \R$ are defined as
\begin{equation}\label{eq:HeSpace}
	x^{\delta}H_e^{\ell }(M)=\{u\in x^{\delta}L_e^2(M)\colon V_1\cdots V_j u\in x^{\delta}L_e^2(M)\ \forall j\leq \ell, \ V_j\in \mathcal{V}_e(M)\}.
\end{equation}
For general $\ell>0$, the definition is obtained by interpolation and for $\ell<0$ by duality.
Upon locally choosing vector fields spanning $\mathcal V_e(M)$ over $C^\infty(M)$ and using a partition of unity, we can turn $x^\delta H_e^\ell(M)$ into Hilbert spaces; different choices give rise to equivalent norms.
It is not hard to check that $x^{\delta}H_e^{\ell }(M)$ is a local $C^\infty(M)$ module, implying that if $E\to M$ is a vector bundle, the space  $x^{\delta}H_e^{\ell }(M;E)$ of $x^{\delta}H_e^{\ell }$-regular sections of $E$ is well defined.
Moreover, for $\ell\geq 0$ the space $C_{c}^{\infty}(M^\circ;E)$ of compactly supported sections of $E$ is dense in $x^{\delta}H_e^{\ell }(M;E) $ (this follows similarly to \cite[Lemma 3.9]{Lee2006}).
The space of \textbf{edge differential operators}, that is, of differential operators consisting of locally finite sums of at most $m$-fold products of edge vector fields will be denoted by $\Diff_e^{m}(M)$.
Any $P\in \Diff_e^{m}(M)$ is bounded $x^{\delta}H_e^{\ell }(M)\to x^{\delta}H_e^{\ell -m}(M)$. 
We write $\Diff_e^*(M)\coloneqq\bigcup_{m\in \mathbb{N}_0}\Diff_e^m(M)$.

If the base of the fibration in the edge space above is a point, the corresponding edge structure consists of the vector fields which are tangent to the boundary, denoted in the introduction with $\mathcal{V}_b(M)$ and locally spanned by $x\partial _x$, $\partial_{z^i}$. 
On the other hand, if the fiber of the fibration is trivial, the edge structure consists of the vector fields on $M$ which vanish at the boundary; we previously denoted it by $\mathcal{V}_0(M)$ and it is locally spanned by $x\partial_x$, $x\partial_{y^{j}}$.
In both cases, the discussion regarding the associated bundles, metrics, and Sobolev spaces follows \textit{mutatis mutandis}, replacing the ``$e$'' superscripts with ``$b$'' and ``0'' respectively.
Note that $L_0^2(M)$ corresponds exactly to the $L^2$ space defined with respect to the measure induced by an AH metric on $M^{\circ}.$

\subsection{Conormality, Polyhomogeneity, \texorpdfstring{$\mathrm{b}$}{b}-Maps}\label{sec:polyhomongeneity}

More details on the material of this section can be found in \cite{Melrose1992}. 
Let $Y$ be a compact manifold with corners, with embedded, connected  boundary hypersurfaces.
Such manifolds are locally modeled on $\R^{n}_{k}=[0,\infty)^{k}\times \R^{n}$, where $k$ corresponds to the \textbf{codimension of the corner} at the origin.
As in \S\ref{sub:edge_structures}, denote  by $\mathcal{V}_b(Y)$ the Lie 
algebra of vector fields that are tangent to all boundary hypersurfaces of $Y$, and by ${}^bTY$ the vector bundle whose smooth sections are the vector fields in $\mathcal{V}_b(Y)$.
If $x_{j}, $ $j\in \{1,\dots, J\}$, is a collection of defining functions for the boundary hypersurfaces $H_1,\dots, H_J$ of $Y$ and $\mathfrak{t}\coloneqq(t_1,\dots, t_J)\in \mathbb{R}^{J}$ is a multi-weight, we write $x^{\mathfrak{t}}:=\prod _{j=1}^{J}x_{j}^{t_j}$.
If $t\in \mathbb{R}$, writing $x^{t}$ will indicate that $t_j=t$ for all $j$.
We denote by $\Omega_b(Y)\to Y$ the b-density bundle, trivialized by $x^{-1}\mu$, where $\mu$ is a smooth positive density on $Y$.
In accordance with \S\ref{sub:edge_structures}, $L^2_b(Y)$ is the space of measurable functions which are square integrable with respect to a positive section of $\Omega_b(Y)$. 
Define the $L^{2}$-based conormal spaces as
\begin{equation}
	x ^{\mathfrak{t}} H_b^{\infty}(Y)\coloneqq\{u\in x^{\mathfrak{t}} L^2_b(Y):V_1\cdots V_ku\in x^{\mathfrak{t}}L^2_b(Y)\quad \forall k\in \mathbb{N}_0, \text{where}\ V_j\in \mathcal{V}_b(Y) \},
\end{equation}
where the vector fields act  distributionally in the interior.
The $L^{\infty}$-based conormal spaces are
\begin{equation}
	\mathcal S^{\mathfrak{t}}(Y)\coloneqq\{u\in x^{\mathfrak{t}} L^{\infty}(Y):V_1\cdots V_ku\in x^{\mathfrak{t}}L^\infty(Y)\quad \forall k\in \mathbb{N}_0, \text{where}\ V_j\in \mathcal{V}_b(Y) \}.
\end{equation}
Notice that we use the opposite sign convention from \cite{Melrose1992}.
Then, for $\mathfrak{s}\in \R^{J}$ let 
\begin{equation}\label{eq:conormal}
	\mathcal{A}_-^\mathfrak{s}(Y)\coloneqq\bigcap_{\mathfrak{t}<\mathfrak{s}}x ^{\mathfrak{t}} H_b^{\infty}(Y)=\bigcap_{\mathfrak{t}<\mathfrak{s}}\mathcal{S} ^{\mathfrak{t}} (Y),
\end{equation}
where $ \mathfrak{t}<\mathfrak{s}$ means $t_j<s_j$ for all $j=1,\dots, J$ and the second equality follows from \cite[(18)]{Melrose1992}.

A \textbf{smooth index set} is a discrete subset $E\subset \mathbb{C}\times \mathbb{N}_0$ satisfying the following three requirements:
\begin{align}
 	|(s_j,p_j)|\to \infty&\implies \Re (s_j)\to \infty,\label{eq:req1}\\
 	 (s,p)\in E&\implies (s+k,p)\in E\quad \forall k\in \mathbb{N}_0,\label{eq:req2}\\ 
 	 (s,p)\in E&\implies (s,k)\in E\quad \forall k\in \mathbb{N}_0 \text{ with }0\leq k<p.\label{eq:req3}
\end{align} 
Given $(s,p)\in \mathbb{C}\times \mathbb{N}_0$, we denote by $\overline{\{(s,p)\}}$ the smallest smooth index set containing $(s,p)$ and we assume always that index sets are smooth.
If $a\in \Rm$, we write $\Re E>a$  if $\Re(s)>a$ for all $(s,p)\in E$ and $\Re {E}\geq a^-$ if $\Re s\geq a$ for all $(s,p)\in {E}$.
We write $\Re {E}\geq a$  if $\Re {E}\geq a^-$ and there is no element $(s,p)$ in ${E}$ with $\Re(s)=a$ and $p\neq 0$.
If $a\in \mathbb{C}$, $E+a\coloneqq \{(s+a,p):(s,p)\in E\}$, and we also use $a$ as a shorthand for $\overline{\{(a,0)\}}$.
Finally, let $\inf\Re(E)\coloneqq \inf\{\Re(s):(s,p)\in E\}$. If $E\neq \emptyset$, the infimum is attained.

We define polyhomogeneous conormal functions.
If $Y$
is a compact manifold with boundary and $x$ a bdf,  we say that $u\in \Aphg^{E}(Y)$, where $E $ is a smooth index set, if in terms of a product decomposition $[0,\epsilon)_x\times \partial Y$ of a collar neighborhood of $\partial Y$ it admits an asymptotic expansion
\begin{equation}
	u\sim \sum_{\Re (s_m)\to \infty}\sum_{k=0}^{p_m}a_{m,k}x^{s_m}\log^{k}(x),\quad a_{m,k}\in C^{\infty}(\partial Y), \quad (s_m,p_m)\in E.\label{eq:phg}
\end{equation}
See \cite{Melrose1992} for a more detailed explanation of ``$\sim$''.
By \eqref{eq:req2}-\eqref{eq:req3}, the definition is independent of the choice of $x$ and the product decomposition of the collar neighborhood.
Note that $\Re(E)\geq 0\implies \Aphg^{E}(Y)\subset L^\infty(Y)$.
If $ Y$ is a manifold with corners up to codimension $k$ and  $\mathcal{E}=(E_1,\dots,E_J)$ is a family of index sets corresponding to its boundary hypersurfaces $H_1,\dots,H_J$, the space $\Aphg^{\mathcal E}(Y)$ is defined inductively in the maximal codimension of the faces: if the maximal codimension is 0, $\Aphg^{\mathcal E}(Y)=C^\infty(Y)$. 
For the inductive step, assuming the definition is known for manifolds with corners up to codimension $k-1$, one demands that for each $j$ we have an expansion near $H_j$ of the form
\begin{equation}
	u\sim \sum_{\Re (s_m)\to \infty}\sum_{k=0}^{p_m}a_{m,k}x_{j}^{s_m}\log^{k}(x_{j}),\quad a_{m,k}\in \Aphg^{\mathcal{E}'_j}(H_j), \quad (s_m,p_m)\in E_{j},
\end{equation} 
where $x_j$ is a defining function of $H_j$ and $\mathcal{E}'_j$ is a collection of index sets that assigns to each  boundary hypersurface $G$ of $H_j$ the index set in $\mathcal{E}$ corresponding to the unique boundary hypersurface $F$ of $Y$ such that $G$ is a component of $F \cap H_j$.
Note that if $f\in \Aphg^{\mathcal{E}}(Y)$, where $\mathcal{E}=(E_1,\dots,E_J),$ and we write $\mathfrak{e}\coloneqq (e_1,\dots,e_J)$ with $e_j=\inf \Re E_j$, then 
$f\in \mathcal{A}_-^\mathfrak{e}(Y)$. 

Conormal and polyhomogeneous conormal functions behave well under pullback and pushforward with respect to maps that are well adapted to the geometry.

\begin{definition}\label{def:b-subm}
 	 Suppose  that $Y$, $Y'$ are manifolds with corners that have embedded boundary hypersurfaces and let $x_j$, $x'_i$, $j\in \{1,\dots,J\}$, $ i\in \{1,\dots,I\}$ be collections of defining functions for the boundary hypersurfaces of $Y$, resp. $Y'$.
 	 A smooth map $f:Y\to Y'$  is called an \textbf{interior b-map} if for every $i\in \{1,\dots,I\}$ there exist  integers $e(i,j)\geq 0$ and an $h_i\in C^{\infty}(Y)$ such that
\begin{equation}\label{eq:bmap}
 	f^{*}x'_i = e^{h_{i}}\prod_{j=1}^{J}x_{j}^{e_f(i,j)}.
\end{equation} 
A b-map $f:Y\to Y'$ is a \textbf{b-fibration} if the following two conditions hold:
\begin{enumerate}
 	\item The b-differential $\d^{b} f\big|_p:{}^{b}T_pY\to {}^{b}T_{f(p)}Y'$ (defined by extension of the usual differential from the interior) is surjective for every $p\in Y$ (that is, $f$ is a \textbf{b-submersion})
 	\item\label{item2} for each $j$ there exists at most one $i$ such that $e_f(i,j)\neq 0$ (that is, no boundary hypersurface of $Y$ is mapped into a corner of $Y'$).\footnote{One can check that this requirement is equivalent to a b-map being b-normal in the terminology of \cite{Melrose1992}.} 
\end{enumerate}
\end{definition}

Let $Y$, $Y'$ be as in Definition \ref{def:b-subm}, $\mathfrak{t}=(t_1,\dots, t_I)$ be a multi-weight associated with the boundary hypersurfaces of $Y'$, and let $f$  be an interior b-map.
Then we define a multi-weight for $Y$ by
\begin{equation}
	f^{\sharp}\mathfrak{t}=\mathfrak{s}=(s_1,\dots,s_J), \quad \text{where}\quad s_{j}=\sum_{i=1}^{I}e_f(i,j)t_i.
\end{equation}
Further, let  $\mathcal{F}$ be a collection of index sets for the boundary faces of $Y'$.
Then let 
\begin{equation}\label{eq:pullback_is}
 	f^\sharp(\mathcal{F})=\mathcal{E},\quad \text{ where }\quad E_j =\overline{\big \{(S,P):S=\sum\nolimits_{i=1}^Ie_f(i,j)s_i,\ P=\sum\nolimits_{e_f(i,j)\neq 0}p_i, \quad (s_i,p_i)\in F_i\big\}}.
\end{equation}
In \eqref{eq:pullback_is}, by convention $P=0$ if $e_f(i,j)=0$ for all $i$ (which means that the $j$-th boundary hypersurface is not mapped into any boundary hypersurface). 
So if $e_f(i,j)=0$ for all $i$, $E_j=\overline{\{0,0\}}$.

\begin{theoremofothers}[Pullback theorem, {\cite[Theorem 3]{Melrose1992}}]\label{thm:pullback}
	If $f:Y\to Y'$ is an interior b-map between manifolds with corners, pull-back gives maps
	\begin{align}
			f^{*}: \mathcal S^{\mathfrak{t}}(Y')\to  \mathcal S^{f^{\sharp}\mathfrak{t}}(Y),\quad   f^{*}:\mathcal{A}_-^{\mathfrak{t}}(Y')\to \mathcal{A}_-^{f^{\sharp}\mathfrak{t}}(Y), \quad 
			f^{*}:\Aphg^{\mathcal{F}}(Y')\to \Aphg^{f^\sharp(\mathcal{F})} (Y).
	\end{align}
\end{theoremofothers}

If $f$ is a b-fibration and $\mathcal{E}$ is a collection of index sets corresponding to the boundary faces of $Y$, we define the pushforward index set $f_\sharp (\mathcal E)$ as the family of index sets $\mathcal{F}=(F_1,\dots,F_I)$, where
	\begin{equation}
		F_i=\overline{\bigcup}_{j:e_f(i,j)\neq 0}\Big(\Big(\frac{s}{e_f(i,j)},p\Big):(s,p)\in E_j\Big),
	\end{equation}
	and the extended union of two index sets is defined as 
	\begin{equation}
		    E\overline{\cup}E'\coloneqq E\cup E'\cup\big\{(s,p+p'+1):\text{there exist }(s,p)\in E \text{ and } (s,p')\in E'\big\}.
	\end{equation}

\begin{theoremofothers}
	[Pushforward theorem, {\cite[Theorems 4 \& 5]{Melrose1992}}]\label{thm:pushforward}
	Let $f:Y\to Y'$ be a b-submersion between compact manifolds with corners.
	Suppose $\mathfrak{s}=(s_1,\dots, s_J)$ is a multiweight for $Y$ and $ \mathfrak{t}=(t_1,\dots,t_I)$ is a multiweight for $Y'$ such that 
	\begin{equation}
		s_j>0 \text{ for every }j \text{ such that  } e_f(i,j)=0\quad  \forall i=\{1,\dots ,I\}; \quad \text{and }\mathfrak{s}\geq f^{\sharp}\mathfrak{t}.
	\end{equation}
	Then the pushforward yields a map
	\begin{equation}
		f_{*}:\mathcal{A}_-^{\mathfrak{s}}(Y;\Omega_b)\to \mathcal{A}_-^{\mathfrak{t}}(Y';\Omega_b).
	\end{equation}
Moreover, let $\mathcal{E}=(E_1,\dots,E_I)$ be a family of index sets on $Y$ such that $\Re E_j>0$ for every $j$ satisfying $e_f(i,j)=0$ for all $i$ and suppose in addition that $f$ is a b-fibration.
Then the pushforward yields a map
	\begin{align}
			f_{*}:\Aphg^{\mathcal{E}}(Y;\Omega_b)\to \Aphg^{f_\sharp(\mathcal{E})} (Y';\Omega_b).
	\end{align}
\end{theoremofothers}

\subsection{Blow-ups}\label{sub:manifolds_with_corners_p_submanifolds_blowups}

A detailed exposition on blow-ups  can be found in  \cite{MelroseDAMC} and an introductory text is  \cite{Grieser2001}.
The process of blowing up a submanifold essentially amounts to introducing polar coordinates about it.
Formally, given a manifold with corners $Y^{m}$, a boundary 
 $p$-submanifold $S$ of codimension $l$ is a subset $S$ with the property that for every $p\in S$ there exists a corner chart $(x_1,\dots , x_m) \subset [0,\infty)^k\times \Rm^{m-k}$ for $Y$ defined on a neighborhood $\mathcal{U}$ of $p$, such that in terms of this chart, $S\cap \mathcal{U}$ is given as the zero set of $l$ of the coordinates, with at least one of them being a defining function for a boundary hypersurface.
If $S\subset Y$ is a closed p-submanifold, the \textbf{blow-up space} $[Y;S]$ is then $(Y\setminus S)\sqcup((N_+S\setminus 0)/\R^+)$, where $N_+S$ is the inward pointing (algebraic) normal bundle to $S$, 0 denotes its zero section, and the $\R^{+}$-quotient is with respect to the radial action on the fibers.
$[Y;S]$ has a natural smooth structure as a manifold with corners such that the blowdown map $\beta :[Y;S]\to Y$, which is equal to the identity on $Y\setminus S$ and projection to the base point on $(N_+S\setminus 0)/\R^+$, becomes smooth. 
 The new boundary face generated in this process (namely $(N_+S\setminus 0)/\R^+$) is called the \textbf{front face}, denoted $\mathrm{ff}([Y;S])$. 
The subspace of $\mathcal{V}_b(Y)$ consisting of vector fields which are additionally tangent to $S$ lift to span, over $C^\infty([Y;S])$, the space $\mathcal{V}_b([Y;S])$ (in particular,  vector fields tangent to $S$ lift to be tangent to $\mathrm{ff}([Y;S])$).
 The lift $\beta^{* }Z$ of a closed subset  $Z$ of $Y$ is defined as $\beta^{-1}(Z)$ if $Z\subset S$ and as $ \overline{\beta^{-1}(Z\setminus S)}$ if $Z=\overline{Z\setminus S}$, with the closures in the respective spaces.

\subsection{0-Pseudodifferential Operators
}
\label{ssec:0pdo}

We briefly introduce notations regarding 0-pseudodiffe- rential operators on a compact manifold with boundary $M$.
They were introduced in \cite{Mazzeo1987} and further developed in \cite{Mazzeo1991}, also see \cite{Hintz2021}.
 We refer the reader to these works or \cite{Eptaminitakis2022} for a more detailed exposition.
0-pseudodifferential operators are defined via the  behavior of their Schwartz kernel in the \textbf{0-stretched space} $M_0^2\coloneqq [M^2;\partial \mathrm{diag}]$, which is obtained by blowing up the boundary of the diagonal in the product $M^2=M\times M$.
$M_0^2$ is a manifold with corners up to codimension two with  three boundary hypersurfaces, called the left, right, and front face.
We denote by $\Psi_0^k(M)$ the operators  of order $k$ in the \textbf{small calculus}, acting on 0-half densities:
their Schwartz kernels are sections of the kernel half density bundle $\mathrm{KD}_0$ over $M_0^{2}$ (see \cite{Hintz2021}) which are conormal to the lift of the diagonal and vanish to infinite order at the left and right face.
The operators in the \textbf{large calculus} are denoted by $\Psi_0^{k,\mathcal{E}}(M)$, where $\mathcal E=(E_\ell,E_r,E_f)$;
their Schwartz kernels are sections of $\mathrm{KD}_0$ which are conormal of order $k $ to the lifted diagonal in $M_0^2$ and have polyhomogeneous expansions with index sets $E_\ell$, $E_r$, $E_f$ respectively at the left, right, and front face of $M_0^2$. 
We still use the shorthand notations for index sets and omit $E_f$ if it is $\overline{\{(0,0)\}}$, so the notation $\Psi_0^{k,a,b}(M)$ indicates that $E_l=\overline{\{(a,0)\}}$, $E_r=\overline{\{(b,0)\}}$, $E_f=\overline{\{(0,0)\}}$.
	
\newcommand{\etalchar}[1]{$^{#1}$}

\end{document}